\documentclass[english,12pt,leqno]{article}
\usepackage{tikz}
\usetikzlibrary{decorations.markings}
\usetikzlibrary{calc}

\usepackage[latin1]{inputenc}

\usepackage{amsxtra}
\usepackage{amsmath}
\usepackage{amssymb}
\usepackage{amsfonts}
\usepackage[all]{xy}
\usepackage{mathrsfs}
\usepackage{amsthm}
\usepackage{float}
 \usepackage{caption}

\usepackage{babel}

\newtheorem{thm}[equation]{Theorem}

\newtheorem{cor}[equation]{Corollary}
\newtheorem{lem}[equation]{Lemma}
\newtheorem{prop}[equation]{Proposition}

\newtheoremstyle{example}{\topsep}{\topsep}%
     {}
     {}
     {\bfseries}
     {.}
     {2pt}
     {\thmname{#1}\thmnumber{ #2}\thmnote{ #3}}

   \theoremstyle{example}
   
   \newtheorem{Defi}[equation]{Definition}
   
   \newtheorem{rem}[equation]{Remark}
   \newtheorem{rems}[equation]{Remarks}
   
   \newtheorem{exas}[equation]{Examples}
   \newtheorem{ex}[equation]{Example}

\newtheorem{exa}[equation]{Example}
\newtheorem{conj}[equation]{Conjecture}

\newtheoremstyle{example}{\topsep}{\topsep}%
     {}
     {}
     {\bfseries}
     {.}
     {2pt}
     {\thmname{#1}\thmnumber{ #2}\thmnote{ #3}}

   \numberwithin{equation}{section}

\def\VV{\mathbb{V}}

\def\CC{\mathbb{C}}

\def\PP{\mathbb{P}}
\def\RR{\mathbb{R}}
\def\ZZ{\mathbb{Z}}
\def\QQ{\mathbb{Q}}

\def\Ps{\mathscr{P}}
\def\scrP{\mathscr{P}}
\def\Qs{\mathscr{Q}}

\def\een{\mathfrak{e}}
\def\gen{\mathfrak{g}}

\def\hen{\mathfrak{h}}

\def\men{\mathfrak{m}}

\def\Pen{\mathfrak{P}}
\def\Sen{\mathfrak{S}}

\def\wen{\mathfrak{w}}

\def\Zen{\mathfrak {Z}}

\def\Ac{\mathcal{A}}

\def\Cc{\mathcal{C}}

\def\Dc{\mathcal{D}}

\def\Fc{\mathcal{F}}

\def\Lc{\mathcal{L}}
\def\Mc{\mathcal{M}}
\def\Nc{\mathcal{N}}

\def\Pc{\mathcal{P}}

\def\Rc{\mathcal{R}}
\def\Sc{\mathcal{S}}
\def\Tc{\mathcal{T}}

\def\Xc{\mathcal{X}}

\def\Sen{\mathfrak{S}}

\def\<{\langle}
\def\>{\rangle}

\def\Aff{\on{Aff}}
\def\Ass{{\Ac ss}}
\def\Aut{\on{Aut}}

\def\be{\begin{equation}}
\def\ee{\end{equation}}

\def\bef{\begin{figure}[H]\centering}
\def\enf{\end{figure}}

\def\btp{\begin{tikzpicture}}
\def\etp{\end{tikzpicture}}

\def\bgen{{\overset\bullet\gen}}
\def\bR{{\overset\bullet R}}

\def\bphi{{\boldsymbol\phi}}
\def\bpsi{{\boldsymbol\psi}}

\def\coh{{\on{coh}}}
\def\Coker{\on{Coker}}
\def\Conv{\on{Conv}}

\def\Der{\on{Der}}

\def\FS{\mathcal{F}\mathcal{S}}

\def\Hom{\on{Hom}}

\def\Imm{\on{Im}}
\def\Id{{\on{Id}}}

\def\k {\mathbf k}
\def\Ker{\on{Ker}}

\def\LG{\on{LG}}
\def\Lie{{\Lc ie}}

\def\Mod{\on{Mod}}

\def\lra{\longrightarrow}

\def\ol{\overline}
\def\on{\operatorname}
\def\op{{\on{op}}}
\def\orr{\on{or}}

\def\Ree{\on{Re}}
\def\RS{\mathcal{R}\mathcal{S}}

\def\vC{{{\overrightarrow C}}{}}
\def\vol{\on{ dVol}}

\def\wt{\widetilde}

\title{  Algebra of the infrared and secondary polytopes}

\author{ M. Kapranov, M. Kontsevich, Y. Soibelman}

  \newenvironment{dedication}
        {\vspace{3ex}\begin{quotation}\begin{center}\begin{em}}
        {\par\end{em}\end{center}\end{quotation}\vspace{5ex}}

\begin{document}


 \maketitle

\begin{dedication} To the memory of Andrei Zelevinsky
\end{dedication}


 \addtocounter{section}{-1}

\tableofcontents

 \vfill\eject

\section{Introduction}

The words ``algebra of the infrared" in the title refer to the physical paper \cite{GMW} by Gaiotto, Moore and Witten, to
which (or, rather, to a part of which) our article is a mathematical commentary. 

\vskip .2cm

In \cite{GMW} the authors developed an algebraic formalism for  the study  of 
certain $2$-dimensional massive  quantum field
theories with $(2,2)$ supersymmetry. One of the main algebraic structures introduced in the {\em loc.cit.}
 is the {\it $L_{\infty}$-algebra of webs},  associated with a generic finite subset  $A\subset {\bf R}^2$
 in the plane. Physically, elements of $A$ correspond to vacua of the theory. 
 A web is a plane graph with faces marked by elements of $A$ with an additional condition on the direction of edges,
 see \S \ref{sec:GMW} below for a review.
  Further, a choice of a half-plane containing $A$ determines an $A_{\infty}$-algebra (or an $A_\infty$-category, if one introduces
  a  coefficient system). This $A_\infty$-category has an ``upper-triangular structure", i.e. a semi-orthogonal decomposition.

  Using certain ``moduli spaces of $\zeta$-instantons" the authors of {\em loc.cit.}
   describe a class of deformations of the above $A_{\infty}$-category which describe the D-brane categories for
   a particular class of (2,2) supersymmetric theories:  Landau-Ginzburg models. Mathematically, the D-brane $A_\infty$-categories
  corresponding to LG models are known as Fukaya-Seidel categories \cite{HIV},\cite{FSbook}.

  \vskip .2cm

  We reinterpret and develop the algebraic structures proposed in \cite{GMW} in a way that allows a generalization to higher
  dimensions ($\RR^d, d\ge 2$ instead of just $\RR^2$). 
 It turns out that using the dual language of polygons rather than webs,  one  quickly uncovers certain structures well-known
  in toric geometry, most notably, {\em secondary polytopes},  see \cite{GKZ}. One of the  subtle and surprising
   points of the GMW construction is the fact that the differential they define, satisfies $d^2=0$.  In our approach this fact becomes
 obvious: the cellular chain complex of any  polytope, in particular, of the secondary polytope, is indeed a complex, i.e., it has $d^2=0$.

  Since secondary polytopes make sense in any number of dimensions, our dual approach leads naturally to higher-dimensional
  structures of ``extended" topological field theories (TFT). We postpone the study of these structures to a later work. So the
  present paper is just the first step in that general direction. 
  
Let us discuss the contents.
   
   \vskip .2cm
   
   Although we were strongly influenced by \cite{GMW}, we have arranged the exposition logically from scratch,
    implementing our dual interpretation right away. A detailed dictionary relating our approach and terminology with that of \cite{GMW}, is
   provided in Section \ref{sec:GMW}. The reader familiar with the terminology of \cite{GMW} can start from that section and then proceed linearly from the beginning.

Section \ref{sec:rem-sec} provides the general background on secondary polytopes. Roughly, the secondary 
polytope of a given polytope in ${\bf R}^d$ 
has vertices corresponding to decompositions of the initial, ``primary",  polytope into simplices (such decompositions
 are called triangulations). The usefulness of the notion of secondary polytope from the point of view of the approach 
 of \cite{GMW} becomes clear in Section \ref{sec:com-dg-alg}: the chain complex of the secondary polytope looks like a dg-algebra. 
 The main reason for that is explained earlier in Proposition \ref{prop:fact-sigma}: each face of the secondary polytope is itself 
  a product of several secondary polytopes associated to members of a certain polyhedral subdivision of the initial polytope.

In Section \ref{sec:L-inf} we introduce an $L_{\infty}$-algebra structure in the space $\gen_A$ spanned by polytopes with vertices in 
  a given set $A\subset {\bf R}^d$. For that we  define a differential in the symmetric algebra generated  by the (shifted) dual space
  to $\gen_A$. The  remark in the previous paragraph about the chain complex of the secondary polytope explains why our $L_{\infty}$-structure 
  comes naturally from such a chain complex. Further, one can visualize higher Lie brackets in $\gen_A$ geometrically. 
   Roughly speaking, they correspond to the operation of composing a
    convex polytope from smaller polytopes with vertices in $A$.

Section \ref{sec:MC}  is a reminder about Maurer-Cartan elements in $L_{\infty}$-algebras. Together with Section \ref{sec:hoch}, 
which is devoted to a short account on deformation theory, Hochschild complexes etc., it prepares the ground for future applications \`a la \cite{GMW}:  a geometrically defined Maurer-Cartan element gives rise to a deformed $L_{\infty}$-structure.

Section \ref{sec:1dim} is devoted to the most elementary case $d=1$, which explains the reason for the second construction of 
  \cite{GMW}: lifting of a certain  $L_\infty$-algebra (related but not equal to $\gen_A$) to an $A_\infty$-algebra. 
  The possibility of such lifting is the essential 1-dimensionality
  of the situation, and the purely 1-dimensional case provides a clear starting point. 
   This is related to the fact that convex polytopes on the real line are just segments, and  so one can naturally order their vertices. 

In Section \ref{sec:rel-set} we study the situation when one element of $A$ is distinguished and referred to as an ``element at infinity"
(denoted $\infty$). This is a higher-dimensional analog of the choice of a half-plane in \cite{GMW}. In this case
the $L_\infty$-algebra $\gen_A$ contains an ideal $\gen_\infty$ spanned by ``infinite" polytopes (containing $\infty$),
and a subalgebra $\gen_{\on{fin}}$ spanned by ``finite" polytopes (not containing $\infty$). 
We explain in Section \ref{sec:refi} how these structures lead, in the case $d=2$, to an action of $\gen_{\on{fin}}$
on $R_\infty$,  a natural $A_\infty$-lift of the $L_\infty$-algebra $\gen_\infty$. Combining this with the results of
Sections \ref{sec:rel-set}, \ref{sec:hoch}, 
 we arrive at a morphism from $\gen_{\on{fin}}$ to the Hochschild complex of $R_\infty$ (all in the
case of trivial coefficients).

The case of non-trivial coefficients is the subject of Section \ref{sec:coeff-arbitr}. Here we introduce and exploit an important notion of 
a {\em factorizing sheaf} (on the secondary polytope). This notion can be seen as a combinatorial version of factorization
algebras such as used by Lurie and Costello-Gwylliam \cite{lurie-fac, costfac}  to describe higher-dimensional TFT. This strongly suggests that the
constructions of the previous sections admit a generalization to the framework of higher-dimensional TFT, with, for instance,
   $E_{d-1}$-algebras
playing the role of $A_\infty$-algebras such as $R_\infty$.  

 For an $A_\infty$-algebra with an upper-triangular structure, there is a version of the Hochschild complex which 
 controls deformations preserving this structure. We call it the {\em directed Hochschild complex} and 
   discuss it  in Section \label{sec:univ}.
After preparations in Sections \ref {sec:bim}, \ref {sec:anal} we finally prove,  in Section \ref{sec:univ}
a result which we call the Universality Theorem. It concerns the case $d=2$ and the situation of Section \ref{sec:anal}
associated with a choice of a point $\infty$ (or, equivalently, of a held-plane). 
 It says that 
   the natural morphism from the $L_{\infty}$-algebra $\gen_{\on{fin}}$ of finite   polygons to the Hochschild complex
   of $R_\infty$, is a quasi-isomorphism onto the directed Hochschild complex.

This result should be compared with the Deligne conjecture, see \cite{defoperads},  as well as 
with similar results about the Swiss Cheese operad and its higher-dimensional generalizations, see e.g.,  \cite{operads}. 
It certainly admits a generalization to the higher-dimensional case (in which  case one deals with Hochschild complexes of $E_{d-1}$-algebras 
in the spirit of  \cite{operads}). This will be a subject of separate publication.

In Section \ref {sec:MCFS}   we discuss possibilities for the Maurer-Cartan element giving the Fukaya-Seidel category associated with a choice of a half-plane.

Section \ref {sec:spec} contains a list of further applications and speculations.

We should mention here that there are several $L_{\infty}$ and $A_{\infty}$ algebras 
that can be naturally associated with a finite subset $A\subset {\mathbb R}^d$. 
The biggest one involves all marked polytopes. It has non-trivial differential corresponding to the operation 
of ``insertion of  internal points". 
 On the other hand, for the application to Fukaya-Seidel categories one needs a much smaller $L_{\infty}$-subalgebra which involves
  polytopes, each  marked with  the set of all its internal points. 
  The differential is trivial on this subalgebra, and moreover, we do not need to include the set of points as a separate piece of the data. 
    In this way we obtain a smaller $L_{\infty}$-algebra, which we refer to as the {\em geometric} one.
    It is  spanned by convex polytopes with vertices in $A$, with higher Lie brackets given by
     the operation of composing a bigger convex polytope from smaller ones. 
As  we pointed out above, the geometric $L_{\infty}$-algebra appears naturally in the description of the Hochschild complex of the Fukaya-Seidel category. The meaning of the ``big" $L_{\infty}$-algebra is not clear at the moment.

{\it Acknowledgements.} We thank  Davide Gaiotto, Greg Moore and Edward Witten for multiple discussions 
and correspondences about their deep and beautiful paper \cite{GMW} as well as for sending us
preliminary drafts. Y.S. is grateful to IHES for excellent research conditions. 
His work was partially supported by an NSF grant.

\vfill\eject

\section {Reminder on secondary polytopes.}\label{sec:rem-sec}

In this section we recall some basic material from \cite{GKZ}, Ch. 7.  We refer the reader to 
 {\em loc.cit.}  for  details of constructions and proofs of statements.

\vskip .2cm

Let $A\subset\RR^d$ be a finite subset of points. We assume that:

\begin{enumerate}
\item[(1)] The affine span of $A$ is the whole of $\RR^d$. 

\item[(2)] $A$ is in general position, i.e., any $p\leq d+1$ points of $A$ are the vertices of a $(p-1)$-simplex. 

\end{enumerate} 
The first assumption can always be satisfied by passing to the affine span. The second assumption can be removed
at the price of somewhat complicating the discussion and we prefer to keep it throughout the paper for simplicity. 

\vskip .2cm

Let $Q=\Conv(A)$ be the convex hull of $A$, a $d$-dimensional convex polytope. 
By a {\em triangulation} $\Tc$  of $(Q,A)$ we will mean a subdivision of $Q$ into straight simplices of full dimension $d$
with vertices in $A$ so that the intersection of any two simplices is a common face (possibly empty).
Note that it is not required that each element of $A$ appears as a vertex of a simplex of $\Tc$. 
A triangulation  $\Tc$ is called {\em regular}, if there is a continuous convex function $f: Q\to\RR$ such that:
\begin{itemize}
\item $f$ is affine-linear on each simplex of $\Tc$.

\item $f$ is not affine-linear on any subset of $Q$ which is not contained in a simplex of $\Tc$. In other words,
$f$ does indeed  break along each codimension 1 simplex which is a common face of two different $d$-dimensional simplices
of $\Tc$. 
\end{itemize}

 For a more systematic discussion, we need the following concept. 

\begin{Defi} A 
 {\em marked polytope}   is a pair
  $(Q,A)$
where $Q$ is a convex polytope in $\RR^d$ and $A\subset\RR^d$ is a finite subset such that $Q=\Conv(A)$,
i.e., $A$ contains all vertices of $Q$. A {\em marked subpolytope} of $(Q,A)$ is a marked polytope
$(Q', A')$ such that $A'\subset A$. Notation: $(Q', A')\subset(Q,A)$. 

\end{Defi}

\begin{Defi} By a {\em polyhedral subdivision} of $(Q,A)$ we will mean a collection $\Pc = (Q_\nu, A_\nu)$
of marked subpolytopes of $(Q,A)$  which have $\dim\,Q_\nu=\dim\,Q$ such that:
\begin{enumerate}
\item[(1)] The $Q_\nu$ form a polyhedral subdivision of $Q$ (i.e. $Q=\cup_{\nu}Q_\nu$) so that each $Q_\nu\cap Q_{\nu'}$ is
a common face of $Q_\nu$ and $Q_{\nu'}$ (possibly empty).
\item[(2)] In addition, we have $A_\nu\cap(Q_\nu\cap Q_{\nu'}) = A_{\nu'}\cap(Q_\nu\cap Q_{\nu'})$.
\end{enumerate}
As for triangulations, a
polyhedral  subdivision $\Pc=(Q_\nu, A_\nu)$ is called {\em regular}, if there is a continuous
convex function $f: Q\to\RR$ which is affine-linear on each $Q_\nu$ and is not affine linear
on any subset which is not contained in some $Q_\nu$. 

\end{Defi}

\bef
\centering
\btp[scale=.4, baseline=(current  bounding  box.center)]
\node (1) at (0,0){}; 
\fill (1) circle (0.15);

\node (2) at (4,-1){}; 
\fill (2) circle (0.15);

\node (3) at (-1,-3){};
\fill (3) circle (0.15);

\node (4) at (-5,-1){};
\fill (4) circle (0.15);

\node (5) at (-2,2){};
\fill (5) circle (0.15);

\node (6) at (3,3){};
\fill (6) circle (0.15);

\node (7) at (0,1){};
\fill (7) circle (0.15);

\node (8) at (2,1){};
\fill (8) circle (0.15);

\node (9) at (3,0){}; 
\fill (9) circle (0.15);

\draw (0,0) -- (4,-1);
\draw (4,-1) -- (-1,-3);
\draw (0,0) -- (-1,-3); 

\draw (0,0) -- (3,3);
\draw (4,-1) --(3,3);
\draw (-1,-3) --(-5, -1);
\draw (-5,-1) -- (-2,2); 
\draw (-2,2) -- (0,0); 
\draw (-2,2) -- (3,3); 
\draw (0,0) -- (-5,-1);

\etp
\hskip 2cm
\btp[scale=.4, baseline=(current  bounding  box.center)]

\node (10) at (0,0){};
\node (11) at (0,2){};
\fill (11) circle (0.15);
\node (12) at (0,-2){};
\fill (12) circle (0.15);

\draw (0,2) -- (0, -2);
\draw (0,2) -- (-3,3);
\draw (0,-2) -- (-2,-4);
\draw (0,2) -- (2,4);
\draw (0,-2) -- (2,-3);
\node at (-2,0) {$Q_\nu$};
\node at (2,0){$Q_{\nu'}$};

\etp
\caption{A  triangulation and a polyhedral subdivision. }
\enf

Regular polyhedral subdivisions of $(Q,A)$ form a poset (partially ordered set)
$\Rc(Q,A)$ with the order given by {\em refinement}. That is,
\[
\Pc' = (Q'_\mu, A'_\mu)\,\,\,  \leq \,\,\,  \Pc=(Q_\nu, A_\nu),
\]
iff $\Pc'$ induces a (necessarily regular) polyhedral subdivision of each $(Q_\nu, A_\nu)$. 
Thus, minimal elements of $\Rc(Q,A)$ are precisely the regular triangulations, 
while the unique maximal element is the subdivision consisting of $(Q,A)$ alone.

\begin{rem}
Note the essential   role of the choices of markings $A_\nu$
in the definition of polyhedral subdivisions and of the poset $\Rc(Q,A)$. We do not require that $A=\cup_\nu A_\nu$.
E.g. if $(Q,A)$ is a triangle with two points in
 the interior then
$(Q',A')$ is a polyhedral subdivision provided $Q'=Q$ and $A'$ is just one of the interior points.

More generally, if
$A'\subset A$ is the set of vertices of $Q$, then each intermediate subset $A'\subset B\subset A$
gives rise to a 1-element regular polyhedral subdivision $(Q,B)\in\Rc(Q,A)$, and
\[
(Q,B') \, \leq \, (Q, B)\,\,  \text{ iff } \,\, B'\subset B.
\]

\end{rem}

The poset $\Rc(Q,A)$ has an interpretation in terms of the
  {\em secondary fan} of $A$. This fan,  denoted $\Sc(A)$,
is a subdivision of
$\RR^A$ into convex polyhedral cones $C_\Pc$ corresponding to regular polyhedral
subdivisions $\Pc\in\Rc (Q,A)$. To define it, we associate to each $\psi: A\to\RR$  
a convex piecewise-affine function $f_\psi: Q\to\RR$ as follows. We consider the
unbounded polyhedron $G_\psi\subset\RR^{d+1} = \RR\times\RR^d$ given by
\[
G_\psi \,\,=\,\, \Conv \bigl\{ (t, a) \bigl| \, t\in\RR, a\in A\subset\RR^d, \,\, t\geq\psi(a)\bigr\}. 
\]

\begin{figure}
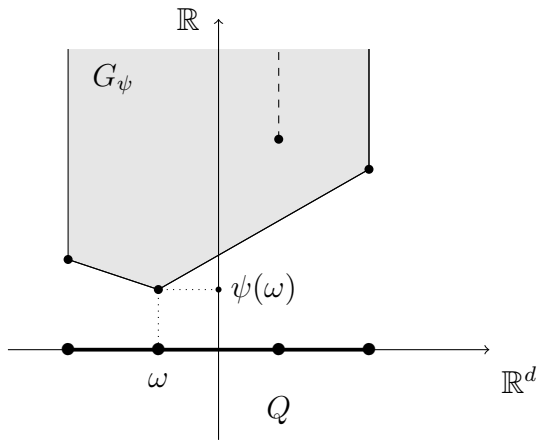

\centering
 \btp[scale=.4, baseline=(current  bounding  box.center)]
 \draw [->] (-7,0) -- (9,0); 
 \node at (10,-1){$\RR^d$}; 
 
 \draw [->] (0,-3) -- (0,11); 
 \node at (-1,11) {$\RR$}; 
 
 \node at (-5,0){$\bullet$};
  \node at (5,0){$\bullet$};
   \node at (-2,0){$\bullet$};
    \node at (2,0){$\bullet$};
    
\node (1) at (-5,3){};
\fill (1) circle (0.15);    

\node (2) at (-2,2){};
\fill (2) circle (0.15);    

\node (3) at (2,7){};
\fill (3) circle (0.15);    

\node (4) at (5,6){};
\fill (4) circle (0.15);    

\draw[line width=.5mm] (-5,0) -- (5,0); 
\node at (2,-2){$Q$}; 
 
 \filldraw[fill opacity=0.1] (-5,10) -- (-5,3) -- (-2,2) --(5,6) -- (5,10); 
 \draw (-5,3) -- (-2,2); 
 \draw (-2,2) -- (5,6); 
 \draw (5,6) -- (5,10); 
 \draw[dashed] (2,7) -- (2,10); 
 
 \node at (-3.5,9){$G_\psi$};
 
 \node at (-2,-1){$\omega$}; 
 \draw[dotted] (-2,0) -- (-2,2); 
  \draw[dotted] (0,2) -- (-2,2); 
  \node (5) at (0,2){};
  \fill(5) circle (0.1);
  \node at (1.5,2){$\psi(\omega)$};

 \etp
\caption{The polyhedron $G_\psi$ and the function $f_\psi$.}
\end{figure}

This image of $G_\psi$ under the projection to $\RR^d$ is $Q$, and the bottom of $G_\psi$
is the graph of a convex piecewise-affine function $f_\psi: Q\to\RR$. The function $f_\psi$
gives rise to unique (necessarily regular) polyhedral subdivision 
  $\Pc(\psi)$  consisting of {\em domains of affine linearity} of $f_\psi$. 
  The cone $C_\Pc\subset\RR^A$ is, by definition, the set of all $\psi$ such that $\Pc_\psi=\Pc$. 
  
  Note that each $C_\Pc$ is invariant under translations by those $\psi$ which come from
global affine linear functions on $\RR^d$.
We denote the space of such functions by $\Aff(\RR^d)$. 
 Then $C_\Pc$ gives rise to the cone  $\overline C_\Pc$ in the quotient space
$\RR^A/\on{Aff}(\RR^d)$. The collection of these cones is
 called the {\em reduced secondary fan} of $A$ and denoted by
$\overline\Sc(A)$.

  \begin{prop}
  The correspondence $\Pc\mapsto C_\Pc$ is an order reversing isomorphism between the poset $\Rc(Q,A)$
  and the poset of cones of $\Sc(A)$  (ordered by inclusion). \qed
  
  \end{prop}
  
 Along with $\Sc(A)$ we will consider the dual object, the {\em secondary polytope}
 $\Sigma(A)$ whose faces correspond to subdivisions $\Pc\in\Rc(Q,A)$
 in an order-preserving, not reversing way. The construction of $\Sigma(A)$ depends on
 the choice of a translation invariant measure $\on{Vol}$ on $\RR^d$, a different choice
 leading to a rescaling of $\Sigma(A)$. For convenience of the reader we recall
 two (equivalent) definitions.
 
 \begin{Defi}
 For a regular triangulation $\Tc$ of $(Q,A)$ define the vector $\phi_\Tc\in\RR^A$ by
 \[
 \phi_\Tc(\omega) \,\,=\,\,\sum_{\stackrel{\Delta\in\Tc}{\omega\in\on{Vert}(\Delta)} 
 }
 \,\on{Vol}(\Delta). 
 \]
 The polytope $\Sigma(A)\subset\RR^A$ is defined to be the convex hull of the vectors
 $\phi_\Tc$ for all regular triangulations $\Tc$ of $(Q,A)$. 
 \end{Defi}
 
 \begin{Defi}\cite{BS}
 Consider the standard simplex with the set of vertices $A$:
 \[
 \Delta^A \,\,=\,\,\biggl\{ (p_\omega)_{\omega\in A} \in\RR^A \bigl| \,\,\,
 p_\omega\geq 0, \,\,\sum_\omega p_\omega=1 \biggr\}, 
 \]
 so that, by definition of the convex hull, we have a surjective projection
 \[
 \pi: \Delta^A \lra Q=\Conv(A), \,\,\, (p_\omega)_{\omega\in A} \,\mapsto \, \sum_\omega p_\omega\cdot\omega.
 \]
 The polytope $\Sigma(A)$ is defined as the set of vector integrals
 \[
 \int_{q\in Q} s(q) \vol \,\,\in \,\,\RR^A
 \]
 for all continuous sections $s:Q\to\Delta^A$ of the projection $\pi$. 
 \end{Defi}

  To explain the relation of the two definitions, 
  note that a triangulation $\Tc$ of $(Q,A)$ defines a ``tautological" piecewise linear section
  $s_\Tc: Q\to\Delta^A$, and $\phi_\Tc = \int_Q s_\Tc(q) \vol$.

  \vskip .2cm

 To each face $F=F_\Pc$ of $\Sigma(A)\subset\RR^A$ we can associate its {\em normal cone}
$N_F\subset (\RR^A)^*$ consisting of those linear functionals $l: \RR^A\to\RR$ which achieve
the maximum on $F$ and are constant on $F$. The cones $N_F$ form  a decomposition
(fan) of the space $(\RR^A)^*$ called the {\em normal fan} of $\Sigma(A)$. 
Let us  now identify $\RR^A$ with its dual vector space by means of the standard pairing
\[
(\phi, \psi) \,\,=\,\,\sum_{\omega\in A} \phi(\omega) \psi(\omega). 
\]

Under this identification we have:

\begin{prop}
   The normal fan of 
 $\Sigma(A)$ is identified with the secondary fan $\Sc(A)$. In particular:
 
 \vskip .2cm
 
 (a) The face lattice of $\Sigma(A)$ is isomorphic to $\Rc(Q,A)$ in an order-preserving way. 
 We denote by $F_\Pc$ the face corresponding to a  subdivision $\Pc\in\Rc(Q,A)$. 
 
  \vskip .2cm
 
 (b) The codimension of $F_\Pc$ is equal to the dimension of the cone $\overline C_\Pc$. 
 
  \vskip .2cm
 
   (c) $\dim \Sigma(A) = |A|-d-1 = \dim(\RR^A/\Aff(\RR^d))$.

\end{prop}

\begin{exas}\label{ex-faces-sigma}
(a) A triangulation of $(Q,A)$ can be seen as a polyhedral subdivision $\Tc=(Q_\nu, A_\nu)$ such that
each $Q_\nu$ is a simplex and $A_\nu$ consists exactly of the vertices of $Q_\nu$. In this case
$\Sigma(A_\nu)$ is a point and so the face $F_\Tc$ is a vertex.

\vskip .2cm

(b) Edges of $\Sigma(A)$ correspond to {\em flips} (elementary modifications) of triangulations.
A flip is based on a {\em circuit}, a subset $Z\subset A$ which has precisely one, up to a scalar, affine dependency
\[
\sum_{\omega\in Z} a_\omega \cdot \omega =0, \quad a_\omega\in\RR, \,\,\,\sum a_\omega = 0. 
\]
In this case the convex hull $\Conv(Z)$ has precisely two triangulations $T_+, T_-$ with vertices in $Z$, so
$\Sigma(Z)$ is an interval.  Each 
 edge
$[\phi_\Tc, \phi_{\Tc'}]$ of $\Sigma(A)$  corresponds to a pair $\Tc, \Tc'$ of triangulations of $Q$ which coincide outside
$\Conv(Z)$ for some circuit $Z$,  and restrict to $T_+, T_-$ inside $\Conv(Z)$.

\begin{figure}[H]
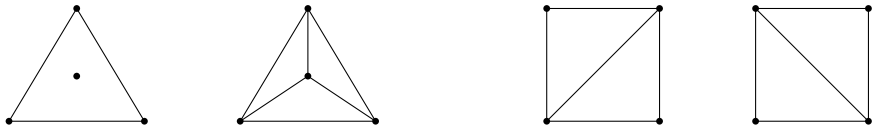

\centering
 \btp[scale=.3, baseline=(current  bounding  box.center)]
 \node (1) at (-3,0){};
 \fill (1) circle (0.15);
 
 \node (2) at (3,0){};
 \fill (2) circle (0.15);
 
 \node (3) at (0,5){};
 \fill (3) circle (0.15);
 
 \node (4) at (0,2){};
 \fill (4) circle (0.15);
 
 \draw (-3,0) --(3,0) -- (0,5) -- (-3,0); 
 \etp
 \hskip 1cm
 \btp[scale=.3, baseline=(current  bounding  box.center)]
 \node (1) at (-3,0){};
 \fill (1) circle (0.15);
 
 \node (2) at (3,0){};
 \fill (2) circle (0.15);
 
 \node (3) at (0,5){};
 \fill (3) circle (0.15);
 
 \node (4) at (0,2){};
 \fill (4) circle (0.15);
 
 \draw (-3,0) --(3,0) -- (0,5) -- (-3,0); 
 \draw (-3,0) -- (0,2) -- (0,5);
 \draw (0,2) -- (3,0); 
 \etp
 \hskip 2cm
 \btp[scale=.3, baseline=(current  bounding  box.center)]
\node (1) at (0,0){};
\fill (1) circle (0.15); 

\node (2) at (5,0){};
\fill (2) circle (0.15);

\node (3) at (5,5){};
\fill (3) circle (0.15);

\node (4) at (0,5){};
\fill (4) circle (0.15);

\draw (0,0) -- (0,5) -- (5,5) -- (5,0) -- (0,0);
\draw (0,0) -- (5,5);

\etp
\hskip 1cm
 \btp[scale=.3, baseline=(current  bounding  box.center)]
\node (1) at (0,0){};
\fill (1) circle (0.15); 

\node (2) at (5,0){};
\fill (2) circle (0.15);

\node (3) at (5,5){};
\fill (3) circle (0.15);

\node (4) at (0,5){};
\fill (4) circle (0.15);

\draw (0,0) -- (0,5) -- (5,5) -- (5,0) -- (0,0);
\draw (0,5) -- (5,0);

\etp
\caption{Circuits and flips}
\end{figure}

\vskip .2cm

 (c) Codimension 1 faces correspond to {\em coarse subdivisions} (see \cite{GKZ}). In the dual language used in
 \cite{GMW} in case $d=2$ they are called {\em taut webs}.
 By definition, a subdivision $\Pc$ is coarse, if the reduced normal cone $\overline C_\Pc$ has
 dimension 1, i.e., there exists only one,  modulo scaling and adding affine functions,
 convex $\Pc$-piecewise-affine function $f:  Q\to\RR$. The following are the most important examples
 of coarse subdivisions.
 
 \vskip .1cm
 
 (c1) Let $\omega\in A$ be not a vertex of $Q$. Then $(Q, A-\{\omega\})$ is a coarse subdivision of $(Q,A)$.
 
 \vskip .1cm
 
 (c2) Let $H$ is a hyperplane in $\RR^d$ subdividing $Q$ into two subpolytopes $Q_1, Q_2$.
 Suppose that  the vertices of both these polytopes lie in $A$. Put $A_\nu = A\cap Q_\nu$. Then
  $(Q_1, A_1)$ and $(Q_2, A_2)$
 form a coarse subdivision of $(Q,A)$.

 \begin{figure}[H]
 \centering
 \btp[scale=.4, baseline=(current  bounding  box.center)]

\node (1) at (0,0){};
\fill (1) circle (0.15);

\node (2) at (2,4){};
\fill (2) circle (0.15);

\node (3) at (6,2){};
\fill (3) circle (0.15);

\node (4) at (7,-1){};
\fill (4) circle (0.15);

\node (5) at (3,-3){};
\fill (5) circle (0.15);

\draw (0,0) -- (2,4) -- (6,2) -- (7,-1) -- (3,-3) -- (0,0); 

\draw (-3,-1) -- (12,4);
\node at (11,5){$H$};
\node at (2.5,2){$Q_1$};
\node at (3,-1){$Q_2$};

\etp
 \caption {A 2-part coarse subdivision. }
 \enf
\end{exas}

The importance of $\Sigma(A)$ for our purposes comes from the following.

\begin{prop}[{\bf Factorization property}]\label{prop:fact-sigma}
 The face   $F_\Pc$ corresponding to a regular polyhedral subdivision $\Pc = (Q_\nu, A_\nu)$, has the form
$F_\Pc = \prod_\nu \Sigma(A_\nu)$, i.e., it is itself a product of several secondary polytopes. 
\end{prop}

\begin{rem} Note that if we drop the assumption of $A$ being in general position, then 
  $F_\Pc$ will be 
not the full product but the
 fiber product of the $\Sigma(A_\nu)$ over the secondary polytopes of intermediate faces.
This one of the complications we wanted to avoid in this paper. 
\end{rem}

For future convenience, we introduce the following terminology.

\begin{Defi}\label{def:geom}
A marked subpolytope $(Q', A')\subset (Q,A)$ is called {\em geometric}, if $A'=A\cap Q'$. A polyhedral subdivision
$\Pc=\{(Q'_\nu, A'_\nu)\}$ of $(Q,A)$ is called geometric, if each $(Q'_\nu, A'_\nu)$ is a geometric marked subpolytope.
A face $F$ of $\Sigma(A)$ is called geometric, if it corresponds to a geometric subdivision. 
\end{Defi}

We will use the term {\em subpolytope} of $(Q,A)$ for a polytope $Q'\subset Q$ with vertices from $A$.
Notation: $Q'\subset (Q,A)$. Such subpolytopes are in bijection with geometric marked subpolytopes
$(Q', A\cap Q')$. 

The following property will be important for the construction of \S \ref{sec:L-inf} to make sense. Its proof is obvious
and is left to the reader.

\begin{prop}\label{prop:geom}
Geometric faces form a right ideal in the partially ordered set of all faces of $\Sigma(A)$. That is, if $F\subset F'$
and $F$ is geometric, then $F'$ is geometric. \qed
\end{prop}

We introduce the {\em geometric part} of the secondary polytope to be the union of the interiors of
geometric faces (of all dimensions): 
\be
\Sigma^{\on{geom}}(A) \,\,=\,\,\bigcup_{
\stackrel
{F\subset \Sigma(A)}  { \text{ geometric}} 
} 
 \on{Int}(F).
 \ee
 Proposition \ref{prop:geom} can be reformulated by saying that $\Sigma^{\on{geom}}(A)$ is an 
 open subset of $\Sigma(A)$. 

\vfill\eject

\section{A commutative dg-algebra from the chain complex of $\Sigma(A)$. }
\label{sec:com-dg-alg}

The factorization property of $\Sigma(A)$ implies that its chain complex looks like a
(part of) a multiplicative differential in  a dg-algebra. 

More precisely, let $\k$ be a field of characteristic 0. For any convex polytope $P$  we denote by $C_\bullet(P)$
its  cellular chain complex with coefficients in $\k$ graded so that $i$-chains are positioned in degree $(-i)$
and the differential raises the degree by $+1$. Thus, as a graded vector space
\[
C_\bullet (P) \,\,=\,\, \bigoplus_{F\subset P}  \orr(F) [\dim(F)],
\]
the direct sum over all faces. Here
\[
\orr(F) \,\,=\,\, H_c^{\dim(F)} (\text{Interior of } F, \k) 
\]
is the 1-dimensional orientation space of $F$. 

We apply this to $P=\Sigma(A)$ and form the graded vector space
\be\label{eq;space-V}
V \,\,=\,\, \bigoplus_{ (Q', A')\subset (Q,A)} V_{A'}, \quad V_{A'} =   \orr(\Sigma(A')) [\dim\Sigma(A')],
\ee
the direct sum over all marked subpolytopes of full dimension $d$. 
Here, the shift means that $\orr(\Sigma(A'))$ is positioned in degree $-\dim\Sigma(A')$. Note that if $\k=\RR$, then
\[
\orr(\Sigma(A')) \,\,=\,\,\Lambda^{\max} \bigl( \RR^{A'}/\on{Aff}(\RR^{d}).
\]
Now, take the symmetric algebra $S^\bullet(V)$ with the graded commutative product
denoted $\odot$. Each regular polyhedral subdivision 
$\Pc''=\{(Q''_\nu, A''_\nu)\}$ of each marked subpolytope $(Q', A')$ of $(Q,A)$  produces a 1-dimensional
subspace
\be\label{eq:V-P}
V_{\Pc''} \,\,=\,\,\bigodot_\nu \,\, V_{A''_\nu}   \,\,\subset \,\, S^\bullet(V). 
\ee
The chain differentials in all the $C_\bullet(\Sigma(A'))$ are compatible with each other and give a differential
$d$ of degree $+1$
\be
d: S^\bullet(V)\lra S^\bullet(V),  
\ee
satisfying the Leibniz rule. 

More precisely, let $\Pc''=\{(\Pc''_\nu, Q''_\nu)\}$ be a coarse polyhedral subdivision of $(Q', A')$, $A'\subset A$.
Using the factorization property  (Proposition \ref{prop:fact-sigma}), we see that
  the top degree part of the chain differential in $C_\bullet(\Sigma(A'))$ gives a map
\[
d_{\Pc''}: V_{A'}  \lra V_{\Pc''} \,\,\subset S^\bullet(V). 
\]
We define $d$ to be given, on generators, by
\[
d|_{  V_{A'} } \,\,=\,\,\sum_{ \substack {\Pc'' \text{ coarse} \\ \text{subdiv. of } (Q', A')}}
d_{\Pc''},
\]
and then extend it to the whole of $S^\bullet(V)$ by the Leibniz rule. 

\begin{prop}\label{prop:d2=0}

The differential $d$ thus defined, satisfies  $d^2=0$ and so makes $S^\bullet(V)$ into a commutative
dg-algebra. 
\end{prop}

 \noindent {\sl Proof:} In virtue of the Leibniz rule, it is enough to check that $d^2=0$ on generators. To see this, note that  on 
 $V_{A'}$ the picture for $d^2$ is embedded,  again by  
 Proposition \ref{prop:fact-sigma}, 
into the picture for $d^2$ in
 $C_\bullet(\Sigma(A'))$. 
  But for any polytope $P$,  the differential in the chain complex $C_\bullet(P)$ satisfies
$d^2=0$. \qed

\begin{rem}
We take the {\em symmetric} algebra of $V$ because, for a subdivision $\Pc''=\{(Q''_\nu, A''_\nu)\}$ of
$(Q', A')$, different marked polytopes $(Q''_\nu, A''_\nu)$ have no apparent order. But this can be refined
depending on  the dimension $d$, as we will see later. 
\end{rem}

\vfill\eject

\section{The $L_\infty$-algebra.}
\label{sec:L-inf}

If $\hen$ is a finite-dimensional dg-Lie algebra over $\k$, then its Chevalley-Eilenberg cochain complex is
a commutative dg-algebra
\[
C^\bullet_{\on{Lie}}(\hen) \,\,=\,\, S^\bullet(\hen^*[-1]). 
\]

For any graded vector space $V$ we denote by $S_+^{\bullet}(V)\subset S^{\bullet}(V)$ the maximal ideal $\oplus_{n\ge 1}V^{\otimes n}$.

As well known, a general algebra differential $d$  on $S^\bullet_+(\hen^*[-1])$ (i.e. an algebra derivation satisfying $d^2=0$), makes
$\hen$ into an {\em $L_\infty$-algebra}, see, e.g.,  \cite{defbook}, Ch.3, Section 2.3. 
Equivalently we can say that  {\em $L_\infty$-algebra} on $\hen$ is given by the algebra differential of $S^\bullet(\hen^*[-1])$ which preserves $S^\bullet_+(\hen^*[-1])$.

\vskip .2cm

Thus, in the situation of \S \ref{sec:com-dg-alg},
we get an $L_\infty$-structure on the vector space
\[
\begin{gathered}
\bgen = \bgen_A = V^*[-1] \,\,=\,\,
\bigoplus_{ (Q', A')\subset (Q,A)}  E_{A'}, \\
E_{A'} = V_{A'}^*[-1] \,\,=\,\, \orr(\Sigma(A'))^* [-1-\dim\Sigma(A')].
\end{gathered}
\]
In other words, the 1-dimensional $\k$-vector space $\orr(\Sigma(A))^*$
(which we can canonically identify with its dual $\orr(\Sigma(A))$), is
now put in the positive degree $1+\dim\Sigma(A')$. 

\vskip .2cm

The $L_\infty$-algebra $\bgen$ is, while natural, too big for our purposes. Indeed, a subpolytope
$Q'\subset (Q,A)$ with vertices in $A$ can give rise to many summands $E_{A'}\subset\bgen$
corresponding to the choice of $A'$ sandwiched between $\on{Vert}(Q')$ and $A\cap Q'$. 
The dot in the notation $\bgen$ is supposed to symbolize this freedom in choosing the set of points $A'$. 
Taking only the ``geometric
 summand corresponding to $A'=A\cap Q'$, we define the subspace
 \[
 \gen=\gen_A=\bigoplus_{(Q', A') \text{ geom. }} E_{A'} \,\,=\,\, \bigoplus_{Q'\subset (Q,A)} E_{A\cap Q;} \,\,\subset \,\,\bgen. 
 \]
 
 \begin{rem} We should warn the reader that in the Sections 3-8 we discuss $L_\infty$-algebras which are 
  toy models of the ``realistic" $L_\infty$-algebras appearing in the ``nature" (in particular, in 
 the study of Fukaya-Seidel categories). For the latter situation one has to introduce a ``coefficient system" on 
 the secondary polytope. This will be done in Section 9. Nevertheless the case of ``trivial coefficient system" 
 discussed in Sections 3-8 has most of the features of the general case, so we decided to keep it for
  pedagogical reasons.

 \end{rem}

 \begin{prop}\label{prop:geom-subalg}
 $\gen_A$ is an $L_\infty$-subalgebra in $\bgen_A$, with trivial differential.
 \end{prop}
 
 The $L_\infty$-algebra $\gen_A$ (together with its various generalizations) will be the primary
 object of study in this paper. We will call it the {\em geometric $L_\infty$-algebra} associated to $A$. 

\vskip .2cm

\noindent {\sl Proof of Proposition \ref{prop:geom-subalg}:}
We denote by
\[
\lambda_n: {\bgen}^{\otimes n}\lra\bgen, \,\,\, n=1,2,\cdots, 
\]
the components of the $L_\infty$-structure. Thus $\lambda_1=d$ is the differential,
$\lambda_2(x\otimes y) = [x,y]$ is the bracket, etc. 

It follows from our construction that $\lambda_n$ corresponds to coarse subdivisions of marked subpolytopes
of $(Q,A)$ into precisely $n$ marked subpolytopes. More precisely, choose a generator (canonical up to sign)
$e_{A'}\in\orr(\Sigma(A'))$ for each $A'\subset A$, $|A'|=d+1$. Thus, the $e_{A'}$ form a basis in $\gen$. 
Then each coarse subdivision
 $\Pc''=\bigl\{ (Q''_\nu, A''_\nu)\}_{\nu=1}^n$ of each $(Q', A')$  
 contributes a matrix element
 \be\label{eq:mat-el}
 \bigl\langle e_{A''_1} \otimes\cdots \otimes e_{A''_n} \bigl| \,\, \lambda_n\,\, \bigl| e_{A''}\bigr\rangle \,\,=\,\,\pm 1, 
 \ee
 and all non-zero matrix elements of $\lambda_n$ are obtained in this way. Now, if each $(Q''_\nu, A''_\nu)$ is
 geometric, i.e., $A''_\nu = A\cap Q''_\nu$, the $A'=A\cap Q'$ as well, i.e., $(Q', A')$ is geometric,
 cf. Proposition \ref {prop:geom}. This means that $\lambda_n(\gen^{\otimes n}) \subset \gen$, i.e., that
 $\gen$ is an $L_\infty$-subalgebra. 
 The fact that the differential in $\gen$ vanishes, follows from Example \ref{ex:3.3}(a) below.\qed

\begin{ex} Recall that we have assumed  that $A$ is in affinely general position, i.e.  each $(d+1)$-element subset
generates a $d$-simplex. Hence 
 the low degree part of  the $L_\infty$-algebra $\gen$ looks as follows. We have $\gen^{\leq 0} = 0$, while
$\gen^1$ is spanned by $d$-dimensional simplices, i.e., by marked subpolytopes
 $(Q', A')$ with $|A'|=d+1$. The degree 2 part $\gen^2$
is spanned by circuits ($|A'|=d+2$) and so on. 
\end{ex}

 \begin{exas}\label{ex:3.3}
 (a) Differential in $\bgen$. It corresponds to Example \ref{ex-faces-sigma}(c1) of coarse subdivisions. That is, given
 $A'\subset A$ and $\omega\in A'$ which is not a vertex of $Q'=\Conv(A')$, we have, specializing 
 \eqref{eq:mat-el}, that
 \[
 \bigl\langle e_{A'-\{\omega\}} \bigl| \,\, d \,\, \bigl| e_{A'} \bigr\rangle \,\,=\,\,\pm 1. 
 \]
 Writing $B$ for $A'-\{\omega\}$ (which can be arbitrary) and collecting the matrix elements (i.e., adding, not subtracting one point
 in all admissible ways), we have, cf. Fig. \ref{fig:d-bgen}:
 \[
 d(e_B) \,\,=\,\,\sum_{\omega\in (A\cap Q')-B} \pm e_{B\cup\{\omega\} }, \quad Q'=\on{Conv}(B). 
 \]
 
 \bef
 \[ 
 \btp[scale=0.4, baseline=(current  bounding  box.center)]
\node (1) at (0,0){};
\node (2) at (3,0){};
\node (3) at (5,1){};
\node (4) at (2,6){};
\node (5) at (-1,2){};
\fill (1) circle (0.15); 
\fill (2) circle (0.15); 
\fill (3) circle (0.15); 
\fill (4) circle (0.15); 
\fill (5) circle (0.15); 
 
 \draw (0,0) -- (3,0) -- (5,1) -- (2,6) -- (-1,2) -- (0,0);
  
\etp 
\hskip 1cm 
\btp[scale=0.4, baseline=(current  bounding  box.center)]
\draw[->] (0,3) --(3,3); 
\node at (1.5, 4){$d$};
\etp \hskip 1cm 
\sum_\bullet \quad \pm \quad
\btp[scale=0.4, baseline=(current  bounding  box.center)]
\node (1) at (0,0){};
\node (2) at (3,0){};
\node (3) at (5,1){};
\node (4) at (2,6){};
\node (5) at (-1,2){};
\fill (1) circle (0.15); 
\fill (2) circle (0.15); 
\fill (3) circle (0.15); 
\fill (4) circle (0.15); 
\fill (5) circle (0.15); 
 
 \draw (0,0) -- (3,0) -- (5,1) -- (2,6) -- (-1,2) -- (0,0);
 \node at (2,2){$\bullet$};
 
\etp
\]
 \caption {Differential in $\bgen$. 
}
\label{fig:d-bgen}
\enf
In other words, fixing a subpolytope $Q'\subset (Q,A)$, all $e_B$ with $\Conv(B)=Q'$, span a subcomplex in $(\bgen, d)$. 
This subcomplex is isomorphic to the augmented cochain complex of the simplex $\Delta^I$, whose set of vertices is
\[
I=(A\cap Q') - \on{Vert}(Q'). 
\]
In particular, it is exact, if $I\neq\emptyset$, i.e., if $Q'$ contains points of $A$ other than its vertices. This exactness
is another reason why $\gen$ is more important for us than $\bgen$. 

Note that if $B=A\cap Q'$, then $d(e_B)=0$ (there are no more points to add), so the differential in $\gen$ vanishes. 
 
 \vskip .2cm
 
 (b) The binary bracket in $\gen$ corresponds to  Example \ref{ex-faces-sigma}(c2) of coarse subdivisions. 
 That is, if $(Q',A')$ is subdivided by a hyperplane $H$ into $(Q''_1, A''_1)$ and $(Q''_2, A''_2)$, then
 \[
 [e_{A''_1}, e_{A''_2}] = \pm e_A'
 \]
 
 (c) Consider the particular case  $d=2$ and assume that $A\subset\RR^2$ is {\em in convex position},
 i.e., each element of $A$ is a vertex of $Q=\Conv(A)$. Then $\gen=\bgen$ and only the binary bracket is present, so
 $\gen$ is a graded Lie algebra in a more familiar sense (no differential, no higher $\lambda_n$). 
 As an illustration, consider the further particular case when $A$ consists of 4 points, forming the vertices
 of a convex 4-gon. Then $\gen^1$ has dimension 4, with the basis vectors corresponding to the 4 triangles
 $a,b,c,d$ in Fig.\ref{fig:4-gon}. The space $\gen^2$ is 1-dimensional, with the basis vector   corresponding to
 the 4-gon itself, while $\gen^{\geq 3}=0$. Thus $\gen$ is  has five basis vectors:  $e_a, e_b, e_c, e_d$
 of degree 1 and  $e_A$ of degree 2,   with the only non-zero brackets being
  
 \[
 [e_a, e_b]=[e_c ,e_d]=e_A.
 \]

 \bef
 \centering
 \btp[scale=.3, baseline=(current  bounding  box.center)]
\node (1) at (0,0){};
\fill (1) circle (0.15); 

\node (2) at (5,0){};
\fill (2) circle (0.15);

\node (3) at (5,5){};
\fill (3) circle (0.15);

\node (4) at (0,5){};
\fill (4) circle (0.15);

 \etp
\hskip 1cm
 \btp[scale=.3, baseline=(current  bounding  box.center)]
\node (1) at (0,0){};
\fill (1) circle (0.15); 

\node (2) at (5,0){};
\fill (2) circle (0.15);

\node (3) at (5,5){};
\fill (3) circle (0.15);

\node (4) at (0,5){};
\fill (4) circle (0.15);

\draw (0,0) -- (0,5) -- (5,5) -- (5,0) -- (0,0);
\draw (0,0) -- (5,5); 

\node at (1.5,3.5){$a$};
\node at (3,1.5){$b$};

\etp
\hskip 1cm
 \btp[scale=.3, baseline=(current  bounding  box.center)]
\node (1) at (0,0){};
\fill (1) circle (0.15); 

\node (2) at (5,0){};
\fill (2) circle (0.15);

\node (3) at (5,5){};
\fill (3) circle (0.15);

\node (4) at (0,5){};
\fill (4) circle (0.15);

\draw (0,0) -- (0,5) -- (5,5) -- (5,0) -- (0,0);
\draw (0,5) -- (5,0); 

\node at (1.5,1.5){$c$};
\node at (3.5, 3.5){$d$};

\etp

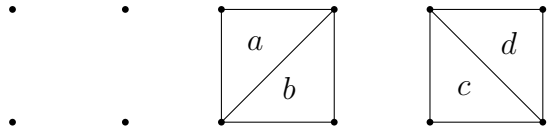
 \caption {Four points in convex position}
 \label{fig:4-gon}
 \enf
 
 (d) Take $A\subset\RR^2$ to be a 4-element circuit consisting of the 3 vertices of a triangle $Q$ and one
 more point $\omega'$ inside $Q$. Marked subpolytopes of $(Q,A)$ include 3 small triangles $a,b,c$ in Fig. 
 \ref{fig:4-pt-nonconv}, one large triangle $f=(Q, A-\{\omega'\})$, and $(Q,A)$ itself. In this case $\bgen$
 has two non-zero  components: $\bgen^1$, with basis $e_a, e_b, e_c, e_f$ and $\bgen^2$ with basis $e_A$.
  The only non-trivial parts of the $L_\infty$-structure are the differential and the ternary bracket whose
  nontrivial matrix elements are given by:
  \[
  d(e_f)=e_A, \,\,\, \lambda_3(e_a, e_b, e_c) = e_A,
  \]
  so the cohomology (w.r.t. $d$) of $\bgen$ has all the multiplications $\lambda_n$ being $0$.
   On the other hand, $\gen$ is spanned by
  $e_a, e_b, e_c, e_A$ and has $\lambda_3\neq 0$. 
 
 \bef
 \centering

 \btp[scale=.4, baseline=(current  bounding  box.center)]
 \node (1) at (-3,0){};
 \fill (1) circle (0.15);
 
 \node (2) at (3,0){};
 \fill (2) circle (0.15);
 
 \node (3) at (0,5){};
 \fill (3) circle (0.15);
 
 \node (4) at (0,2){};
 \fill (4) circle (0.15);
 
 \draw (-3,0) --(3,0) -- (0,5) -- (-3,0); 
 \draw (-3,0) -- (0,2) -- (0,5);
 \draw (0,2) -- (3,0); 
 
 \node at (-1,2.2){$a$};
 \node at (1,2.2){$b$}; 
 \node at (0,1){$c$};
 \etp
  \hskip 1cm
 \btp[scale=.4, baseline=(current  bounding  box.center)]
 \node (1) at (-3,0){};
 \fill (1) circle (0.15);
 
 \node (2) at (3,0){};
 \fill (2) circle (0.15);
 
 \node (3) at (0,5){};
 \fill (3) circle (0.15);
 
 \node (4) at (0,2){$f$};
 
 \draw (-3,0) --(3,0) -- (0,5) -- (-3,0); 
 \etp
  \hskip 1cm
 \btp[scale=.4, baseline=(current  bounding  box.center)]
 \node (1) at (-3,0){};
 \fill (1) circle (0.15);
 
 \node (2) at (3,0){};
 \fill (2) circle (0.15);
 
 \node (3) at (0,5){};
 \fill (3) circle (0.15);
 
 \node (4) at (0,2){};
 \fill (4) circle (0.15);
 
 \draw (-3,0) --(3,0) -- (0,5) -- (-3,0); 
 \node at (-1,1){$A$};
 \etp

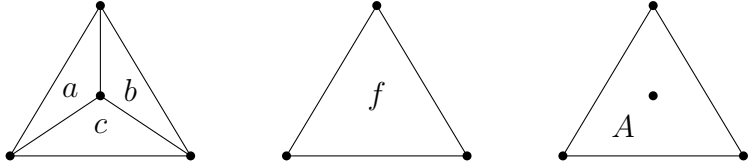
 \caption {Four points in non-convex position. }
 \label{fig:4-pt-nonconv}
 \enf
 
 Recall \cite{getzler} that an $L_\infty$-algebra $\Lc$ is called {\em nilpotent}, if there exists $r_0>0$ such that
 all $r$-ary iterated superpositions of the $\lambda_n$ are identically zero as maps $\Lc^{\otimes r}\to\Lc$
 for $r>r_0$.

 \begin{prop}
 The $L_\infty$-algebra $\bgen$ (and therefore  $\gen$) is nilpotent. 
 \end{prop}
 
 \noindent {\sl Proof:} Since the matrix elements of 
 $\lambda_n$ corresponds to {\em coarse} subdivisions of marked subpolytopes $(Q', A')\subset (Q,A)$
 into $n$ marked subpolytopes, the matrix elements of 
 any $r$-ary iterated superposition of the $\lambda_n$ correspond to
 some (not necessarily coarse) subdivisions into $r$ marked subpolytopes and therefore vanish for $r\geq |A|$. \qed
  
  \end{exas}. 
 
 \vfill\eject

 \section{Maurer-Cartan elements in $\gen$.}\label{sec:MC}
 
 Let $\Lc$ be a nilpotent $L_\infty$-algebra over $\k$. Recall that a 
    {\em Maurer-Cartan element} in $\Lc$   is an  element $\gamma\in\Lc^1$
 such that
 \be\label{eq:MC}
 \sum_{n=0}^\infty\,\,  \frac {1} {n!} \lambda_n (\gamma, \cdots, \gamma) \,\,=\,\, 0 \,\, \in \,\,\Lc^2. 
 \ee
 Maurer-Cartan elements form a subset $\on{MC}(\Lc)$ in
  $\Lc^1$. If $\dim\,\Lc^1<\infty$, we can consider $\Lc^1$ as an affine space over $\k$ and view $\on{MC}(\Lc)$
  as an affine subscheme in $\Lc^1$.

 If $\Lc$ is a general $L_\infty$-algebra, one can define the {\em formal scheme of Maurer-Cartan elements}
 $\widehat{\on{MC}}(\Lc)$ which is the functor on local  Artin $\k$-algebras $\Lambda$ defined by
 \[
  \Lambda \mapsto \on{MC}(\Lc\otimes_k\men),
 \]
 where $\men$ is the maximal ideal of $\Lambda$. 
 
 \vskip .2cm
 Such a functor represents  a geometric object called a {\em formal pointed dg-manifold} in \cite{defbook}, Ch.3.
 
 We are interested in Maurer-Cartan elements in the  finite-dimensional nilpotent $L_\infty$-algebra $\gen = \gen_A$. 
 In this section we assume that $A$ is in general position. 
 As we saw, $\gen^1$ is spanned by basis vectors $e_\sigma$ corresponding to all possible marked  $d$-simplices
$\sigma = (\Conv(A'), A')$ with $A'\subset A$, $|A'|=d+1$. Such a marked simplex is  nothing but a straight $d$-simplex
on vertices from $A$. 
Thus we can view
\[
\gamma \,\,=\,\,\sum \gamma_\sigma \cdot e_\sigma, \quad \gamma_\sigma\in\k
\] 
as a ``$d$-cochain", i.e., a datum of numbers $\gamma_\sigma \in\k$ attached to all marked simplices $\sigma$. 

Now, $\gen^2$ is spanned by basis vectors $e_Z$ corresponding to all possible circuits $Z\subset A$, $|Z|=d+2$.
The coefficient of the LHS of \eqref{eq:MC} at $e_Z$ is then equal to the sum of two terms corresponding to the 
two triangulations $T_+, T_-$ of $\Conv(Z)$, as in Example \ref{ex-faces-sigma}(b). More precisely, this coefficient is
equal to
\[
\prod_{\sigma\in T_+} \gamma_\sigma \,\,-\,\, \prod_{\sigma\in T_-} \gamma_\sigma. 
\]
Therefore we obtain:

\begin{prop}\label{prop:mc-g}
 $\gamma$ is a Maurer-Cartan element in $\gen$, if and only if it is a ``cocycle" in the multiplicative sense:
for any circuit $Z$ with triangulations $T_+, T_-$ of $\Conv(Z)$ we have
\[
\prod_{\sigma\in T_+} \gamma_\sigma \,\,=\,\, \prod_{\sigma\in T_-} \gamma_\sigma. \quad\quad\quad \qed
\]

\end{prop}

This means that for any marked subpolytope $(Q',A')$ of $(Q,A)$ we have an element $\gamma_{Q',A'}\in\k$
defined in terms of any chosen regular triangulation $\Tc'$ of $(Q',A')$:
\[
\gamma_{Q',A'} \,\,=\,\,\prod_{\sigma\in\Tc'} \gamma_\sigma,
\]
with the RHS being independent of the choice of $\Tc'$. 

\begin{ex}
Let $\k=\CC$ and let $\Omega$ be any smooth  $d$-form on $Q$. Then putting
\[
\gamma_\sigma = e^{\beta_\sigma}, \quad \beta_\sigma = \int_\sigma \Omega, 
\]
we get a Maurer-Cartan element in $\gen$. 
\end{ex}

One way of interpreting Proposition \ref{prop:mc-g}  is by saying that the subscheme $\on{MC}(\gen)$
is given inside the affine space $\gen^1$ by  a system of {\em binomial equations}: each equation had the form
of equality of exactly two monomials. This means that the scheme $\on{MC}(\gen)$ is in fact a toric variety.

More precisely, let $\Sen_d(A)$ be the set of all $d$-simplices on vertices from $A$. Thus
$\Sen_d(A)$ is a subset in the set $A\choose d+1$ of all $(d+1)$-element subsets in $A$,
and $\Sen_d(A) = {A\choose d+1}$ when $A$ is in affinely general position. 
Let $\ZZ^{\Sen_d}(Q,A)$ be the free abelian group of {\em combinatorial $d$-chains}  in $(Q,A)$,
i.e., of integer linear combinations of $d$-simplices. We write an element $\beta\in \ZZ^{\Sen_d}(Q, A)$
as a system $\beta = (\beta_\sigma)$ where $\sigma\in\Sen_d(A)$ is a $d$-simplex, and $\beta_\sigma\in\ZZ$. 
Let $\Zen_d(Q,A)\subset \ZZ^{\Sen_d(Q,A)}$ be the  sublattice of {\em $d$-cocycles}, i.e., of $\beta$ as above such
that for each circuit $Z$ as in Proposition \ref{prop:mc-g}, we have
\[
\sum_{\sigma\in T_+} \beta_\sigma \,\,=\,\, \sum_{\sigma\in T_-} \beta_\sigma.  
\]
Let $\Gamma_A = \Zen_d(Q,A)\cap \ZZ_+^{\Sen_d(Q,A)}$. This is a finitely generated abelian semigroup.
The monomial equations for $\on{MC}(\gen)$ translate into the following. 

\begin{cor}
The scheme $\on{MC}(\gen)$ is isomorphic to $\on{Spec}\, \k[\Gamma_A]$, the spectrum of the
semigroup algebra of $\Gamma_A$. \qed
\end{cor}

\vfill\eject

\section {1-dimensional case: refinement to an $A_\infty$-algebra.}
\label{sec:1dim}

In this section we consider the simplest case $d=1$. 
So we assume that $A=\{\omega_1 < \cdots < \omega_r\} \subset\RR$. Then $Q=\Conv(A)$ is the interval 
$[\omega_1, \omega_r]$.  
Let us form the vector space  $V$ as in \eqref {eq;space-V}.

If $(Q', A')$ is a marked subpolytope  of $(Q,A)$ and
$\Pc''= \{ (Q''_\nu, A''_\nu)\}$ is a polyhedral subdivision of $(Q', A')$, then  each $Q''_\nu$ is a sub-interval in $Q'$
and so {\em there is a natural
order} on the set of the intervals $Q''_\nu$, induced by the order on $\RR$  (i.e., geometrically, from left to right).
This means that we can lift the differential $d$ from $S^\bullet(V)$ (free commutative algebra generated by $V$)
 to a differential in  $T^\bullet(V)$, the tensor algebra (free associative algebra) generated by $V$. More precisely, 
  we realize the vector space $V_{\Pc''}$  from \eqref{eq:V-P}
   as a subspace in $T^\bullet(V)$:
\[
V_{\Pc''} \,\,=\,\,\bigotimes_\nu \,\,  V_{A'_\nu} \,\,\subset \,\, T^\bullet(V)
\]
(embedding given by the tensor multiplication in the order from left to right). We then define the action of 
$d$ on the generators  of $T^\bullet(V)$ by making it act on a summand $V_{A'} \subset V$ as
  the top degree part of the chain differential in $C_\bullet(\Sigma(A'))$, now  considered as a map $\sum_{\Pc''} d_{\Pc''}$,
  where $\Pc''$ runs over all coarse subdivisions of $(Q', A')$ and
  \[
  d_{\Pc''}:  V_{A'} \lra V_{\Pc''} \,\,\subset \,\, T^\bullet(V)
  \]
  is induced by the chain differential in $C_\bullet(\Sigma(A'))$.
  Then we extend $d$ to the whole $T^\bullet(V)$ by the Leibniz rule. Similarly to 
  Proposition 
\ref{prop:d2=0}, we then get:  

  \begin{prop}
The differential $d$ in $T^\bullet(V)$ satisfies $d^2=0$ and so makes $T^\bullet(V)$ into an
associative dg-algebra. The canonical commutativization homomorphism $(T^\bullet(V), d)
\to S^\bullet(V), d)$ is a morphism of associative dg-algebras.  \qed
 \end{prop}
 
 Clearly $d$ preserves the ideal $T^\bullet_+(V):=\oplus_{n\ge 1}V^{\otimes n}\subset T^\bullet(V)$.

 This means that the $L_\infty$-algebra $\bgen$ lifts to  an $A_\infty$ algebra  with the
 same underlying vector space by the $A_\infty$-analog of the formula
 \[
[x,y] = xy - yx.
\]
  We denote this $A_\infty$-algebra by $\bR$. As before, we have a similar lifting of the subalgebra $\gen$.
  We denote this lifting by $R$ and it will be more fundamental for us than $\bR$.

  Thus $R$ has basis $e_{ij}$, $1\leq i<j\leq r$, corresponding to the geometric marked subpolytopes
  \[
  \bigl( [\omega_i,\omega_j],\,\, \{\omega_i, \omega_{i+1}, \cdots, \omega_j\}\bigr). 
  \]
  The degree of $e_{ij}$ is $j-i$, and the binary multiplication is given by
  \[
  e_{ij} \cdot e_{jk} = e_{ik}, \quad i<j<k,
  \]
  and all other binary products as well as other $A_\infty$-data vanish. 
In other words, $R$ is the associative algebra without unit, isomorphic to $T_n(\k)$, the algebra
of strictly upper triangular $r$ by $r$ matrices with zeroes on and below the diagonal. This triangular
structure is the simplest instance of the relation between our formalism and deformations of
semi-orthogonal collections, see below.

\vfill\eject

\section{Relative setting: one point ``at infinity". }\label{sec:rel-set}

Let $(Q,A)$ be a marked polytope in $\RR^d$, as above. Let us add to $A$ one more point
$\infty\in\RR^d$, forming $\widetilde A=A\cup \{\infty\}$ and $\widetilde Q=\Conv(\wt A)$.
We assume that $\infty$ is a new vertex of $\wt Q$, i.e., that it lies ``far outside of $A$". 
We also assume that the set $\wt A$ is in general position. 
Let $\wt\gen = \gen_{\wt A}$ be the geometric $L_\infty$-algebra corresponding to $\wt A$. 

\begin{figure}[h!]
\centering
 \btp[scale=.4, baseline=(current  bounding  box.center)]
 
 \node (1) at (0,0){};
 \fill (1) circle (0.15); 
 
  \node (2) at (4,0){};
 \fill (2) circle (0.15); 
 
  \node (3) at (6,3){};
 \fill (3) circle (0.15); 
 
  \node (4) at (3,5){};
 \fill (4) circle (0.15); 
 
  \node (5) at (-2,3){};
 \fill (5) circle (0.15); 
 
 \draw (0,0) -- (4,0) -- ( 6,3) -- (3,5) -- (-2,3) -- (0,0); 
 
 \draw (-2,3) -- (-2,12);
 
 \draw (6,3) -- (6,12); 
 
 \node at (4,10.5){$\wt Q$};
 
 \node at (2,13){$\cdots \infty\cdots$};
 
 \node at (2,2){$Q$};

 \etp
 \hskip 2cm
 \btp[scale=.4, baseline=(current  bounding  box.center)]
 
 \node (1) at (0,0){};
 \fill (1) circle (0.15); 
 
  \node (2) at (4,0){};
 \fill (2) circle (0.15); 
 
  \node (3) at (6,3){};
 \fill (3) circle (0.15); 
 
  \node (4) at (3,5){};
 \fill (4) circle (0.15); 
 
  \node (5) at (-2,3){};
 \fill (5) circle (0.15); 
 
 \draw (0,0) -- (4,0) -- ( 6,3) -- (3,5) -- (-2,3) -- (0,0); 
 
 \node (6) at (2,12){};
 \fill (6) circle (0.15);
 \draw (-2,3) -- (2,12);
 \draw (6,3) -- (2,12);

 \node at (4,10.5){$\wt Q$};
 
 \node at (2,13){$ \infty$};
 
 \node at (2,2){$Q$}; 
 
 \etp

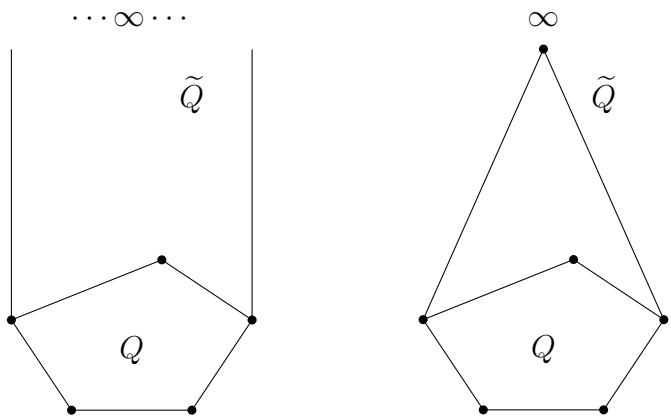
\caption{$\infty$ as a point at the projective infinity vs.  a finite point far away.}
\label{fig:extended-polyt}
\enf

\begin{rem} We can, if we want, imagine $\infty$  to be an actual ``infinite"  point of the projective completion
$\RR \PP^{d}= \RR^d \cup \RR \PP^{d-1}$, say
 the point of $\RR \PP^{d-1} $ with homogeneous coordinates $[0:0: \cdots : 1]$.  
 Note that the concepts of hyperplanes, triangulations, convexity etc. are invariants under projective transformations
 of $\RR \PP^{d}$. So there is no combinatorial difference between a bounded polytope $\wt Q$ on the right and
 the unbounded polyhedron $\wt Q$ on the left of Fig. \ref{fig:extended-polyt}. 
 \end{rem}
 
 A subpolytope $Q' \subset (\wt Q, \wt A)$  will be called {\em finite}, if $\on{Vert}(Q')\subset A$, and
 {\em infinite}, otherwise, i.e., if $\infty\in Q'$. The $L_\infty$-algebra $\wt\gen$ splits, as a vector space, into a direct sum
 \[
 \wt\gen = \gen\oplus \gen_\infty, 
 \]
 where $\gen=\gen_A$ is spanned by the summands
 \[
E_{A'} \,\,=\,\, \orr(\Sigma(Q'\cap A))^*[-1-\dim\Sigma(Q'\cap A)]
 \]
 for finite $Q'$,  and $\gen_\infty$ is spanned by similar summands for infinite $Q'$.

 \begin{prop}\label{prop:gentilde}
 Both $\gen$ and $\gen_\infty$ are $L_\infty$-subalgebras of $\wt\gen$, and, moreover,
 $\gen_\infty$ is an $L_\infty$-ideal. In other words. $\wt\gen$ is the semi-direct product
 $\wt\gen = \gen \ltimes \gen_\infty$.

 \end{prop}
 
 \noindent {\sl Proof:}  The $L_\infty$-brackets are given by combining  subpolytopes, say,
 $Q''_1, \cdots, Q''_n$,
 together into a coarse subdivision. If neither of the $Q''_i$ contains $\infty$, then the combined marked 
 polytope $Q'$ does not contain it either. This means that $\gen$ is an $L_\infty$-subalgebra
 (it is nothing but $\gen_A$, the geometric $L_\infty$-algebra associated to $A$ alone). 
 If at least one of the $Q''_i$ does contain $\infty$, then so does $Q'$. This means that
 $\gen_\infty$ is a $L_\infty$-ideal. \qed
 
 \vskip .2cm
 
\begin{rems}\label{rems:olA}
 (a) 
 Let us choose a Euclidean subspace $\RR^{d-1}\subset\RR^d$ not passing through $\infty$
 and let
 \[
 p: \RR^d-\{\infty\} \lra \RR^{d-1}
 \]
 be the projection obtained by drawing straight lines through $\infty$ and intersecting them with $\RR^{d-1}$.
 In the ``projective infinity" picture on the left of Fig. \ref{fig:extended-polyt}, we can think of  $p: \RR^d\to\RR^{d-1}$
 as being the projection forgetting the last coordinate. 
 
 We then have the point configuration $\ol A= p(A)$ in $\RR^{d-1}$, and by our assumption,
 $\ol Q=\Conv(\ol A)$ coincides with the image $p(Q)$. 
 It is instructive to compare
 the $L_\infty$-algebra $\gen_\infty$  with $\gen_{\ol A}$, the $L_\infty$-algebra associated to
 $\ol A\subset \RR^{d-1}$.

More precisely,   each infinite subpolytope $Q'\subset (\wt Q, \wt A)$
gives a subpolytope $p(Q')\subset (\ol Q, \ol A)$ via the projection $p$.
However, not every subpolytope of  $(\ol Q, \ol A)$ is obtained in this way.
Similarly, each coarse subdivision of $Q'$ gives a coarse subdivision of $p(Q')$
but not the other way around. Nevertheless, the $L_\infty$-algebra $\gen_\infty$
has, in many respects, $(d-1)$-dimensional nature. 
 
 \vskip .2cm

(b) The two point configurations $A, \ol A$ associated to $\wt A$ and an element $\infty\in\wt A$
 are precisely the two {\em minors} of $\wt A$ studied by Billera, Gelfand and Sturmfels \cite{BGS}.
 More precisely, $A$ is the {\em minor by deletion} and $\ol A$ is the {\em minor by contraction},
 in the terminology of \cite{BGS}. Among other things, it was proved in \cite{BGS}
 that $\Sigma(\ol A)$ is a face of $\Sigma(\wt A)$, and $\Sigma(A)$ is a Minkowski summand
 of $\Sigma(\wt A)$. 
\end{rems}

 \begin{rem}
 
 If we assume that $\wt\gen$ is an ordinary dg-Lie algebra, i.e. if $\lambda_{\geq 3} =0$,
 then the Lie algebra structure on $\wt\gen$ gives a morphism of dg-Lie algebras
 \[
 \gen\lra \Der(\gen_\infty), \quad x\mapsto \on{ad}(x): \gen_\infty\to\gen_\infty. 
 \]
 In general, what we get is a ``derived" analog of such a morphism, as we explain in the next section. 
 
 \end{rem}
 
 \vfill\eject

 \section{Hochschild complexes and derived derivation spaces. }\label{sec:hoch}
 
 In this section we recall  some general principles of deformation theory and of Koszul duality.  
 
 Let $\scrP$ be any dg-operad (over $\k$), and $\Ac$ be any $\Ps$-algebra (in the category of
 dg-vector spaces over $\k$).  Automorphisms of $\Ac$ form a   group $\Aut_\Ps(\Ac)$, and derivations
 of $\Ac$ (of all degrees) form a dg-Lie algebra $\Der_\Ps(\Ac)$. One can form the derived functor
 of $\Ac\mapsto \Der_\Ps(\Ac)$, defining
 \[
 R\Der_\Ps(\Ac) \,\,=\,\,\Der_\Ps(\wt\Ac), 
 \]
 where $\wt\Ac\to\Ac$ is a ``cofibrant resolution", i.e., a quasi-isomorphism of $\Ps$-algebras
 with $\wt\Ac$ being free as a $\Ps$-algebra without differential.
 RHS of this formula can be understood (at least in $char(\k)=0$ case) as an $L_{\infty}$-algebra. There is a canonical choice of a cofibrant resolution called the Boardmann-Vogt resolution in \cite{defoperads}. It was used in  {\em loc. cit.} for the description  of general deformation theory of algebras over operads.
 
 \begin{rem}
 If $\Ps\to\Qs$ is a morphism of (dg-)operads, then any $\Qs$-algebra $\Ac$ can be also considered as a $\Ps$-algebra
 (for example, we have $\Ac ss\to\Cc om$, so a commutative algebra can be considered as a particular
 case of an associative algebra). In this case $\Der_\Ps(\Ac)=\Der_\Qs(\Ac)$, but
  $R\Der_\Ps(\Ac)$  can be different from $R\Der_\Qs(\Ac)$. 
  \end{rem}
 
 According to the general philosophy of deformation theory (see e.g. \cite{defbook}) the aim of formal deformation theory in the case when $char({\bf k})=0$ is the construction of formal pointed dg-manifold (equivalently, $L_{\infty}$-algebra) which describes the local structure of the moduli space of deformations. More precisely, if ${\mathcal L}$ is such an $L_{\infty}$-algebra, then the  commutative topological algebra generated by ${\mathcal L}^{\ast}[-1]$  carries a structure of a dg-algebra. Its formal spectrum (as a graded topological algebra) can be thought of as a formal graded scheme, endowed with a (graded) vector field $Q$ of degree $+1$ such that $[Q,Q]=0$ (``homological vector field"). The formal scheme of zeros $Z(Q)$ admits a foliation generated by graded vector fields of the form $[Q, \bullet]$ (all schemes are understood as functors on commutative Artin algebras). The ``space of leaves" of the above foliation should be thought of as the moduli space of the deformation problem controlled by ${\mathcal L}$. The reader can look at \cite{defbook} for the details and examples and at
 \cite{defoperads} for the implementation of this approach in the case of deformations of algebras over operads.
 In the more conventional language of stacks one has the following corollary of those considerations.
 
 \begin{prop}[{\bf  Local structure of the moduli stack}]\label{prop:DT}
  
The formal germ of the  deformation stack of $\Ac$ as a $\Ps$-algebra is isomorphic to
the quotient stack
\[
\widehat {\on{MC}}(R\Der_\Ps(\Ac))\bigl/\hskip -.15cm \bigl/ \exp(R\Der^0_\Ps(\Ac)).
\]
Here $ \widehat{\on{MC}}(R\Der_\Ps(\Ac))$ is the formal   scheme of
 Maurer-Cartan elements in the dg-Lie algebra $R\Der_\Ps(\Ac)$, and 
 $ \exp(R\Der^0_\Ps(\Ac))$ is the formal group associated to the Lie algebra $R\Der^0_\Ps(\Ac)$,
 which naturally acts on the formal scheme  $ \widehat{\on{MC}}(R\Der_\Ps(\Ac))$. \qed
 \end{prop}
 
 So we can say that $R\Der_\Ps(\Ac)$ ``governs deformations of $\Ac$ as a $\Ps$-algebra". 
 
 \begin{exa}
 Let $\Ps=\Ass$, and $\Ac$ be an associative algebra. Then $R\Der(\Ac)=R\Der_{\Ass}(\Ac)$ is the shifted and truncated
 Hochschild cochain complex of $\Ac$:
 \[
 \begin{gathered}
C^{\geq 1}(\Ac, \Ac)[1] :=C_{\on{\Ass}}^{\geq 1}(\Ac, \Ac)[1] \,\,=\\
\,\,\biggl\{  \Hom_k(A,A) \buildrel\delta_0\over\lra \Hom_k(A^{\otimes 2}, A) 
 \buildrel \delta_1\over\lra \Hom_k(A^{\otimes 3}, A)\to \cdots\biggr\}. 
 \end{gathered}
 \]
 Thus $\Ker(\delta_0)=\Der(\Ac)$ is the usual space of derivations. 
 The Lie algebra structure on $R\Der(\Ac)$ is given by the ``brace formula". That is, for
 $p\in\Hom(A^{\otimes m}, A)$ and $q\in\Hom(A^{\times n}, A)$ and $i=1, ...,m$ we define
 the {\em $i$th brace} of $p$ and $q$  by
 $p\circ_i q\in \Hom(A^{\otimes(m+n-1)}, A)$ by
 \[
 \begin{gathered}
 (p\circ_i q)(a_1, \cdots , a_{m+n-1}) \,\,=\,\, \\
 =\,\, p\bigl( a_1, \cdots, a_{i-1}, q(a_i, a_{i+1}, \cdots, a_{i+n-1}), a_{i+n}, \cdots a_{m+n-1}\bigr).
  \end{gathered}
 \]
 Then
 \[
 [p,q] = \sum_{i=1}^m \pm  p\circ_i q \,\, -\,\, \sum_{j=1}^n \pm q\circ_j p. 
 \]
 \end{exa}
 
 \begin{exa}
 Let $\Ps=\Lie$ and $\Ac=\Lc$ be a Lie algebra. Then $R\Der_{\Lie}(\Lc)=  $ is the shifted and truncated
 Chevalley-Eilenberg cochain complex of $\Lc$ with coefficients in the adjoint representation:
 \[
C^{\geq 1}_{\on{Lie}}(\Lc, \Lc)[1]\,\,=\,\,
\biggl\{  \Lc^*\otimes \Lc\buildrel \delta_0\over\lra \Lambda^2\Lc^* \otimes \Lc
 \buildrel \delta_0\over\lra \Lambda^3\Lc^* \otimes \Lc\to \cdots\biggr\}
 \]

 \noindent {\bf General principles of Koszul duality:}   Let $\Ps$ be a quadratic Koszul operad
 with Koszul dual operad $\Ps^!$ (see e.g. \cite{GiKa} about basics on Koszul duality for operads). A
 {\em weak $\Ps$-algebra structure} (also called a {\em $\Ps_\infty$-algebra structure}) on a graded
 vector space $\Ac$ is the same as a $\Ps^!$-algebra differential $d$ on the free $\Ps^!$-algebra
 $F_{\Ps^!}(\Ac^*[-1])$, satisfying $d^2=0$. 

We remark that we define the  free $\Ps^!$-algebra
 $F_{\Ps^!}(\Ac^*[-1])$ in such a way that it does not have $0$-ary operations.

     In particular:
 
 \vskip .2cm
 
 \noindent {\bf (a)} An $L_\infty$-structure on a graded vector space $\Lc$ is a multiplicative differential on the (completed) symmetric
 algebra $S^\bullet_+(\Lc^*[-1])$.
 
 \vskip.2cm
 
  \noindent {\bf (b)} An $A_\infty$-structure on a graded vector space $R$ is a multiplicative differential on the (completed) tensor
 algebra $T^\bullet_+(R^*[-1])$.
   \end{exa}

Although technically one should consider completed tensor products (since corresponding structures are defined in dual terms of the differential cofree coalgebras), we will often skip the word ``completed", since our examples are nilpotent in the sense that only finitely many ``higher" operations are non-trivial.
 
 Most important for us will be the following principle.
 
 \begin{prop}\label{prop:AL-mor}
 A datum consisting of:
 \begin{enumerate}
 \item[(1)] A $\Ps_\infty$-algebra $\Ac$.
 
 \item[(2)] An $L_\infty$-algebra $\Lc$.
 
 \item[(3)] An $L_\infty$-morphism $\alpha: \Lc\to R\Der_{\Ps_\infty}(\Ac)$
 
 \end{enumerate}
  is the same as  a datum consisting of:
  \begin{enumerate}
  
 \item[(i)]  A differential $d_{\Lc}$ in the (completed) algebra $S^\bullet(\Lc^*[-1])$ which preserves the maximal ideal $S^\bullet_+(\Lc^*[-1])$.
  
 \item[(ii)]  A differential $d$  in the (completed)  graded space
  \[
  S^\bullet(\Lc^*[-1]) \otimes_\k F_{\Ps^!}(\Ac^*[-1])
  \]
  which is a $\Ps^!$-algebra derivation, satisfies $d^2=0$  as well as the following 
  {\em Leibniz formula}:
  \[
  d(a\otimes x)= d_{\Lc^*[-1]}(a)\otimes x+(-1)^{|a|}a\otimes d(1\otimes x),\quad a\in  S^\bullet(\Lc^*[-1]), 
  \,\,\, x\in F_{\Ps^!}(\Ac^*[-1]). 
  \]

 \end{enumerate}

   \end{prop}

   \vskip .2cm
   
   \noindent {\sl Proof:} We explain a more conceptual argument and then  spell  out additional
   algebraic details in the examples we will be interested in. Suppose data (i)-(ii) are given. Then
   (i) implies that 
   $\Lc$ is an $L_\infty$-algebra, thus giving the datum (2). 
   Further, considering the restriction $d|_{\Ac^*[-1]}$ and projecting the image on
$1\otimes F_{\Ps^!}(\Ac^*[-1])$, we get a $\Ps^!$-algebra differential in $F_{\Ps^!}(\Ac^*[-1])$,
i.e., we get a $\Ps_\infty$-algebra structure on $\Ac$ which we will call the ``initial" structure.
This gives the datum (1). 
Now, by general principles of deformation theory, $\on{Spec}\, S^\bullet(\Lc^*[-1])$ is a formal commutative
dg-manifold (functor on local artinian commutative dg-algebras)
   whose underlying ordinary (non-dg)  formal commutative manifold is $\widehat{\on{MC}}(\Lc)$,
   the formal scheme of Maurer-Cartan elements in $\Lc$. So we denote
   \[
   R\widehat{\on{MC}}(\Lc) \,\,=\,\, \on{Spec}\, S^\bullet(\Lc^*[-1]).
   \]
   Thus our differential $d$ in the tensor product is the same as a family of $\Ps_\infty$-structures on $\Ac$
   parametrized by $R\widehat{\on{MC}}(\Lc)$ which gives the datum (3) by 
   Proposition \ref{prop:DT}. The argument in the opposite direction is similar.  \qed
   
   \begin{ex}[(Particular case $\Ps=\Ass$)] \label{ex:case-ass}
   A datum consisting of: 
    \begin{enumerate}
 \item[(1)] An $A_\infty$-algebra $R$.
 
 \item[(2)] An $L_\infty$-algebra $\Lc$.
 
 \item[(3)] An $L_\infty$-morphism $\alpha: \Lc\to R\Der(R)= {C}^{\geq 1}(R,R)[1]$
 
 \end{enumerate}
is the same as an algebra differential d  in the (completed) tensor product $S^\bullet(\Lc^*[-1])\otimes T^\bullet(R^*[-1])$
which preserves $S^\bullet_+(\Lc^*[-1])\otimes 1$ as well as $1\otimes T^\bullet_+(R^*[-1])$. 

The relations of the data (1) and (2) to the differential $d$ has already been explained in the general proof of
Proposition \ref{prop:AL-mor}. Let us explain how to obtain the datum (3) explicitly. 

We recall, first of all, the 
general  meaning of the concept of an 
$L_\infty$-morphism $\beta: \Lc \to\Lc'$ between two $L_\infty$-algebras $\Lc$ and $\Lc'$ whose
$L_\infty$-operations we denote by $\lambda_n, \lambda_n'$ respectively. 
By definition, $\beta$ is   the same as a  morphism of commutative dg-algebras
\be\label{eq:beta*}
\beta^*: S^\bullet(\Lc'{}^*[1]) \lra S^\bullet (\Lc^*[1]), 
\ee
the datum which spells out to a collection
 of morphisms of graded vector spaces
\be\label{eq:beta_n}
\beta_n: S^n(\Lc[-1]) \lra \Lc'[-1], \quad n\geq 0
\ee
satisfying a chain of compatibility conditions. More precisely, 
  any collection of the $\beta_n$  in \eqref{eq:beta_n} defines a unique morphisms
 of {\em commutative graded algebras} $\beta^*$  as in 
 \eqref{eq:beta*}
  and the conditions on the $\beta_n$
 say that $\beta^*$ thus defined, actually commutes with the differentials. 
 Note, in particular, that $\beta_1$ can be seen, after shift, as a morphism
 of graded vector spaces
 $
 \beta_1[1]: \Lc \lra \Lc'
 $
 and the first in the chain of conditions on the $\beta_n$ says that $ \beta_1[1]$ commutes with the
 differentials $\lambda_1$ and $\lambda_1'$ in $\Lc$ and $\Lc'$. This is the ``underlying morphism"
 of the $L_\infty$-morphism $\beta$, and higher $\beta_n$ can be seen as higher homotopies for its
 compatibility with the brackets. 
 
 We now specialize to our situation: $\Lc'=C^{\geq 1}(R,R)[1]$. Suppose we are given a differential $d$ in 
$S^\bullet(\Lc^*[-1])\otimes T^\bullet(R^*[-1])$
which preserves $S^\bullet(\Lc^*[-1])$. We describe the corresponding $L_\infty$-morphism $\alpha$ 
in terms of its components
\[
\alpha_n: S^n (\Lc[-1]) \lra C^{\geq 1}(R,R)[1] \,\,=\,\,\bigoplus_{m\geq 1} \Hom_k(R^{\otimes m}, R) [1-m]. 
\]
More precisely, the matrix element 
\[
\alpha_{n,m}: S^n(\Lc[-1]) \to \Hom_\k(R^{\otimes m}, R)[1-m]
\]
of $\alpha_n$ with the respect to the above direct sum decomposition of the RHS
is obtained as the dual of the matrix element of the restriction $d|_{R^*[-1]}$ which we denote
\[
d_{n,m}: R^*[-1] \lra S^n(\Lc^*[-1]) \otimes T^m(R^*[-1]), \quad T^m(R^*[-1]) = (R^*)^{\otimes m} [-m]. 
\]
 \end{ex}
 
 \begin{ex}[(Particular case $\Ps=\Lie$)] 
 A datum consisting of: 
    \begin{enumerate}
 \item[(1)] An $L_\infty$-algebra $\Lc_\infty$.
 
 \item[(2)] An $L_\infty$-algebra $\Lc$.
 
 \item[(3)] An $L_\infty$-morphism $\alpha: \Lc\to R\Der_{\Lie}(R)= C_{\on{Lie}}^{\geq 1}(\Lc, \Lc)[1]$
 
 \end{enumerate}
is the same as an algebra differential $d$   in the (completed) symmetric algebra
\[
S^\bullet(\Lc^*[-1]) \otimes S^\bullet(\Lc_\infty^*[-1]) \,\,=\,\, S^\bullet\bigl( (\Lc\oplus\Lc_\infty)^*[-1]\bigr), 
\]
preserving $S^\bullet_+(\Lc^*[-1])\otimes 1$ as well as $1\otimes S^\bullet_+(\Lc_\infty^*[-1])$. The correspondence between the matrix elements of $d$ and
the matrix elements of the components  $\alpha_n$ is completely analogous to Example
\ref{ex:case-ass}, and we leave it to the reader.

Note that the data of a differential $d$ as above,   is the same as an $L_\infty$-structure on $\Lc\oplus\Lc_\infty$
such that $\Lc$ is an $L_\infty$-subalgebra and $\Lc_\infty$ is an $L_\infty$-ideal. 
This is precisely the kind of situation we had in \S \ref{sec:rel-set}. So we obtain the following.
 \end{ex}
 
 \begin{cor}\label{cor:L-infty-mor} 
 In the situation of \S \ref{sec:rel-set}, we have a natural $L_\infty$-moprhism
 \[
 \alpha: \gen \lra R\Der_\Lie(\gen_\infty).\quad\quad\quad \qed
 \]

 \end{cor}
 
 \vfill\eject

 \section{ Refinement in $d=2$: relative setting. }\label{sec:refi}
 
 Consider the relative setting of \S \ref{sec:rel-set} of which we keep the notation: $\wt A= A\cup\{\infty\}$ etc. 
 
 In this case the $L_\infty$-algebra $\gen_\infty$ can be upgraded to an $A_\infty$-algebra $R_\infty$,
 similar to \S \ref{sec:1dim}. This is because, 
 for each purely infinite coarse
 subdivision $\Pc''=\{(Q''_\nu, A''_\nu)\}$ of an infinite marked subpolygon $(Q', A')$, $\infty\in A'$, the
 infinite marked polygons $Q''_\nu$, are
 naturally ordered using the chosen orientation of $\RR^2$ and the projection $p: \RR^2-\{\infty\}\to\RR$
 as in Remark \ref{rems:olA}(a). 

\bef
\centering
\btp[scale=.4]
\usetikzlibrary{calc}
\def\centerarc[#1](#2)(#3:#4:#5)  { \draw[#1] ($(#2)+({#5*cos(#3)},{#5*sin(#3)})$) arc (#3:#4:#5); }

\node (1) at (0,0){};
 \node (2) at (4,0){};
 \node (3) at (6,2){};
 \node (4) at (2,9){};
 \node (5) at (-2,2){};
 
 \fill (1) circle (0.15); 
  \fill (2) circle (0.15);
   \fill (3) circle (0.15);
    \fill (4) circle (0.15);
     \fill (5) circle (0.15);
     
     \draw (0,0) -- (4,0) -- (6,2) -- (2,9) -- (-2,2) -- (0,0); 
     
     \draw (0,0) -- (2,9) -- (4,0); 
     
 \node at (-.5,2.5){$1$}; 
 \node at (2,2){$2$}; 
 \node at (4.5, 2.5){$3$}; 
 \node at (2,10){$\infty$}; 
 
 \centerarc[->](2,9)(200:350:2)

\etp 
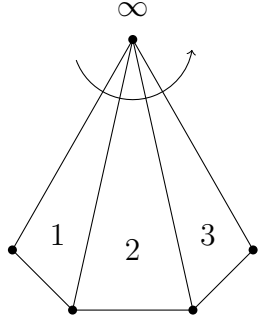
\caption{Ordering of infinite polygons in $\RR^2$. }
\enf

So if we put 
\[
V_\infty \,\,=\,\, \bigoplus_{(Q', A') \text{ infinite}} 
 V_{A'}, 
\]
then the formulas for the chain differential in the $C_\bullet(\Sigma(A'))$ define an algebra differential
$d$ in $T^\bullet(V_\infty)$ with $d^2=0$. This means that $\bR_\infty= V_\infty^*[-1]$ is an $A_\infty$-algebra.
 
Similarly to \S \ref{sec:1dim}, $\bR_\infty=\bgen_\infty$ as a graded vector space. 
Further, the only $A_\infty$-operations in $R_\infty$ are the differential and the binary bracket, so $R_\infty$
is an associative dg-algebra.  

As before, wee define the geometric subalgebra $R_\infty\subset \bR_\infty$ spanned by geometric infinite
marked polygons $(Q', A\cap Q')$. Thus $R_\infty=\gen_\infty$ as a vector space and the differential on $R_\infty$
is trivial, so $R_\infty$ is a graded associative algebra.

 The product of two infinite marked subpolygons $a,b$ is zero unless their union is again
an infinite marked subpolygon (in particular, the union is convex), with $a$ being on the left of $b$
and meeting $b$ along an edge, see Fig. \ref{fig:mult-r-infty}.
\bef
  \btp[scale=0.4, baseline=(current  bounding  box.center)]
\node (1) at (0,0){};
\node (2) at (3,0){};
\node (3) at (5,1){};
\node (4) at (2,6){};
\node (5) at (-1,2){};
\fill (1) circle (0.15); 
\fill (2) circle (0.15); 
\fill (3) circle (0.15); 
\fill (4) circle (0.15); 
\fill (5) circle (0.15); 
 
 \draw (0,0) -- (3,0) -- (5,1) -- (2,6) -- (-1,2) -- (0,0);
 \draw (0,0) -- (2,6); 
 
 \node at (0,2){$a$};
 \node at (2,2){$b$}; 
 \node at (2,7) {$\infty$}; 
\etp 
\hskip 1cm 
\btp[scale=0.4, baseline=(current  bounding  box.center)]
\draw[->] (0,3) --(3,3); 
\node at (1.5, 4){mult.};
\etp \hskip 1cm 
 \btp[scale=0.4, baseline=(current  bounding  box.center)]
\node (1) at (0,0){};
\node (2) at (3,0){};
\node (3) at (5,1){};
\node (4) at (2,6){};
\node (5) at (-1,2){};
\fill (1) circle (0.15); 
\fill (2) circle (0.15); 
\fill (3) circle (0.15); 
\fill (4) circle (0.15); 
\fill (5) circle (0.15); 
 
 \draw (0,0) -- (3,0) -- (5,1) -- (2,6) -- (-1,2) -- (0,0);
 \node at (2,2){$ab$};
  \node at (2,7) {$\infty$}; 

\etp
\caption{ Multiplication on $R_\infty$. }

\label{fig:mult-r-infty}
\enf
Now, looking at all coarse subdivisions of all marked subpylgons of $(\wt Q, \wt A)$, both finite and infinite,
we get an algebra differential in
\[
S^\bullet(V) \otimes T^\bullet(V_\infty), \quad 
V\,\,=\,\, \bigoplus_{(Q', A') \text{ finite}} 
 V_{A'}. 
\]
This differential preserves $S^\bullet(V)$ and gives there the differential defining $\gen_\infty$.
So we obtain, by Proposition \ref{prop:AL-mor} and, more particularly, Example \ref{ex:case-ass},
 $L_\infty$-morphisms
\be
\overset \bullet \phi: \bgen_\infty \lra R\Der_\Ass(\bR_\infty), \quad  \phi: \gen_\infty \lra R\Der_\Ass(R_\infty),
\ee
of which the second one will be more important for us. 
This is the setting of Gaiotto-Moore-Witten 
\cite{GMW} in the ``case of trivial coefficients". In the next sections we explain how can one
generalize our setting to  fully include that of \cite{GMW}. 

\vfill\eject

\section{Introducing coefficients: factorizing sheaves on secondary polytopes.}
\label{sec:coeff-arbitr}

Since we assumed that $A\subset\RR^d$ is in general position,
 all  proper faces of $Q$ are simplices. Each subset of $d$ points $\{\omega_1, \cdots, \omega_d\}\subset A$
gives a $(d-1)$-simplex $\sigma = \Conv\{\omega_1, \cdots, \omega_d\}$.

By an {\em oriented simplex} we will mean a simplex together with a numbering of vertices
 considered up to an even permutation. We will use the notation $\sigma
 = \< \omega_1, \cdots, \omega_d\>$  to denote the oriented $(d-1)$-simplex  given by the
 numbering of the verices $\omega_1, \cdots, \omega_d$. We write $\ol\sigma$ for ``$\sigma$
 with the opposite orientation", i.e, $\ol\sigma = \< \omega_{s(1)}, \cdots, \omega_{s(d)}\>$,
 where $s$ is any odd permutation.   
 
 \begin{Defi}\label{def:sys-coef}
  A {\em system of coefficients} for $A$ is a datum $N=(N_\sigma)$ of cochain complexes,
 one for each oriented $(d-1)$-simplex $\sigma$ with vertices in $A$, so that $N_{\ol\sigma} = N_\sigma^*$. 
 \end{Defi}
 
 We fix an orientation of $\RR^d$. Then for each marked subpolytope $(Q', A')\subset (Q,A)$,
 each codimension 1 face of $Q'$ is naturally an oriented $(d-1)$-simplex. We denote
 \be
 N_{A'} \,\,=\,\,\bigotimes_{\sigma\subset\partial Q'} N_\sigma.
 \ee
 Let $\Pc=\{(Q_\nu, A_\nu)\}$ be a regular polyhedral subdivision of $(Q,A)$. We define
 \be
 N_\Pc \,\,=\,\, \bigotimes_\nu N_{A'_\nu}.
 \ee
 If $\Pc'$ is a refinement of $\Pc$, then we have morphisms of cochain complexes
 \[
 \gamma_{\Pc' \Pc}: N_{\Pc'}\to N_{\Pc}
 \]
obtained by taking traces $N_{\sigma}\otimes N_{\ol\sigma} \to\k$ over those $(d-1)$-faces
of $\Pc'$ which are {\em internal}, i.e., are not faces of $\Pc$ and thus belong to exactly
two polytopes from $\Pc'$. 

\begin{prop}\label{prop:gamma}
 The maps $\gamma_{\Pc'\Pc}$ are transitive for chains of refinements: if $\Pc''$ refines $\Pc'$ and
$\Pc'$ refines $\Pc$, then
\[
\gamma_{\Pc''\Pc} = \gamma_{\Pc'\Pc}\circ\gamma_{\Pc''\Pc'}. \qed
\]
\end{prop}

We now interpret this proposition in terms of constructible complexes of sheaves on secondary polytopes.
Let $X$ be a CW-complex with a cellular stratification $(X_\alpha)_{\alpha\in I}$. In other words, $(I,\leq)$ is a poset
together with a strictly monotone function $d: I\to\ZZ_+$, 
each $X_\alpha$ is a locally closed subspace of $X$ homeomorphic to the Euclidean space $\RR^{d(\alpha)}$,
and the condition $\alpha\leq\beta$ is equivalent to the $X_\alpha$ being contained in the closure of $X_\beta$.
 A sheaf $\Fc$ of $\k$-vector spaces on $X$   constant along each $X_\alpha$ is, as well known
  \cite{curry}\cite {gelfand-macpherson} 
 described by the following data:
 \begin{enumerate}
 \item[(1)] Vector spaces $\Fc_\alpha = H^0(X_\alpha, \Fc)$ (stalk  of $\Fc$).
 
 \item[(2)] {\em Generalization maps} $\gamma_{\alpha\beta}: \Fc_\alpha \to\Fc_\beta$ defined for any $\alpha\leq\beta$
 and satisfying the transitivity conditions, i.e., forming a representation of the poset $I$ in $\k$-vector spaces. 
 
 \end{enumerate}

 Proposition \ref{prop:gamma}  implies that the data $(N_\Pc, \gamma_{\Pc',\Pc})$ give rise to a complex of sheaves $\Nc_A$
 on the secondary polytope  $\Sigma(A)$ which is constant on each open face
$F^\circ_\Pc$ with stalk $N_\Pc$ and with  the generalization maps 
 being the  $\gamma_{\Pc'\Pc}$. We then have the  following counterpart of
 Proposition \ref {prop:fact-sigma}.

 \begin{prop} The complex $\Nc_A$ satisfies the {\em factorization  property}:
the restriction of $\Nc_A$ to any closed face $F_\Pc= \prod_{(Q'_\nu, A'_\nu)\in\Pc} \Sigma(A'_i)$
is identified with the exterior tensor product $\boxtimes_\nu \Nc_{A'_\nu}$. \qed
\end{prop}

The fact that the generalization maps for a cellular sheaf go from the stalks at smaller cells to the stalks at bigger
cells, is a part of the general phenomenon that  sheaves have {\em cohomology} rather than {\em homology}. 
In particular, the cellular cochain complex of a sheaf $\Fc$ given by $(\Fc_\alpha, \gamma_{\alpha\beta})$
has the form
\[
C^\bullet_{\on{cell}}(X, \Fc) \,\,=\,\,\biggl\{ \bigoplus_{\dim\, X_\alpha =0} \orr(X_\alpha) \otimes\Fc_\alpha \lra
\bigoplus_{\dim\, X_\alpha =1} \orr(X_\alpha) \otimes\Fc_\alpha \lra\cdots\biggr\}
\]
where $\orr(X_\alpha) = H^{d(\alpha)}_c(X_\alpha, \k)$ is the 1-dimensional orientation space of the cell
$X_\alpha$ and the differentials are induced by the $\gamma_{\alpha\beta}$. 

So sticking to the cohomological picture, we form the graded vector space
\[
V_N^{{\coh}} \,\,=\,\, \bigoplus_{\stackrel{A'\subset A} {|A'|\geq d+1}} V_{A'}^\coh, \quad 
V_{A'}^\coh=
N_{A'}\otimes \orr(\Sigma(A'))[-\dim\Sigma(A')]. 
\]
Then the cellular {\em cochain} complexes of the factorizing sheaves $\Nc_{A'}$ unite to give a {\em coalgebra
differential} $d$, $d^2=0$,  in the symmetric coalgebra $S^\bullet(V^\coh_N)$. This is an alternative way of saying that:

\begin{prop}  (a) The cellular cochain differentials of the $\Nc_{A'}$ give rise to
an $L_\infty$-algebra structure on
\[
\bgen_N= \bgen_{A,N} \,\,=\,\, \bigoplus_{ (Q', A')\subset (Q,A)} E_{A'}, \quad
E_{A'}=
N_{A'}\otimes \orr(\Sigma(A'))[-\dim\Sigma(A')-1]. 
\]
The $L_\infty$-algebra $\bgen_N$ is nilpotent.

\vskip .1cm

(b) The subspace
\[
\gen_N=\gen_{A,N} \,\,=\,\,\bigoplus_{Q'\subset (Q,A)} E_{A\cap Q'}
\]
is an $L_\infty$-subalgebra in $\bgen_N$. 
\qed
\end{prop} 

As before $\gen_N$ will be more important for us than $\bgen_N$.

\vskip .2cm

Thus, the component $\gen_N^1$ (where Maurer-Cartan elements live) is
\[
\gen_N^1 \,\,=\,\, \bigoplus_{ Q'\subset (Q,A) } N_{A'}^{-\dim\Sigma(A')} \otimes \orr(\Sigma(A')), \quad A'=A\cap Q'. 
\]
In other words, an element of $\gen_N^1$ is a rule associating to each marked $d$-simplex $(Q',A')\subset (Q,A)$ an
element of $N_{A'}^0$, to each circuit $(Q',A')$ an element of $N_{A'}^{-1}$ and so on.

 \vskip .2cm
 
 \noindent {\bf Relative case.} Suppose $A$ is extended to $\wt A=A\cup \{\infty\}$, as in \S \ref{sec:rel-set}. 
 Let $\wt N$ be a coefficient system for $\wt A$.  We then generalize Propositions 
 \ref{prop:gentilde} 
 and Corollary \ref{cor:L-infty-mor} 
 to: 
 \begin{prop}
 (a) We have a semi-direct product decomposition
 \[
 \gen_{\wt N} = \gen_N \ltimes \gen_{\wt N, \infty}
 \]
 where $\gen_{N,\infty}$ is the direct sum of summands $E_{A'}$
  corresponding to infinite marked subpolytopes $(Q',A')$.
 Therefore we have an $L_\infty$-morphism $\alpha: \gen_{\wt N}\to R\Der_\Lie (\gen_{\wt N, \infty})$. 
 
 (b) If $d=2$, then $\gen_{\wt N,\infty}$ lifts to an $A_\infty$-algebra $A_{\wt N, \infty}$ and we have an $L_\infty$-morphism
 $\psi: \gen_{\wt N}\to R\Der_\Ass(R_{\wt N, \infty})$. \qed
  \end{prop}
  
  \vfill\eject

  \section{Coefficients in bimodules: $d=2$}\label{sec:bim}
  
The systems of coefficients considered in \S \ref{sec:coeff-arbitr} consisted in data associated
to codimension 1 simplices only. It would be   natural to consider more general
coefficient systems which would associate some data to simplices of all dimensions 
from $0$ up to $(d-1)$. In this section we explain how to do this in the case $d=2$.

Let us denote elements of $A$ as $i,j,k$ and so on. 

\begin{Defi}\label{def:ex-sys-co} 
An {\em extended system of coefficients} for $A$ is a following system of data:
\begin{enumerate}
\item[(1)]  For each $i\in A$,  an
associative dg-algebra $S_i$.

\item[(2)] For each oriented edge (1-simplex)  $(i,j)$, $i,j\in A$, a (differential graded)
 $(S_i, S_j)$-bimodule
$N_{ij}$ which we assume projective, of finite rank, 
 over  the graded algebra underlying $S_i\otimes_\k S^\op_j$.

\item[(3)] For each $(i,j)$, a  pairing 
\[
\beta_{ij}: N_{ij}\otimes_\k  N_{ji} \lra S_i\otimes_\k S_j, 
\]
which is a morphism of $(S_i\otimes_\k S_j, S_i\otimes_\k S_j)$-bimodules. It is required that:

\item[(4)] The pairing $\beta_{ij}$ is non-degenerate, i.e., the morphism of $(S_i, S_j)$-bimodules
\[
\beta^t_{ij}: N_{ij} \to \Hom_{S_j\otimes_k S_i^\op} (N_{ji}, S_j\otimes_k S_i^\op)
\]
induced by $\beta_{ij}$, is an isomorphism.

\item[(5)]  The diagram
\[
\xymatrix{
N_{ij}\otimes_\k N_{ji} \ar[r]^{\beta_{ij}}
\ar[d]_{\on{perm.}}
 & S_i\otimes_\k S_j 
 \ar[d]^{\on{perm.}}
\\
N_{ji}\otimes_\k N_{ij} \ar[r]^{\beta_{ji}} & S_j\otimes_\k S_i
}
\]
is commutative. 
\end{enumerate}
\end{Defi}
Thus the setting of Definition \ref{def:sys-coef}, specialized to $d=2$, corresponds to the case when
all $S_i=\k$. 

\begin{rem}\label{rem:Si-Fij}
The data of an extended  system of coefficients can be interpreted as follows. The datum (1) means that we have 
 a  (pre-)triangulated dg-category $\Cc_i$. More precisely, let us assume that $S_i$ is smooth and proper
 \cite{toen-vaquie}, in particular, that $\dim_\k S_i <\infty$. Then we take 
   $\Cc_i = D\Mod^{<\infty}_{S_i}$  to be the derived category  formed by
finite-dimensional 
right dg-modules over $S_i$. Note that $\Cc_i$ has a Serre functor \cite{BK}, i.e., a covariant exact functor
$\Sen_i: \Cc_i\to\Cc_i$ equipped with natural isomorphisms
\[
\Hom_{\Cc_i} (M, M')^* \simeq \Hom_{\Cc_i}(M', \Sen_i(M)).
\]
Here and elsewhere $*$ means dualization over $\k$. Explicitly, $\Sen_i$ is given by the Nakayama formula
\be\label{eq:nakayama}
\Sen_i(M) =  \bigl(R\Hom_{S_i}(M, S_i)\bigr)^*.
\ee
The datum (2) means that we have an exact functor $F_{ij}: \Cc_i\to \Cc_j$ (given by the tensor product with $N_{ij}$). 
Finally,  the datum (3) means that $F_{ij}$ is left
adjoint to $\Sen_i F_{ji}$ (or, equivalently, that $F_{ji}$ is left adjoint to $\Sen_j F_{ij}$). This follows from
\eqref{eq:nakayama}. 

Note that for a pair of exact functors 
$F: \Cc\to\Dc$, $G: \Dc\to\Cc$ between triangulated categories $\Cc, \Dc$ with Serre functors
$\Sen_\Cc, \Sen_\Dc$ the condition  ``$F$ is left adjoint to $\Sen_\Cc G$" is symmetric in $F$ and $G$
(unlike the usual condition of adjointness). Indeed, such ``twisted adjointness" is expressed via natural
isomorphisms
\[
\Hom_\Dc(F(X), Y) ^* \simeq \Hom_\Cc (G(Y), X), \quad X\in\Cc, Y\in\Dc.
\]
 \end{rem}
 
 \vskip .2cm

We fix an orientation of $\RR^2$. 
Suppose given  an extended system of coefficients $(S_i, N_{ij})$.  Let 
$Q'\subset(Q,A)$ be a  subpolygon. Let us label the vertices of $Q'$ cyclically, in the counterclockwise
order as $i_0, \cdots, i_m$. We can then form the {\em cyclic tensor product}
\[
N_{Q'} \,\,=\,\, \biggl( N_{i_0i_1}\otimes_{S_{i_1}} N_{i_1i_2}\otimes_{S_{i_2}} \cdots
\otimes_{S_{i_{m-1}}} N_{i_{m-1}i_m} \biggr) \otimes_{S_{i_0^{\on{op}}}\otimes S_{i_m}} N_{i_m i_0},
\]
which can be pictorially represented in Figure \ref{fig:cyc-prod}, which makes evident its cyclic symmetry
(independence on the choice of the starting point $i_0$). Note that $N_{Q'}$ is just a cochain complex:
all module structures have been used in the formation of the tensor products. 
\usetikzlibrary{decorations.text}

\bef
\begin{tikzpicture}
\draw [ double distance=10mm,
        rotate=120,
       postaction={
            decorate,
            decoration={
                  raise=-1ex,
                  text along path, 
                  reverse path,
                  text align={fit to path stretching spaces},
                  text={ 
                  |\large|
                  {$N_{i_0i_1}$}
                   {$\otimes_{S_{i_1}}$} 
                    {$N_{i_1i_2}$} 
                    {$\otimes_{S_{i_2}}$} 
                  {$N_{i_2i_3}$}
                   {$\otimes_{S_{i_3}}$} 
                   {$\cdots \otimes$} 
                   {$\cdots$}
                  {$\otimes_{S_{i_{n-1}}}$} 
                  {$N_{i_{n-1}i_n}$} 
                  {$\otimes_{S_{i_n}} $}
                  {$N_{i_n i_0}$}
                  {$\otimes_{S_{i_0}}$}
                   }
            }
       }
    ] (0,0) circle (2.7cm);
\end{tikzpicture}
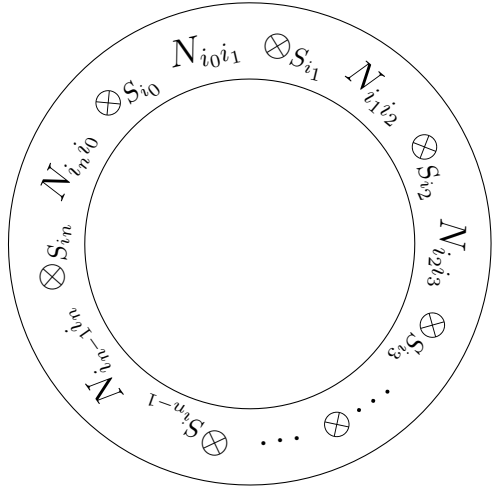
\caption{ The cyclic tensor product.} 
\label{fig:cyc-prod}
\enf

 \begin{rem}\label{rem:cyc-prod}
 More generally, we can define a cyclic tensor product $N_\sigma$ for any closed
 oriented edge path $\sigma = (i_0, i_1, \cdots, i_n, i_0)$, not nevessarily coming from
 the boundary of a convex polygon. As before,  $N_\sigma$ depends only   on $\sigma$ as a closed path
 and not on the choice of a starting point.  
 \end{rem}
 
 Let
 \[
 \sigma_1 = (\cdots, h,i,j,k, \cdots), \,\,\, \sigma_2 = ( \cdots, p,j,i,q,\cdots)
 \]
 be two closed oriented edge paths having one edge $[i,j]$ in common, with opposite orientations. Then we
 can concatenate $\sigma_1$ and $\sigma_2$ along this edge, erasing it and obtaining a new closed
 oriented edge path
 \[
 \sigma_1 *_{[i,j]}\sigma_2 = (\cdots, p,j,k, \cdots, h,i,q,\cdots),
 \]
 see Fig. \ref{fig:concat}. 
 
 \bef
\begin{tikzpicture}[scale=.3]
\node (i) at (0,0){};
\fill (i) circle (0.1); 

\node (h) at (-7,1){};
\fill(h) circle (0.15);

\node (q) at (3,-5){};
\fill(q) circle (0.15); 

\node (j) at (6,7){};
\fill(j) circle (0.1); 

\node (p) at (13,9){};
\fill(p) circle (0.15); 

\node (k) at (3,11){};
\fill(k) circle (0.15); 

\draw[->>, line width=0.8] (-7,1) -- (0,0); 
\draw[->>, line width=0.8] (0,0) -- (3,-5); 
\draw[->>, line width=0.8] (6,7) -- (3,11); 
\draw[->>, line width=0.8] (13,9) -- (6,7); 
\draw[<<->>, dotted, line width=1] (0,0) -- (6,7); 

\node at (-3,6){\huge$\sigma_1$}; 

\node at (9,1){\huge$\sigma_2$}; 

\node at (-7,0){\large$h$};

\node at (-1,-1){\large$i$}; 

\node at (6.5,8){\large$j$}; 

\node at (13,10){\large$p$}; 

\node at (3,12){\large$k$}; 

\node at (3,-6){\large$q$};

\draw[dashed] (-9,4) -- (-7,1); 
\draw[dashed] (0,12) -- (3,11); 
\draw[dashed] (13,9) -- (16,7); 
\draw[dashed] (3,-5) -- (7,-7);

\end{tikzpicture}
\caption{ Concatenation.} 
\label{fig:concat}
\enf

\noindent In this situation, the pairing $\beta_{ij}$ gives rise to the {\em concatenation map}
\be\label{eq:comp-map}
\gamma_{\sigma_1, \sigma_2}^{[i,j]}: N_{\sigma_1}\otimes_\k N_{\sigma_2} \lra N_{ \sigma_1 *_{[i,j]}\sigma_2 }.
\ee
To define it, consider two decomposable tensors
\[
\begin{gathered}
n_{\sigma_1} = \,\,\cdots \otimes n_{hi}\otimes n_{ij}\otimes n_{jk} \otimes \cdots \,\, \, \in \,\, N_{\sigma_1},\\
n_{\sigma_2} = \,\, \cdots  \otimes n_{pj}\otimes n_{ji} \otimes n_{iq}\otimes \cdots \,\,\, \in \,\, N_{\sigma_2},
\end{gathered} 
\]
where $n_{ij}\in N_{ij}$ etc.  Suppose
\[
\beta_{ij}(n_{ij}\otimes n_{ji}) \,\,=\,\,\sum_\nu s'_\nu \otimes s''_\nu, \quad s'_\nu\in S_i, \,\, s''_\nu\in S_j. 
\]
Then we define
\[
\gamma_{\sigma_1, \sigma_2}^{[i,j]} (n_{\sigma_1}\otimes n_{\sigma_2}) 
\,\,=\,\,  \sum_\nu \bigl( \,\, \cdots \otimes n_{pj}s''_\nu \otimes n_{jk} \otimes \cdots \otimes n_{hi}s'_\nu \otimes n_{iq}\otimes \cdots\bigr).
\]
It is immediate to check that this gives a well defined map of cyclic tensor products. 
 
 \vskip .2cm

If  now $\Pc$ is a polygonal subdivision of $Q'=\Conv(A')$ into $\{Q''_\nu\}$, we define
\[
N_\Pc \,\,=\,\,\bigotimes_\nu N_{Q''_\nu}
\]
  (tensor product over $\k$). Note that the concatenation maps \eqref{eq:comp-map} corresponding to various intermediate edges of
  $\Pc$, commute with each other, so the result of applying them all in any order, is independent on the
  order.  We call  this result the {\em composition map} and denote it
  \be\label{eq:comp-map-2}
  \gamma_{\Pc}:  N_\Pc\lra  N_{A'}. 
  \ee
    Further, if $\Pc'$ is a subdivision refining $\Pc$,
  the  composition maps \eqref{eq:comp-map-2} corresponding to all polygons of $\Pc$, give rise to a map
  \be
  \gamma_{\Pc'\Pc}: N_{\Pc'}\lra N_{\Pc}
  \ee
  which we call the {\em generalization  map} to emphasize that we want to use it to construct cellular sheaves on secondary
  polytopes.  More precisely, we have the following statement whose proof is straightforward.
  
\begin{prop}
The data $(N_\Pc, \gamma_{\Pc'\Pc})$ define a factorizing system  of complexes of sheaves $\Nc_{A'}$ on the $\Sigma(A')$,
with the stalk of $\Nc_{A'}$ on the face $F_{\Pc}\subset\Sigma(A')$ being $N_\Pc$.
 Denote
\[
E_{A'} \,\,=\,\, N_{A'}\otimes\orr(\Sigma(A'))[-\dim\Sigma(A')-1].
\]
Then the differentials in the cellular cochain complexes of the $\Nc_{A'}$
unite to make the graded vector space $\bgen=\bigoplus_{(Q', A')\subset (Q,A)} E_{A'}$ into a nilpotent $L_\infty$-algebra.
The subspace $\gen = \bigoplus_{Q'\subset (Q,A)} E_{A\cap Q'}$ is an $L_\infty$-subalgebra in $\bgen$. 
 \qed
\end{prop}

\noindent {\bf Relative situation (as in \S \ref {sec:rel-set}). } $\wt A=A\cup\{\infty\}$, $\wt Q=\Conv(\wt A)$.
Let $(S_i, N_{ij})$ be an extended coefficient system for $A$ (not $\wt A$)
For any infinite marked polytope $(Q', A')$ with vertices $\infty, i_0, \cdots, i_m$ in counterclockwise order,
we define the {\em linear tensor product}
\[
L_{A'} \,\,=\,\, N_{i_0i_1}\otimes_{S_{i_1}} N_{i_1i_2} \otimes_{S_{i_2}}\cdots
\otimes _{S_{i_{m-1}}} N_{i_{m-1} i_m}. 
\]
This is an $(S_{i_0}, S_{i_m})$-bimodule. If $\Pc$ is a subdivision of a (possibly infinite)
subpolygon $(Q', A')$ into $(Q''_\nu, A''_\nu)$, then we define
\[
N_{\Pc} \,\,=\,\, \bigotimes_{A''_\nu \not\ni\infty}N_{A''_\nu} \,\,\,\otimes \,\,\,\bigotimes_{A''_\nu\ni\infty} L_{A''_\nu}. 
\]
and these complexes, again, give a factorizable system of  complexes of sheaves on the $\Sigma(A')$.
Restricting to the geometric summands as before, we obtain:
\begin{prop}
Denoting for   $A'\ni\infty$,
\[
F_{A'} \,\,=\,\, L_{A'}\otimes\orr(\Sigma(A'))[-\dim\Sigma(A')-1],
\]
we get the following data:
\begin{enumerate}
\item[(1)] An $L_\infty$-algebra $\gen = \bigoplus_{Q'\not\ni\infty} E_{A\cap Q'}$.

\item[(2)] An $A_\infty$-algebra $R_\infty=\bigoplus_{Q'\ni\infty} F_{A\cap Q'}$.

\item[(3)] An $L_\infty$-morphism $\psi: \gen\to R\Der(R_\infty)$. 

\end{enumerate} \qed

\end{prop}

\begin{rem}
Study   of Fukaya-Seidel categories   (see
\S \ref{sec:MCFS} below for more background) leads naturally to a structure somewhat more general 
and symmetric
  than an extended system of coefficients. In that structure,
   for each point $i\in A$ we have not just a category, but a 
   {\em local system of triangulated categories} (together with pre-triangulated dg-enhancements)
    ${\mathcal C}_i$ over a small circle $S_i^1$ centered at $i$. It is required that:
    
    \begin{enumerate}

\item[(a)]  Each fiber category ${\mathcal C}_{i, \theta}, \theta\in S^1_i$ has  Serre functor $\Sen_i$.

\item[(b)]  The monodromy around the circle $S^1_i$ gives rise to the functor isomorphic to $\Sen_i\circ [d]$, where $[d]$ denotes the functor of shift of the grading by $d\in {\bf Z}$ (for the Fukaya-Seidel category associated with the pair $(X,W)$ the number $d$ is the dimension of the complex manifold $X$).

\end{enumerate}

Further, for each pair of different points $i, j\in A$ we have functors $F_{ij}: {\mathcal C}_{j,{\theta_j}}\to {\mathcal C}_{i,\theta _i}$ and 
$F_{ji}: {\mathcal C}_{i,{\theta_i}}\to {\mathcal C}_{j,\theta _j}$
between the fiber categories over the points $\theta_i\in S^1_i$ and $\theta_j\in S^1_j$ which are intersection points of the straight segment $ij$ with circles $S^1_i$ and $S^1_j$. These functors are required to be ``twisted adjoint" to each
other  in the sense of the Remark \ref{rem:Si-Fij}.

\vskip .2cm

A choice of a  ``generic" point at infinity (equivalently, of a generic half-plane) breaks the local rotational symmetry,
 since each circle $S^1_i$ contains the intersection point with the ray from $i$ to infinity. This allows us to replace each 
 $S^1_i$ by a segment. After this,  the more general structure described above  reduces effectively to 
 an extended system of coefficients.

\end{rem}


\vfill\eject

\section{Analysis of the algebra $R_\infty$ and the morphism $\psi$ for $d=2$.}\label{sec:anal}

Since we assumed that $\wt A=A\cup\{\infty\}$ is in   general position,  $A$ has a total order by
the slope of the line $(i\infty)$ from left to right. We get a decomposition
\[
R_\infty \,\,=\,\,\bigoplus_{\stackrel {i,j\in A}{i<j}} R_{ij}, 
\]
where $R_{ij}$ is the direct sum of the summands $F_{A\cap Q'}$ corresponding to infinite polygons $Q'$ whose edges 
through $\infty$ are $(i\infty)$ and $(j\infty)$. Note that $R_{ij}$ is an $(S_i, S_j)$-bimodule. 

\begin{prop}\label{prop:hi-comp}
(a) The higher compositions $m_n, n\geq 3$, in $R_\infty$, are trivial. The only non-trivial part of the binary compositions
$m_2$ are the maps
\[
\mu_{ijk}: R_{ij}\otimes R_{jk}\to R_{ik}, \quad i<j<k.
\]
Thus $R_\infty$ is a strictly upper-triangular dg-algebra (without unit). 

\vskip .1cm

(b) The multiplication $\mu_{ij}$ is $S_j$-bilinear with respect to the right $S_j$-module structure
on $R_{ij}$ and the left $S_j$-module structure on $R_{jk}$. Therefore, putting $R_{ii}=S_i$,
we obtain a non-strictly triangular associative dg-algebra with unit
\[
R=\bigoplus_{i\leq j} R_{ij}. 
\]
\end{prop}

\begin{rem}
Note that part (b) means that $\Mod_R$ has a semi-orthogonal decomposition  with quotients
\[
\Mod_{S_{i_1}}, \Mod_{S_{i_2}}, \cdots, \Mod_{S_{i_r}},
\]
where $i_1 < i_2 <\cdots i_r$ are all the elements of $A$ in our order. 

\end{rem}

\noindent {\sl Proof of  Proposition \ref{prop:hi-comp}:} (a) The only coarse subdivisions of an infinite polygon into 
all infinite subpolygons contain either 2 parts (and so account for $m_2$), or 1 part (and so
account for $m_1=d$, as in Examples \ref{ex-faces-sigma} (c2) and (c1). The subdivisions into 1 part are excluded
in the geometric subalgebra $R_\infty\subset\bR_\infty$. 
This shows part (a).
Part (b) is clear. \qed

\vskip .2cm

We now describe the morphism 
\[
\psi: \gen\lra R\Der(R_\infty) = C^{\geq 1}(R_\infty, R_\infty)[1]
\]
 from Proposition 10.5
explicitly. Because of the direct sum decompositions
\[
\gen = \bigoplus_{A'\subset A \text{ geom.}} E_{A'}, \quad R_\infty = \bigoplus_{B\ni\infty \text{ geom.}} F_B, \quad |A'|, |B|\geq 3, 
\]
$\psi$ is given by its matrix elements with respect to the induced direct sum decompositions of the source and target.
These matrix elements are morphisms of dg-vector spaces
\be\label{eq:mat-el-psi}
\psi_{A'}^{(B_1, \cdots, B_m|C)}: E_{A'} \lra\Hom_\k (F_{B_1}\otimes_\k \cdots\otimes_\k F_{B_m}, F_C),
\ee
given for all finite geometric $A'\subset A$ and all infinite geometric $B_1,\cdots, B_m, C\subset \wt A$ of cardinality $\geq 3$.
By the analysis of  Example \ref {ex:case-ass}, $\psi_{A'}^{(B_1, \cdots, B_m|C)}\neq 0$ only if $Q'=\Conv(A')$ together with the
$P_\nu = \Conv(B_\nu)$, forms a coarse regular subdivision of $P=\Conv(C)$. 
In other words, nonzero matrix elements of $\psi$ correspond to coarse regular subdivisions of infinite subpolygons in $(\wt Q, \wt A)$
which contain exactly one finite polygon (the rest being infinite). Let us call such subdivisions {\em 1-finite}. We will now
describe all 1-finite subdivisions explicitly, starting with the following class of examples.

\begin{ex}\label{ex:1-finite-1}
Let $Q'$ be a finite  subpolygon of $(Q,A)$,
and let 
\[
A'=A\cap Q',\,\, \,\wt Q'=\Conv(A'\cup\{\infty\}), \,\,\,
\wt A'=\wt A\cap\wt Q'.
\]
Since we assumed that no 3 points of $\wt A$ lie on a line, the boundary $\partial Q'$ splits into two
distinct parts meeting at two points:
\begin{itemize}

\item The {\em positive boundary}, i.e., the union of sides facing $\infty$, drawn thick in Figure 
\ref {fig:pos-neg-bou} 
and denoted $\partial_+Q'$.

\item The {\em negative boundary}, denoted $\partial_- Q'$ formed by the sides which lie on the
opposite side of $\infty$. 

\end{itemize}

\bef
\btp[scale=0.4]
\node (1) at (0,0){};
\node (2) at (6,4){};
\node (3) at (2,6){};
\node (4) at (-1,6){};
\node (5) at (-6,4){};
\node (6) at (0,15){};

\fill(1) circle (0.15); 
\fill(2) circle (0.15); 
\fill(3) circle (0.15); 
\fill(4) circle (0.15); 
\fill(5) circle (0.15); 
\fill(6) circle (0.15); 

\draw (0,0) -- (6,4) -- (2,6) -- (-1,6) -- (-6,4) -- (0,0); 
\draw (-6,4) -- (0,15) -- (-1,6); 
\draw (2,6) -- (0,15) -- (6,4); 

\draw[line width=0.8mm] (-6,4) -- (-1,6) -- (2,6) -- (6,4); 

\node at (0,16){$\infty$};
\node at (0,3){$Q'$};
\node at (-3,6){$\eta_1$};
\node at (0.5,6.7){$\eta_2 \cdots$};
\node at (4,6){$\eta_m$};
\node at (-2,9){$\Pi_1$};
\node at (0.5,9){$\Pi_2 ...$};
\node at (2.4,9){$\Pi_m$};

\node (d) at (10,7){\large $\partial_+ Q'$};
\node (e) at (6,6){};
\draw[->, line width=0.5mm] (d) -- (e); 

\node (f) at (8,1){\large $\partial_- Q'$};
\node (g) at (4,2){};
\draw[->] (f) -- (g); 

\etp
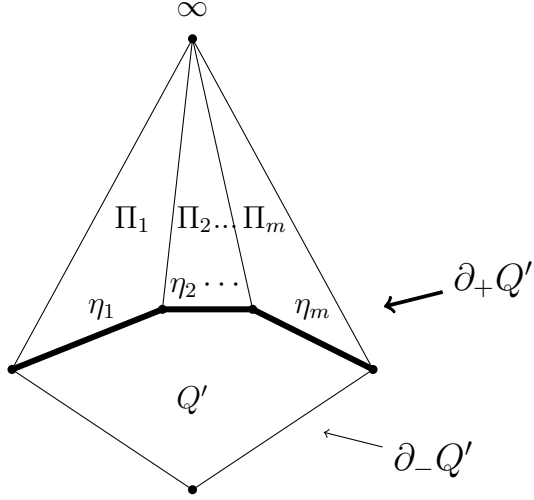
\caption{The positive and negative boundary.}
\label{fig:pos-neg-bou}
\enf
We denote the sides of $\partial_+ Q'$ by $\eta_1, \cdots, \eta_m$ in the natural order from left to right
(counter-clockwise rotation around $\infty$), and let $\Pi_1, \cdots, \Pi_m$ be the triangles with one vertex
at $\infty$ and the opposite sides being $\eta_1, \cdots, \eta_m$. Let also $D_\nu = \wt A\cap \Pi_\nu$.
Then $(\Pi_\nu, D_\nu)$ is an infinite marked subpolygon  of $(Q', A')$, and we have a polygonal
subdivision (see Fig. \ref{fig:pos-neg-bou})
\be\label{eq:coarse-sub-2}
(\wt Q', \wt A') \,\,=\,\, (Q', A') \cup (\Pi_1, D_1)\cup\cdots\cup (\Pi_m, D_m).
\ee
It is easy to see that \eqref {eq:coarse-sub-2} is a 1-finite subdivision. Indeed, the unique (up to scalar factors and  adding global
affine-linear functions) piecewise-affine function $f: \wt Q'\to\RR$ is determined by putting $f(\infty)=1$ and $f|_{A'}=0$.
\end{ex}

However, Example \ref {ex:1-finite-1} does not give the most general form of a 1-finite subdivision. To obtain the general
form, considering the following situation, extending the above one. 

\begin{ex}\label{ex:1-finite-2}
Let $P\subset (\wt Q, \wt A)$ be an infinite subpolygon  and $Q'\subset P$ be a
 finite subpolygon such that $\partial_-Q'$ is contained in $\partial P$
and is the union of edges from $\partial P$. The two remaining finite parts of $\partial P$ (not in $\partial_-Q'$)
will be called {\em handles}, see Fig. \ref{fig:gen-1-fin}. We denote the left handle by $\lambda$ and
edges constituting it by
$\lambda_1, \cdots, \lambda_p$. Similarly, we denote the right handle by $\rho$ and the edges constituting
 it by $\rho_1, \cdots, \rho_q$,
in our order from left to right. The handles may be empty, i.e., it may happen that $p=0$ or $q=0$ or both.

\bef
\btp[scale=0.4]
\node (1) at (0,0){};
\node (2) at (6,4){};
\node (3) at (2,6){};
\node (4) at (-1,6){};
\node (5) at (-6,4){};
\node (6) at (0,15){};

\fill(1) circle (0.15); 
\fill(2) circle (0.15); 
\fill(3) circle (0.15); 
\fill(4) circle (0.15); 
\fill(5) circle (0.15); 
\fill(6) circle (0.15); 

\node (7) at (-7,7){};
\fill (7) circle (0.15); 

\node (8) at (8,6){};
\fill (8) circle (0.15); 

\node (9) at (9, 9){};
\fill (9) circle (0.15); 

\draw (0,15) -- (9,9);
\draw[line width=0.8mm, color=gray] (9,9) -- (8,6) -- (6,4);

\draw (0,0) -- (6,4) -- (2,6) -- (-1,6) -- (-6,4) -- (0,0); 
\draw  (0,15) -- (-1,6); 
\draw (2,6) -- (0,15); 

\draw (-7,7) -- (0,15); 
\draw[line width=0.8mm, color=gray] (-7,7) -- (-6,4); 

\draw[line width=0.8mm] (-6,4) -- (-1,6) -- (2,6) -- (6,4); 

\node at (0,16){$\infty$};
\node at (0,3){$Q'$};
\node at (-3.5,6){$\eta_1$};
\node at (0.5,6.7){$\eta_2 \cdots$};
\node at (4,6){$\eta_m$};
\node at (-3,9){$P_1$};
\node at (0.5,9){$P_2 ...$};
\node at (5,9){$P_m$};

 \node (e) at (6,6){};

\node (f) at (8,1){\large $\partial_- Q'$};
\node (g) at (4,2){};
\draw[->] (f) -- (g); 

\node[text width=3cm, left] (LH) at (-10, 6){Left handle $\lambda=(\lambda_1, \cdots, \lambda_p)$
here $p=1$}; 

\node (lh) at (-7,5.5){};
\draw[->, line width=0.8mm, color=gray] (LH) -- (lh); 

\node[text width=3cm, right] (RH) at (12,7){Right handle $\rho= (\rho_1, \cdots, \rho_q)$ here $q=2$};
\node (rh) at (8.5,6){};
\draw[->, line width=0.8mm, color=gray] (RH) -- (rh); 

\node at (5,14){$P$};

\etp
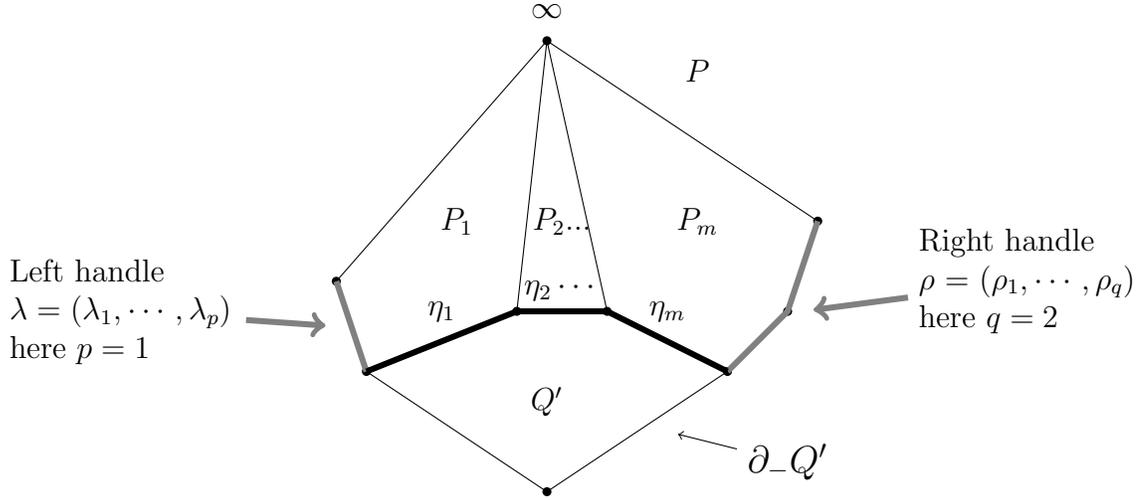
\caption{General 1-finite subdivision.}
\label{fig:gen-1-fin}
\enf
Note that by taking the convex hull $\wt Q'$ of $Q'$ and $\infty$, we get the situation described in Example
\ref {ex:1-finite-1}, from which we retain the notation. 
We then have a decomposition of $P$ into $Q'$ and polygons $P_1, \cdots P_m$, where:
\begin{itemize}
\item $P_\nu=\Pi_\nu$ for $\nu=2, \cdots, m-1$.

\item  The finite part of $\partial P_1$ is composed out of $\lambda$ and $\eta_1$.

\item The finite part of $\partial P_m$ is composed out of $\eta_m$ and $\rho$. 

\end{itemize}
\end{ex}

\begin {prop}
The  subdivision described in Example \ref {ex:1-finite-2}, is 1-finite, and each 1-finite subdivision is
obtained in this way. 
 \end{prop}

\noindent {\sl Proof:} By construction, the subdivision contains exactly one finite subpolygon. Let us show that it is
regular. For this, note that, as in Example \ref {ex:1-finite-1}, the conditions $f(\infty)=1$ and $f_{Q'}=0$
determine a unique continuous fonction $f: P\to\RR$ which is affine on each polygon and breaks exactly
along the subdivision. Further, this function is unique with these properties up to scalar factors and
adding global affine functions. This follows from analyzing the restriction to $\wt Q'=\Conv(Q\cup \{\infty\})$. 
Therefore our subdivision is coarse and therefore is 1-finite. 

Conversely, note that the construction of Example \ref {ex:1-finite-2}
by its very nature produces all
 subdivisions of any infinite polygon $P$  which contain exactly one
finite subpolygon.  \qed

\vskip .2cm

In the situation of Example \ref {ex:1-finite-2}, we put
\[
A'=A\cap Q', \,\,\, B_\nu = \wt A\cap P_\nu, \,\,\, C=\wt A\cap P. 
\]
In order to describe the corresponding matrix element $\psi_{A'}^{(B_1, \cdots, B_m|C)}$, we write it
in the transposed way, as a map
\be\label{eq:psi-trans}
\psi_{A'; B_1, \cdots, B_m}^C: E_{A'}\otimes_\k F_{B_1} \otimes_\k\cdots \otimes_\k F_{B_m}\lra F_C. 
\ee
Let us denote
\[
N_\lambda = N_{\lambda_1}\otimes_S \cdots \otimes_S N_{\lambda_p}, \quad N_\rho = N_{\rho_1}\otimes_S \cdots
\otimes_S N_{\rho_q}
\]
the linear tensor products corresponding to the handles $\lambda$ and $\rho$ considered as oriented edge paths.
Here and in the sequel we write $\otimes_S$ for the tensor product over the algebra $S_i$ corresponding to the
vertex common to the two edges in an edge path. 
Then
\[
\begin{gathered}
F_{B_1} = N_\lambda\otimes_S E_{\eta_1}, \quad F_{B_\nu}=E_{\eta_\nu}, \,\,\, 2\leq\nu\leq m-1, \quad
F_{B_m} = E_{\eta_m}\otimes_S N_\rho,\\
F_C = N_\lambda\otimes_S E_{A'}\otimes_S N_\rho. 
\end{gathered} 
\]

Let
 \be\label{eq:gamma-Q'}
 \gamma_{Q'}: E_{A'}\otimes_\k F_{D_1}\otimes_k\cdots \otimes_k F_{D_m} \lra F_{\wt A'}
 \ee
 be the composition map 
 obtained, similarly to \eqref{eq:comp-map-2}, 
 by substituting the pairings $\beta_{ij}$ over all the intermediate edges $[i,j]$
 in the coarse subdivision \eqref {eq:coarse-sub-2}. 
 Define
\be\label{eq:mat-el-psi-trans}
\psi_{A'; B_1, \cdots, B_m}^C \,\,=\,\, \Id_{N_\lambda} \otimes_S \gamma_{Q'}\otimes_S \Id_{N_\rho}. 
\ee
Then the source and target of $\psi_{A'; B_1, \cdots, B_m}^C$ are as specified in \eqref{eq:psi-trans}.

\begin{prop}\label{prop:comp-phi}
The maps $ \psi_{A'; B_1, \cdots, B_m}^C$ corresponding to all pairs $(P, Q')$ consisting of an
infinite subpolygon $P\subset (\wt Q, \wt A)$ and a finite subpolygon $Q'$ with $\partial_-Q'$ contained
in the finite part of $\partial P$, 
  are precisely all the
matrix elements of the $L_\infty$-morphism $\psi: \gen\to C^{\geq 1}(R_\infty, R_\infty)[1] $
with respect to the direct sum decompositions of the source and target  described above.
\end{prop} 

\noindent {\sl Proof:} This is a consequence of the analysis of Example \ref {ex:case-ass} applied to
our case.\qed 

\vfill\eject

\section {The universality theorem ($d=2$). } \label{sec:univ}

We work in the setting of \S \ref{sec:bim} and  \ref{sec:anal}. So we have the triangular dg-algebra
\[
R=\bigoplus_{i\leq j} R_{ij}, \quad R_{ii}=S_i. 
\]
We are interested in deformations of $R$ which, first,  preserve the triangular structure and, second,  do not change the quotients
of the semi-orthogonal decomposition. The second condition  means   that the algebras $S_i$
are not deformed. Such deformations are governed by the {\em ordered Hochschild complex}
 $\vC^\bullet(R, R)$, with
 \[
 \begin{gathered}
 \vC^n(R, R)\,\,=\\
 =\,\,\bigoplus_{i_0 < i_1 <\cdots < i_n} 
 \Hom_{S_{i_0}\otimes S_{i_n}^{\on{op}}}
 \biggl( R_{i_0i_1}\otimes_{S_{i_1}} R_{i_1 i_2}\otimes_{S_{i_2}} \cdots \otimes_{S_{i_{n-1}}} R_{i_{n-1}}, R_{i_n}, \,\,\,
 R_{i_0 i_n}\biggr)
 \end{gathered} 
 \]
 (strict inequalities under the direct sum sign). 
 
 \begin{rem}
Note that $\vC^\bullet(R, R)$ is a subcomplex in the ordinary Hochschild complex $C^\bullet(R, R)$,
specified by the multi-linearity conditions. Moreover, this subcomplex is closed with respect to
the Lie bracket and is thus a dg-Lie subalgebra. 

\end{rem}

\begin{thm}\label{thm:uni}
The $L_\infty$-moprhism $\psi$ factors through an $L_\infty$-morphism
\[
\Psi: \gen \lra \vC^{\geq 1}(R, R)[1]. 
\]
The morphism $\Psi$ is a quasi-isomorphism. 
\end{thm}

\begin{rems} (a) Note that the differential in $\gen$ is trivial, while that in $\vC^{\geq 1}(R, R)[1]$ is not.

 \vskip .2cm
 
 (b) Theorem \ref{thm:uni} suggests an alternative way to think about the $L_{\infty}$-algebra  $\gen$ of convex polygons. 
  Namely, we can start with the dg-algebra $R$ and define the $L_{\infty}$-structure on $\gen$ by 
  transferring, along the quasi-isomorphism $\Psi$,  the DGLA structure from   $\vC^\bullet(R, R)$ to  a $L_\infty$-structure on
   $\gen$. Explicit formulas for the transfer
    in terms of  sums over planar trees can be written similarly to those from \cite{MS}, Section 4. 
   Note that the transfered structure a priori depends on a choice of $\infty$
    and hence is different from the one described in Section 9. Because of our theorem,  these two $L_{\infty}$ structures are
     equivalent. In particular, the transfered structure does not, up to equivalence,  depend on a choice of $\infty$.
     
     On the other hand, the transfered structure has the advantage that it  can be 
      defined in terms of the  {\em oriented matroid}  generated by the
       set $\wt A$, i.e., by the knowledge of which pairs consisting of  a point $\omega\in\wt A$ and a subset $A'\subset \wt A$, 
       satisfy $\omega\in\Conv(A')$, see
      \cite{OM}. This is a much weaker and more combinatorial datum, in particular, it cannot be used to recover the concepts of a
       regular decomposition or  a convex piecewise affine function. 
    \end{rems}
    
    \vskip .2cm

\noindent {\sl Proof of the theorem:} For simplicity, let us write $\vC^\bullet$ for $\vC^{\geq 1}(R,R)[1]$ in the sequel.

 Let us prove the first statement of the theorem.  By Proposition \ref{prop:comp-phi},
the components of $\psi$ are the maps $\psi_{A'}^{(B_1, \cdots, B_m|C)} $ from
\eqref{eq:mat-el-psi} and \eqref{eq:mat-el-psi-trans} and we keep the corresponding notation. 
These are non-zero only if $F_{B_1}, \cdots, F_{B_m}$ are summands
of $R_{i_0i_1}$,  $R_{i_1 i_2}$, ..., $R_{i_{m-1} i_m}$ and $F_C$ is a summand of $R_{i_0 i_m}$.
They are also bilinear with respect to the intermediate algebras. This proves that $\psi$ factors
through a morphism $\Psi$ with values in  $\vC^\bullet$.

Let us now prove that $\Psi$ is a quasi-isomorphism. It will be convenient to organize the proof into several
steps.

\vskip .3cm

\noindent {\bf A. Interpretation of $\vC^\bullet$ via closed paths.}
Recall that $R_{ij}=\bigoplus_{Q'} F_{A\cap Q'}$ where 
  $Q'$ runs over infinite subpolygons with edges $[\infty, i]$ and $[ j, \infty]$.  
 
 \begin{Defi} A sequence  $P_0, P_1, \cdots, P_n$
  of marked infinite polygons  is called {\em admissible}, if its left and right infinite edges match as in Fig. 
  \ref{fig:adm-seq}, that is, there are $i_0 < \cdots < i_n$
  such that $P_0$ has edges $[\infty, i_0]$ and $[i_n, \infty]$ which $P_\nu$ has edges $[\infty, i_{\nu-1}]$
  and $[i_\nu, \infty]$. 
  Note that it is not required that $P_0$ actually contains any of the
  $P_\nu$.    
  \end{Defi} 
 
  \bef
 \btp[scale=0.35]
 \node (0) at (-1,0){};
  \node (1) at (3.5,0){};
   \node (2) at (8,0){};
    \node (n) at (14,0){};
     \node (i) at (6,10){};
     
 \fill (0) circle (0.15);  
  \fill (1) circle (0.15);  
   \fill (2) circle (0.15);  
    \fill (n) circle (0.15);  
     \fill (i) circle (0.15);  
 
 \draw (-1,0) -- (6,10) -- (3.5,0);
 \draw (8,0) -- (6,10) -- (14,0);
 
 \node at (6,11){$\infty$}; 
 
 \node at (-2,0){$i_0$}; 
 
 \node at (3.4, -1.5){$i_1$}; 
 
 \node at (8, -1.5){$i_2$}; 
 
 \node at (15,0){$i_n$}; 
 
 \node at (1.3,0){\small$P_1$};
 
 \node at (5.7,0){\small$P_2$};
 
 \node at (11,0){\small$P_n$};
 
 \node at (6,-3){\small$P_0$};

 \draw[line width=0.3mm] (-1,0) to[out=-120,in=-100] (3.5,0); 
  \draw[line width=0.3mm] (3.5,0) to[out=-120,in=-100] (8,0); 
   \draw[line width=0.3mm] (8,0) to[out=-120,in=-100] (14,0); 

 \draw[line width=0.3mm] (-1,0) to[out=-120,in=-100] (14,0); 

 \etp
 \btp[scale=0.35]
 \node (0) at (-1,0){};
   
    \node (n) at (14,0){};
     \node (i) at (6,10){};
     
 \fill (0) circle (0.15);  
      \fill (n) circle (0.15);  
     \fill (i) circle (0.15);  
 
 \draw (-1,0) -- (6,10) -- (14,0);
  
 \node at (6,11){$\infty$}; 
 
 \node at (-2,0){$i_0$};

 \node at (15,0){$i_n$};

 \node[align=center] at (6,2){\small $P_0$\\ is this\\big polygon};
 
 \draw[line width=0.3mm] (-1,0) to[out=-120,in=-100] (14,0); 
\etp
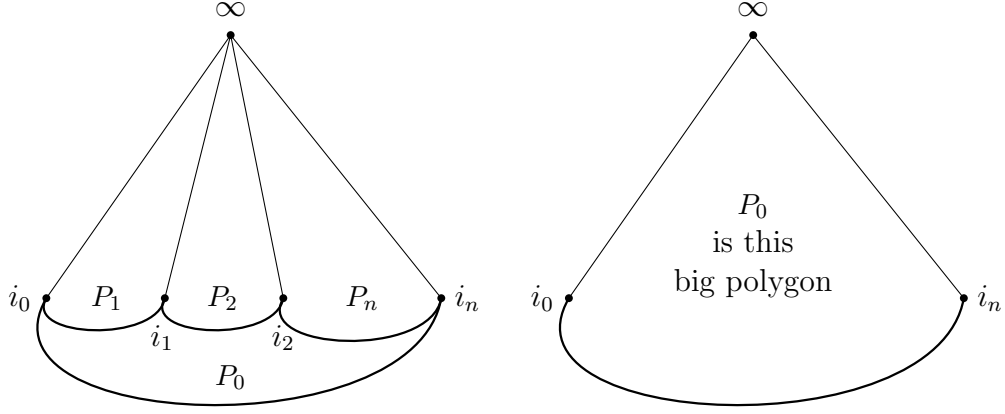
  \caption{ An admissible sequence. }
  \label{fig:adm-seq}
  \enf

   We  denote $B_\nu=\wt A\cap P_\nu$. 
   For an admissible sequence  $P_0, P_1, \cdots, P_n$ we define
  \[
  F_{P_1, \cdots, P_n}^{P_0} \,\,=\,\,\Hom_{S_{i_0}\otimes S_{i_n}^{\on{op}}}
 \biggl( F_{B_1}\otimes_{S_{i_1}} F_{B_2} \otimes_{S_{i_2}} \cdots \otimes_{S_{i_{n-1}}} F_{B_n}, \,\,
 F_{B_0}\biggr). 
  \]
 Then
 \[
 \vC^n \,\,=\,\,\bigoplus_{P_0, P_1, \cdots, P_n \text{ admiss.}} 
  F_{P_1, \cdots, P_n}^{P_0} . 
 \]
 We note further that $ F_{P_1, \cdots, P_n}^{P_0} $ is the same as the cyclic tensor product $N_\sigma$
 (in the sense of Remark \ref{rem:cyc-prod})
  over the closed edge path $\sigma$ obtained by running over the negative boundaries $\partial_- P_1, \partial_- P_2,
 \cdots, \partial_- P_n$ and then over $\partial_- P_0$ in the opposite direction, see Fig. \ref{fig:F-cyc-ten}. 
 
  \bef
  \tikzset{middlearrow/.style={
        decoration={markings,
            mark= at position 0.5 with {\arrow{#1}} ,
        },
        postaction={decorate}
    }
}

 \btp[scale=0.4]
 \node (0) at (-1,0){};
  \node (1) at (3.5,0){};
   \node (2) at (8,0){};
    \node (n) at (14,0){};
         
 \fill (0) circle (0.15);  
  \fill (1) circle (0.15);  
   \fill (2) circle (0.15);  
    \fill (n) circle (0.15);

 \node at (-2,0){$i_0$}; 
 
 \node at (3.4, -1.5){$i_1$}; 
 
 \node at (8, -1.5){$i_2$}; 
 
 \node at (15,0){$i_n$}; 
 
 \node at (1.3,0){\small$\partial_- P_1$};
 
 \node at (5.7,0){\small$\partial_- P_2$};
 
 \node at (11,0){\small$\partial_- P_n$};
 
 \node at (6,-3){\small$\partial_- P_0$};
 
 \draw[middlearrow={>}, line width=0.3mm] (-1,0) to[out=-120,in=-100] (3.5,0); 
  \draw[middlearrow={>},  line width=0.3mm] (3.5,0) to[out=-120,in=-100] (8,0); 
   \draw[middlearrow={>}, line width=0.3mm] (8,0) to[out=-120,in=-100] (14,0); 

 \draw[middlearrow={<}, line width=0.3mm] (-1,0) to[out=-120,in=-100] (14,0); 
 
 \node at (-7,-1){\huge$\sigma=$};
 
 \etp
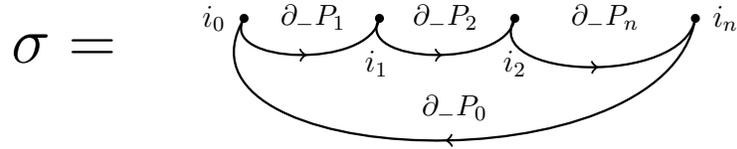
  \caption{ $F_{P_1, \cdots, P_n}^{P_0}$ as a cyclic tensor product $N_\gamma$. }
  \label{fig:F-cyc-ten}
  \enf
  
  \noindent {\bf B. Filtration of $\vC^\bullet$ by handle length.}
  Note that closed edge paths $\sigma$ labeling summands in $\vC^\bullet$, can retrace parts of themselves
  on the left or on the right. Similarly to \S \ref {sec:anal}, we will call such retraced parts {\em handles}
  and denote them by $\lambda$ and $\rho$,  see
  Fig. \ref{fig:handles}. 
  
    \bef
    \tikzset{middlearrow/.style={
        decoration={markings,
            mark= at position 0.5 with {\arrow{#1}} ,
        },
        postaction={decorate}
    }
}
   
     \btp[scale=0.35]
 \node (0) at (-1,0){};
  \node (1) at (3.5,0){};
   \node (2) at (8,0){};
    \node (n) at (14,0){};

 \fill (0) circle (0.15);  
  \fill (1) circle (0.15);  
   \fill (2) circle (0.15);  
    \fill (n) circle (0.15);  
    
    \node at (-4,2){\large$\lambda$}; 
      \node at (17,2){\large$\rho$}; 

\node at (1.5,0){$\eta_1$};  
\node at (6,0){$\eta_2 \,\,\,\cdots$}; 
\node at (11,0){$\eta_m$};

 \node (p) at (17,5){};
 \fill (p) circle (0.15);

 \draw[middlearrow={>}, line width=0.3mm] (-1,0) to[out=-50,in=-100] (3.5,0); 
  \draw[middlearrow={>}, line width=0.3mm] (3.5,0) to[out=-120,in=-100] (8,0); 
   \draw[middlearrow={>}, line width=0.3mm] (8,0) to[out=-120,in=-120] (14,0); 

 \draw[middlearrow={<}, line width=0.3mm] (-1,0) to[out=-50,in=-120] (14,0); 
 
  \draw[line width=0.3mm] (14,0) -- (17,5); 
  
  \node (q) at (-4,5){};
  \fill (q) circle (0.15); 
  
  \draw[line width=0.3mm] (-4,5) -- (-1,0); 
 \etp
   \caption{ A closed path with handles.}
   \label{fig:handles}
   \enf
   
    \begin{lem}
   Put
   \[
   G^l \vC^\bullet  \,\,=\,\,\bigoplus_n \bigoplus F_{P_1, \cdots, P_n}^{P_0},
   \]
   the second sum  running over  admissible sequences $P_0, P_1, \cdots, P_n$ with the sum of edge
   lengths of the two ``handles" on the left and right being $\geq l$. Then $G$ is a decreasing filtration of 
   $\vC^\bullet$ by subcomplexes. 
      \end{lem}
      
      \noindent {\sl Proof:}   The differential in $\vC^\bullet$ is the sum $d+\delta$, where $\delta$ is the Hochschild differential,
   induced by the multiplication in $R$, and $d$ is induced by the differential in $R$ itself. We recall the formula for
   \[
   \delta: \Hom(R^{\otimes n}, R) \lra \Hom(R^{\otimes(n+1)}, R)
   \]
   (in our case, $\Hom$ and $\otimes$ are taken over intermediate algebras, as explained above):
   \[
   \begin{gathered}
   (\delta f)(a_1, \cdots, a_{n+1}) \,\,=\,\, a_1 f(a_2, \cdots, a_{n+1}) \,\, + 
   \\
   +\,\,  \sum_{\nu=1}^n (-1)^\nu f(a_1, \cdots, a_\nu a_{\nu+1}, \cdots, a_{n+1}) \,\,+ \,\,
   (-1)^n f(a_1, \cdots, a_n) a_{n+1}. 
   \end{gathered}
   \]
   Recall also the formula for the multiplication $m_2$ in $R$:
   \[
   m_2|_{F_{B_1}\otimes F_{B_2}} = \begin{cases} \text{The canonical map } F_{B_1}\otimes F_{B_2}\to F_{B_1\cup B_2}, \text{ if }\\
   Q_1\cup Q_2 \text{ is convex and } Q_1  \text{ is on the left of } Q_2 ; \\
   0, \text{ otherwise.} 
   \end{cases}
   \]
   Here $Q_1$ and $Q_2$ are two   infinite subpolygons, and $B_\nu=\wt A\cap Q_\nu$.

   \vskip  .2cm

    Consider a summand $F_{P_1, \cdots, P_n}^{P_0}\subset \vC^n(R, R)$. The first and last terms in
   the formula for $\delta$, applied to this summand, consist in adding an infinite polygon $P$ to the picture on the left
   and on the right so that $\partial_-P$ together with $\partial_- P_0$ form a convex polygon (if the union is not
   convex, the contribution is $0$). See Fig. \ref{fig:delta-flast}. 
   
   \bef

    \btp[scale=0.35]
 \node (0) at (-1,0){};
  \node (1) at (3.5,0){};
   \node (2) at (8,0){};
    \node (n) at (14,0){};
     \node (i) at (6,10){};
     
 \fill (0) circle (0.15);  
  \fill (1) circle (0.15);  
   \fill (2) circle (0.15);  
    \fill (n) circle (0.15);  
     \fill (i) circle (0.15);  
 
 \draw (-1,0) -- (6,10) -- (3.5,0);
 \draw (8,0) -- (6,10) -- (14,0);
 
 \node at (6,11){$\infty$}; 
 
 \node (p) at (-1,5){};
 \fill (p) circle (0.15); 
 \draw (-1,5) -- (6,10);

 \node at (1.3,0){\small$P_1 $};
 
 \node at (5.7,0){\small$\cdots$};
 
 \node at (11,0){\small$P_n$};
 
 \node at (6,-3){\small$P_0$};
 
 \node at (0.5,4.5){\small$P$};
 
 \node at (5,-7){$a_1f(a_2,\cdots, a_{n+1})$};

 \draw[line width=0.3mm] (-1,0) to[out=-80,in=-100] (3.5,0); 
  \draw[line width=0.3mm] (3.5,0) to[out=-120,in=-100] (8,0); 
   \draw[line width=0.3mm] (8,0) to[out=-120,in=-100] (14,0); 

 \draw[line width=0.3mm] (-1,0) to[out=-80,in=-100] (14,0); 
 
  \draw[line width=0.3mm] (-1,5) -- (-1,0);

 \etp
    \btp[scale=0.35]
 \node (0) at (-1,0){};
  \node (1) at (3.5,0){};
   \node (2) at (8,0){};
    \node (n) at (14,0){};
     \node (i) at (6,10){};
     
 \fill (0) circle (0.15);  
  \fill (1) circle (0.15);  
   \fill (2) circle (0.15);  
    \fill (n) circle (0.15);  
     \fill (i) circle (0.15);  
 
 \draw (-1,0) -- (6,10) -- (3.5,0);
 \draw (8,0) -- (6,10) -- (14,0);
 
 \node at (6,11){$\infty$}; 
 
 \node (p) at (14,5){};
 \fill (p) circle (0.15); 
 \draw (14,5) -- (6,10);

 \node at (1.3,0){\small$P_1 $};
 
 \node at (5.7,0){\small$\cdots$};
 
 \node at (11,0){\small$P_n$};
 
 \node at (6,-3){\small$P_0$};
 
 \node at (12.3,4.5){\small$P$};
 
 \node at (5,-7){$f(a_1,\cdots, a_n)a_{n+1}$};

 \draw[line width=0.3mm] (-1,0) to[out=-120,in=-100] (3.5,0); 
  \draw[line width=0.3mm] (3.5,0) to[out=-120,in=-100] (8,0); 
   \draw[line width=0.3mm] (8,0) to[out=-120,in=-100] (14,0); 

 \draw[line width=0.3mm] (-1,0) to[out=-120,in=-100] (14,0); 
 
  \draw[line width=0.3mm] (14,0) -- (14,5);

 \etp
   \caption{The first and last terms in $\delta$.}
   \label{fig:delta-flast}
   \enf

   The new summands in which these terms lie, correspond to closed edge paths with the length of the left or right handle
   increased.  The remaining summands $(-1)^\nu f(a_1, \cdots, a_\nu a_{\nu+1}, \cdots, a_{n+1})$ of $\delta$ send 
   $F_{P_1, \cdots, P_n}^{P_0}$ into summands obtained by {\em all possible splittings} of each of the
   polygons $P_\nu$, $\nu=1, \cdots, n$, into two subpolygons along some edge of the form $[\infty, p]$
   where $p$ is some intermediate vertex. If $P_\nu$ is a triangle and cannot be split, then the corresponding
   contribution is $0$. See Fig. \ref{fig:middle}.    These summands have the same lengths of handles. 
   Therefore $\delta$ preserves the filtration $G$,
   and so does $d+\delta$. \qed

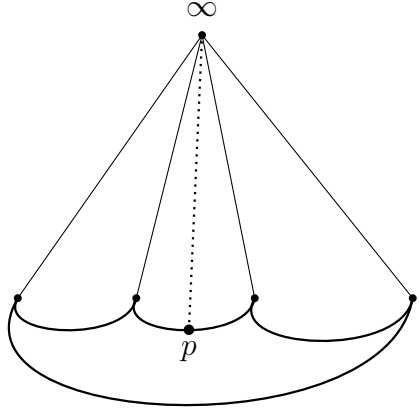
   \bef

 \btp[scale=0.35]
 \node (0) at (-1,0){};
  \node (1) at (3.5,0){};
   \node (2) at (8,0){};
    \node (n) at (14,0){};
     \node (i) at (6,10){};
     
 \fill (0) circle (0.15);  
  \fill (1) circle (0.15);  
   \fill (2) circle (0.15);  
    \fill (n) circle (0.15);  
     \fill (i) circle (0.15);  
 
 \draw (-1,0) -- (6,10) -- (3.5,0);
 \draw (8,0) -- (6,10) -- (14,0);
 
 \node at (6,11){$\infty$};

 \draw[line width=0.3mm] (-1,0) to[out=-120,in=-100] (3.5,0); 
  \draw[line width=0.3mm] (3.5,0) to[out=-120,in=-100] (8,0); 
   \draw[line width=0.3mm] (8,0) to[out=-120,in=-100] (14,0); 

 \draw[line width=0.3mm] (-1,0) to[out=-120,in=-100] (14,0); 
 
 \node (p) at (5.5,-1.2){};
 \fill(p) circle (0.2); 
 \draw[dotted, line width=0.3mm] (5.5,-1.2) -- (6,10); 

\node at (5.5,-2){$p$};

\node at (6,-6){\large $f(a_1, \cdots, a_\nu a_{\nu+1}, \cdots, a_{n+1})$};
 \etp
  \caption{ Intermediate terms in $\delta$.}
   \label{fig:middle}
   \enf

    \noindent {\bf C. $\ol\Psi: \gen\to\on{gr}^0\vC^\bullet$ is an embedding with exact quotient.}
Consider the morphism of complexes
\[
\ol\Psi: \gen \lra \on{gr}^0_G \vC^\bullet = \vC^\bullet/G^1 \vC^\bullet,
\]
induced by $\Psi$. We can think of $\on{gr}^0_G \vC^\bullet$ as the direct sum of cyclic tensor products
$N_\sigma$ corresponding to closed paths $\sigma$ without handles (and coming from admissible sequences).

  On the other hand,  $\gen$ is the direct sum of cyclic tensor products taken over the boundaries of
  (convex) finite subpolygons $Q'$. These boundaries are particular cases of the paths
  we obtain for $\on{gr}^0_G \vC^\bullet$. Further, it follows from Proposition \ref{prop:comp-phi}
  that the only matrix element of $\Psi$ on $E_{Q'}$ not landing in $G^1 \vC^\bullet$, is the transpose of the map
  $\gamma_{Q'}$ from \eqref{eq:gamma-Q'}. This transpose is in fact an isomorphism (identity)
  \[
  \gamma{Q'}^t: E_{A'}\lra N_{\partial Q'} = E_{A'}.
  \]
    This establishes the following.
  
  \begin{lem}
 $\ol\Psi$ is an embedding of complexes, whose image is the direct sum of summands  
 $ F_{P_1, \cdots, P_n}^{P_0} $ such that the negative boundaries  $\partial_- P_\nu$, $\nu=1,\cdots, n$,
 are segments which, together with $\partial_-P_0$, form a convex polygon. \qed
   \end{lem}
   
   So it is enough to prove that $\Coker(\ol\Psi)$ is exact. As a vector space, $\Coker(\ol\Psi)$,
   is simply the direct sum of all  the summands $ F_{P_1, \cdots, P_n}^{P_0}$ without handles, 
   other than those forming a convex polygon as above.

    Now, in $\on{gr}_F$ the first and last
      terms of the formula for $\delta$ disappear (they increase the sum of the lengths of the handles). The remaining
      intermediate terms correspond to splitting each $P_\nu$, $\nu=1,..., n$ along possible intermediate vertices  $p$ of $P_\nu$
      as in Fig. \ref{fig:middle}:
      \[
      P_\nu \,\,\rightsquigarrow \,\, P'_\nu \,\cup \, P''_\nu. 
      \]
      The summand corresponding to the split configuration is identical with the original one:
      \[
      F_{P_1, \cdots, P_{\nu-1}, P'_\nu, P''_\nu, P_{\nu+1}, \cdots,  P_n}^{P_0}\,\,=\,\, F_{P_1, \cdots, P_n}^{P_0}. 
      \]
      
      Thus $\on{gr}_F \vC^\bullet$ is split into a direct sum of complexes of the form
      \[
      F_{P_1, \cdots, P_n}^{P_0} \otimes_\k C^\bullet(\Delta^{r-1}),  
      \]
      where:
      \begin{itemize}
      
      \item[(a)]  $P_0, P_1, \cdots, P_n$ is an admissible sequence such that none of the
      intermediate adjacent pairs  $P_\nu, P_{\nu+1}$, $i+1\leq n$, can be combined 
      together to form a convex polygon (from which they could be obtained by splitting).
      
      \item[(b)] $r$ is the number of intermediate vertices of  $P_1, \cdots, P_n$ .
      
      \item[(c)] $\Delta^{r-1}$ is the $(r-1)$-dimensional combinatorial simplex whose vertices
      correspond to the intermediate vertices (the ways of splitting) in (b), and $C^\bullet(\Delta^{r-1})$
      is the augmented simplicial cochain complex of $\Delta^{r-1}$ with coefficients in $\k$. Note that 
      $C^\bullet(\Delta^{r-1})$ is exact if $r\geq 0$. 
      
      \end{itemize}
      
      So the remaining summands $F_{P_1, \cdots, P_n}^{P_0}$ which will not be tensored with an exact
      complex, correspond to admissible sequences $P_0, P_1, \cdots, P_n$  which:
      \begin{itemize}
      \item[(1)] Cannot be split anywhere, i.e., each $P_\nu$, $\nu=1,\cdots, n$, is a triangle, and so
      $\partial_-P_\nu$ is a segment. 
      
      \item[(2)] Cannot appear in the splitting of something else. This means that the segments
      $\partial_- P_1, \cdots, \partial_- P_n$ form an {\em upwardly convex} broken line from $i_0$ to $i_n$,
      which, together with $\partial_- P_0$ (a {\em downwardly convex} broken line) form a convex
      finite  polygon $Q'$ representing a summand in $\Psi(\gen)$, see Fig. \ref{fig:poly-app},
      where $\partial_-P_0$ depicted curved to emphasize that it can consist of many segments. 
   \end{itemize}
      
        \bef
   \btp[scale=0.4]
 \node (2) at (6,4){};
\node (3) at (2,6){};
\node (4) at (-1,6){};
\node (5) at (-6,4){};
\node (6) at (0,15){};

\fill(2) circle (0.15); 
\fill(3) circle (0.15); 
\fill(4) circle (0.15); 
\fill(5) circle (0.15); 
\fill(6) circle (0.15);

\draw (-6,4) -- (0,15) -- (-1,6); 
\draw (2,6) -- (0,15) -- (6,4); 

\draw[line width=0.4mm] (-6,4) -- (-1,6) -- (2,6) -- (6,4); 

\node at (0,16){$\infty$};
 
\node at (-3,4){$\partial_- P_1$};
\node at (0.5,5){$\partial_-P_2$};
\node at (4,4){$\partial_-P_n$};

  \draw[line width=0.4mm] (-6,4) to[out=-120,in=-100] (6,4); 
  
  \node at (0,1.6){$\partial_- P_0$};

 \etp

   \caption{ A convex polygon appearing.}
   \label{fig:poly-app}
   \enf
   
These are precisely the summands forming $\ol\Psi(\gen)$. This proves that $\ol\Psi$ is a quasi-isomorphism. 

\vskip .3cm

\noindent {\bf D. All $\on{gr}^{\geq 1}\vC^\bullet$ are exact.}    To finish the proof of   Theorem \ref{thm:uni},
it will now be sufficient to show that each $\on{gr}^l_G\vC^\bullet$, $l\geq 1$, is an exact complex.

This is done similarly to Step C. Since the first and last summands in $\delta$
disappear in  $\on{gr}_G$, the differential induced by $\delta$ on  $\on{gr}^l_G \vC^\bullet$
consists in splittings of the polygons $P_\nu$ into two subpolygons, as in Fig. \ref{fig:middle}. 
So $\on{gr}^l_G \vC^\bullet$ splits into a direct sum of sub complexes corresponding to different total
shapes of the pictures, and the differential in each such summand consists of splittings (the $\delta$ part)
plus the differential $d$ induced from $d_R$.

Now, if $l\geq 1$, the picture representing each summand has at least one handle,  $\lambda$ or $\rho$, 
nontrivial. If $\lambda$ is nontrivial, then the edge path $\sigma$ representing the picture
allows for nontrivial splittings of the union of $\lambda$ and $\eta_1$, see Fig. \ref{fig:handles}. 
As in part C, his will exhibit a nontrivial (i.e., contractible) $C^\bullet(\Delta^{r-1})$ as a tensor factor of our summand.
Similarly if $\rho$ is nontrivial. 
 
   Theorem \ref{thm:uni} is proved. 

\vfill\eject

\section{Comparison with Gaiotto-Moore-Witten.} \label{sec:GMW}

In this section  we provide a dictionary between the approach and terminology of \cite{GMW} and the ones
we use in the main body of the paper. 

\vskip .2cm

\noindent {\bf Vacua = elements of $A$.} The construction of \cite{GMW} starts with a finite set $\VV$ of ``vacua"
whose elements are denoted by $i,j,k$, etc., and a ``weight map" $z: \VV\to \CC$ (which we assume injective).
Our set $A$ is the image of this map. So in this section we denote elements of $A$ by $z_i, i\in\VV$. 

\vskip .2cm

\noindent {\bf Webs=polygonal decompositions.}  We recall:

\begin{Defi}(\cite{GMW})\label{def:web}
A {\em plane web} is a graph $\wen$ in the plane $\RR^2=\CC$ together with a labeling of faces
(``countries") by vacua so that: 
\begin{enumerate}
\item[(1)] Different faces have different labels.

\item[(2)] If an edge is oriented so that $i$ is on the left and $j$ on the right, then this edge is parallel to the vector
$z_i-z_j$ (i.e., has the same oriented direction but possibly different nonzero length). 

\end{enumerate}
\end{Defi}

\begin{rem}
More precisely, in \cite{GMW} the labels are allowed to be repeated, with (1) replaced by a weaker condition: labels
differ across edges. Such repeating patterns are accounted for by taking the symmetric algebra as in
\S \ref{sec:com-dg-alg}, and we do not consider them here. 

\end{rem}

\bef
 \btp[scale=0.4]

  \node (a) at (-1,3){};
  \node (b) at (2,4){};
  \node (c) at (4,2){};
  \node (d) at (7,3){};
  
  \fill (a) circle (0.15);  \fill (b) circle (0.15); 
 \fill (c) circle (0.15); 
 \fill (d) circle (0.15); 
 
 \node (1) at (-6,6){};
 \node (2) at (2,9){};
 \node (3) at (8,7){};
 \node (4) at (9,-1){};
 \node (5) at (2,-2){};
 \node (6) at (-4,-1){};
 
 \draw[line width=0.4mm, ->] (-1,3) -- (1);
  \draw[line width=0.4mm, ->] (2,4) -- (2);
   \draw[line width=0.4mm, ->] (7,3) -- (3);
    \draw[line width=0.4mm, ->] (7,3) -- (4);
     \draw[line width=0.4mm, ->] (4,2) -- (5);
      \draw[line width=0.4mm, ->] (-1,3) -- (6);
 
  \draw[line width=0.4mm] (-1,3) -- (2,4);
   \draw[line width=0.4mm] (2,4) -- (4,2);
    \draw[line width=0.4mm] (4,2) -- (7,3);
    
    \node at (-5,2){\large$n$}; 
    \node at (-1,6.5){\large$i$}; 
 \node at (5.3,6.8){\large$j$}; 
 \node at (10,3){\large$k$}; 
\node at (6, -1){\large$l$}; 
\node at (-.5,-.8){\large$m$};

 \etp

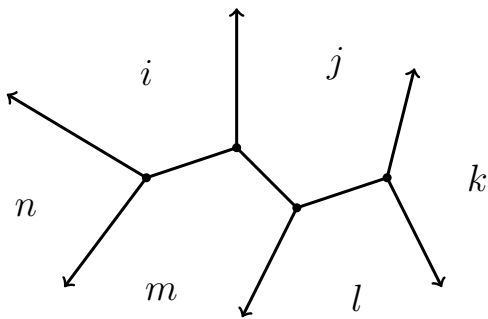
\caption{A plane web. Edges with arrows go to infinity. } 
\label{fig:web}
\enf

We call the set of vacua that appear as labels of the faces of a web $\wen$, the {\em support} of $\wen$. 
We further recall that two webs $\wen$ and $\wen'$ are said to be of the same {\em deformation type}, if they
are topologically equivalent in the following restricted sense: $\wen'$ can obtained from $\wen$ by translation
and stretching (but not rotating) of the edges. Following \cite{GMW}, we denote by $D(\wen)$ the set formed by
all webs of the same deformation type as $\wen$, considered modulo overall translations. This set is called
the {\em moduli space} of webs in the fixed deformation type. 

\begin{prop}\label{prop:websub}
Let $A'\subset A$ and $Q'=\Conv(A')$. Then:
\begin{itemize}

\item[(a)] Deformation types of webs with support $A'$ are in bijection with regular  geometric polygonal decompositions
$\Pc'$ of $(Q', A')$, i.e., with geometric faces  $F_{\Pc'}$ of the secondary polytope $\Sigma(A')$. 

\item[(b)] If $\wen$ is a web whose deformation type corresponds to a decomposition $\Pc'$, then
the moduli space $D(\wen)$ is identified with the reduced normal cone $\ol C_{\Pc'}$ to the face
$F_{\Pc'}\subset\Sigma(A')$, see \S \ref {sec:rem-sec}. 

\end{itemize}
\end{prop}

\noindent {\sl Proof:} Given a polygonal decomposition $\Pc'=\{(Q''_\nu, A''_\nu)\}$ of $(Q', A')$, we form the dual
graph by putting one vertex in each polygon $Q''_\nu$ and joining them as in Fig. \ref{fig:web-dec}:

\bef
 \btp[scale=0.4]
 
 \node (i) at (0,6){};
\node (j) at (6,6){};
 \node (k) at (9,3){};
 \node (l) at (6,0){};
 \node (m) at (0,0){};
 \node (n) at (-3,3){}; 
 
 \fill (i) circle (0.1); 
  \fill (j) circle (0.1);  \fill (k) circle (0.1);  \fill (l) circle (0.1);  \fill (m) circle (0.1);  \fill (n) circle (0.1); 
  
  \draw (0,6) -- (6,6) -- (9,3) -- (6,0) -- (0,0) -- (-3,3) -- (0,6); 
  
  \draw (0,6) -- (0,0) -- (6,6) -- (6,0); 
  
  \node (a) at (-1,3){};
  \node (b) at (2,4){};
  \node (c) at (4,2){};
  \node (d) at (7,3){};
  
  \fill (a) circle (0.15);  \fill (b) circle (0.15); 
 \fill (c) circle (0.15); 
 \fill (d) circle (0.15); 
 
 \node (1) at (-6,6){};
 \node (2) at (2,9){};
 \node (3) at (8,7){};
 \node (4) at (9,-1){};
 \node (5) at (2,-2){};
 \node (6) at (-4,-1){};
 
 \draw[line width=0.4mm, ->] (-1,3) -- (1);
  \draw[line width=0.4mm, ->] (2,4) -- (2);
   \draw[line width=0.4mm, ->] (7,3) -- (3);
    \draw[line width=0.4mm, ->] (7,3) -- (4);
     \draw[line width=0.4mm, ->] (4,2) -- (5);
      \draw[line width=0.4mm, ->] (-1,3) -- (6);
 
  \draw[line width=0.4mm] (-1,3) -- (2,4);
   \draw[line width=0.4mm] (2,4) -- (4,2);
    \draw[line width=0.4mm] (4,2) -- (7,3);
    
    \node at (-4,3){$z_n$}; 
    \node at (-1,6.5){$z_i$}; 
 \node at (6,6.8){$z_j$}; 
 \node at (10,3){$z_k$}; 
\node at (6.5, -1){$z_l$}; 
\node at (-.8,-.8){$z_m$};

 \etp

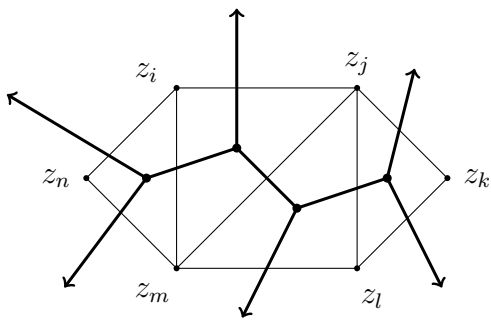
\caption{Webs are dual to decompositions.}
\label{fig:web-dec}
\enf

This gives a plane graph $\Gamma$ with faces labelled by elements of $A'$ but not necessarily satisfying the condition
(2) of Definition \ref{def:web}.  We now show that possible ways of deforming $\Gamma$ to a web, modulo translations,
are in bijection with $\ol C_{\Pc'}$. 

We first recall that $\ol C_{\Pc'}$ consists of convex $\Pc'$-piecewise affine functions $Q'\to\RR$, considered modulo adding global
affine functions on $Q'$, i.e., $\ol C_{\Pc'}=C_{\Pc'}/\Aff(\RR^2)$, where $C_{\Pc'}$ consists of actual 
convex $\Pc'$-piecewise affine functions. 

Next, we note that rotation by $90^\circ$ transforms the condition
(2) of Definition \ref{def:web} into:

\begin{enumerate}
\item[(2$'$)] An edge having $i$ on the left and $j$ on the right, is {\em perpendicular} to the vector $z_i-z_j$
(i.e., forms, together with this vector, a positively oriented orthogonal frame). 
\end{enumerate}

So we will work with (2$'$) instead of (2). Let us identify $(\RR^2)^*$ with $\RR^2$ by means of the standard scalar
product. 

Suppose given a piecewise-linear $f\in C_{\Pc'}$ (an actual function,  not a class modulo adding global affine functions). 
For each polygon $Q''_\nu\in\Pc'$ let $p_\nu\in(\RR^2)^*$ be the slope of the affine function $f|_{Q''_\nu}$.
By the above, we can regard $p_\nu$ as a point of $\RR^2$. If two polygons $Q''_\mu, Q''_\nu\in\Pc'$
have a common edge, the slopes of $f$ must match on this edge, This simply means that $p_\mu-p_\nu$
is perpendicular to the edge. 

So we obtain a plane graph $\wen$ topologically equivalent to $\Gamma$ by joining $p_\mu$ and $p_\nu$
by a straight line internal whenever $Q''_\mu$ and $Q''_\nu$ have a common edge By the above, $\wen$ satisfies
(2$'$), i.e.,  it is a ($90^\circ$  rotation of a) web. In this way the condition (2$'$) expresses the existence of a global
piecewise affine function $f$ with the slope $p_\nu$ on each $Q''_\nu$. The function $f$ is defined by the knowledge of
the $p_\nu$ up to an additive constant. Further, the {\em convexity} of $f$ is expressed by saying that the edges
of $\wen$ (with respect to our choice of orientations) have {\em positive length} (a concave break would result
in the length of an edge counted as negative). Translating $\wen$ as a whole amounts to adding  a global
linear function to $f$. We leave the remaining details to the reader. \qed

\vskip .2cm

We now briefly indicate the meaning, in our terms,  of some further terminology of \cite{GMW}.

\vskip .2cm

\noindent {\bf Cyclic fans of vacua = subpolygons $Q'=\Conv(A')$, $A'\subset A$.} Indeed, according to
{\em loc. cit.}, a cyclic fan of vacua is a cyclically ordered set $I=\{i_1, \cdots, i_n\}$ such that the
rays $(z_{i_k}-z_{i_{k+1}})\RR_+$ are oriented clockwise. This simply corresponds to $z_{i_1}, \cdots z_{i_n}$
being the vertices of a convex $n$-gon, in this cyclic order.

Further, if a web $\wen$ corresponds to a regular polygonal subdivision $\Pc'$ of $(Q', A')$  then
vertices of $\wen$ are labeled, by the above, by polygons $Q''_\nu$ of $\Pc'$, and 
 the local fan of $\wen$ at
a vertex $v$, denoted $I_v(\wen)$, corresponds  the polygons $Q''_\nu$ labeling $v$.
 The fan at
infinity, denoted $I_\infty(\wen)$,  corresponds to $Q'$ itself. 

\vskip .2cm

\noindent {\bf Rigid, taut and sliding webs.}  They correspond to 
faces of a secondary polytope $\Sigma(A'), A'\subset A$ of codimension 0, 1, and 2.
This is a consequence of Proposition \ref{prop:websub}. 

\vskip .2cm

\noindent {\bf Convolution of webs.} It corresponds to further subdivision of one polygon of a polygonal decomposition:

\bef
\btp[scale=0.3]
 
 \node (i) at (0,6){};
\node (j) at (6,6){};
 \node (k) at (9,3){};
 \node (l) at (6,0){};
 \node (m) at (0,0){};
 \node (n) at (-3,3){}; 
 
 \fill (i) circle (0.1); 
  \fill (j) circle (0.1);  \fill (k) circle (0.1);  \fill (l) circle (0.1);  \fill (m) circle (0.1);  \fill (n) circle (0.1); 
  
  \draw (0,6) -- (6,6) -- (9,3) -- (6,0) -- (0,0) -- (-3,3) -- (0,6); 
  
  \draw (0,6) -- (0,0);
  \draw (6,6) -- (6,0); 
   \etp
   \btp[scale=0.3]
   \node at (0,0){};
   \node at (0,6){}; 
   \node at (0,3){$*$}; 
   \etp 
   \btp[scale=0.3]
 
 \node (i) at (0,6){};
\node (j) at (6,6){};
  \node (l) at (6,0){};
 \node (m) at (0,0){};
  
 \fill (i) circle (0.1); 
  \fill (j) circle (0.1);     \fill (l) circle (0.1);  \fill (m) circle (0.1);     
    
  \draw (0,6) -- (0,0) -- (6,6) -- (6,0); 
  \draw (0,6) -- (6,6);
  \draw (0,0) -- (6,0); 
 \etp
   \btp[scale=0.3]
   \node at (0,0){};
   \node at (0,6){}; 
   \node at (0,3){$=$}; 
   \etp 
   \btp[scale=0.3]
 
 \node (i) at (0,6){};
\node (j) at (6,6){};
 \node (k) at (9,3){};
 \node (l) at (6,0){};
 \node (m) at (0,0){};
 \node (n) at (-3,3){}; 
 
 \fill (i) circle (0.1); 
  \fill (j) circle (0.1);  \fill (k) circle (0.1);  \fill (l) circle (0.1);  \fill (m) circle (0.1);  \fill (n) circle (0.1); 
  
  \draw (0,6) -- (6,6) -- (9,3) -- (6,0) -- (0,0) -- (-3,3) -- (0,6); 
  
  \draw (0,6) -- (0,0);
  \draw (6,6) -- (6,0); 
  \draw (0,0) -- (6,6); 
   \etp

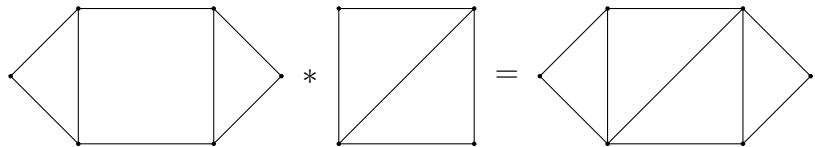
\caption{Convolution of webs in dual language.}
\label{fig:conv-web}

\enf 

\noindent {\bf Half-plane webs.} They correspond, in the dual language, 
to polygonal subdivisions of the extended polytope $\wt Q=\Conv(\wt A)$
where $\wt A = A\cup\{\infty\}$, see Figs.  \ref{fig:extended-polyt} and \ref{fig:half}. The choice of the direction towards $\infty$
corresponds to the choice of a a half-plane $H$ in \cite{GMW}. In Fig. \ref{fig:half}, $H$ is positioned horizontally. 

\bef
 \centering
 \btp[scale=.4, baseline=(current  bounding  box.center)]
 
 \node (1) at (0,0){};
 \fill (1) circle (0.1); 
 
  \node (2) at (4,0){};
 \fill (2) circle (0.1); 
 
  \node (3) at (6,3){};
 \fill (3) circle (0.1); 
 
  \node (4) at (3,5){};
 \fill (4) circle (0.1); 
 
  \node (5) at (-2,3){};
 \fill (5) circle (0.1); 
 
 \draw (-2,3) -- (0,0) -- (4,0) -- (6,3); 
 
 \draw (-2,3) -- (-2,12);
 
 \draw (6,3) -- (6,12); 
 \draw (0,0) -- (3,5) -- (6,3); 
 \draw (3,5) -- (3,12); 
 
 \node (a) at (0,6){};
 \node (b) at (3.5,2){};
 \node (c) at (4.5,6){}; 
 \fill (a) circle (0.15); 
  \fill (b) circle (0.15); 
   \fill (c) circle (0.15); 
 
 \draw[line width=0.4mm] (0,6) -- (3.5,2);  
  \draw[line width=0.4mm] (0,6) -- (4.5,6); 
  \draw[line width=0.4mm] (3.5,2) -- (4.5,6); 
    \draw[line width=0.4mm, ->] (0,6) -- (0,11); 
        \draw[line width=0.4mm, ->] (4.5,6) -- (4.5,11); 
     \draw[line width=0.4mm, ->] (0,6) -- (-2,0); 
         \draw[line width=0.4mm, ->] (3.5,2) -- (1,-3); 
   \draw[line width=0.4mm, ->]  (3.5,2) -- (8, -1); 
    \draw[line width=0.4mm, dotted] (-6,7) -- (10,7);   
  \filldraw[opacity=0.2] (-6,7) rectangle (10,12);  
 
 \node at (2,13){$\cdots \infty\cdots$};
 
 \node at (9,8){$H$};

 \etp
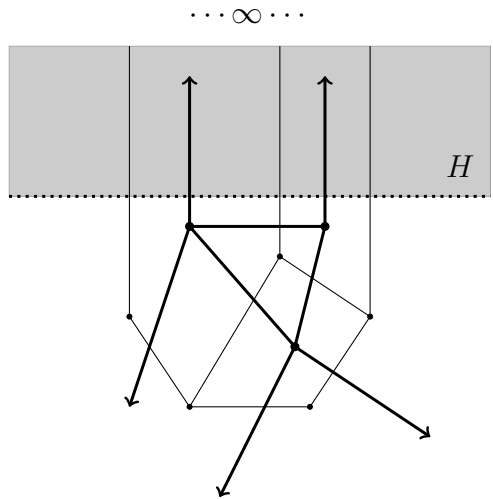
 \caption{A half-plane web in the dual language.}
\label {fig:half}
\enf

\vfill\eject

 \section{Maurer-Cartan elements for the Fukaya-Seidel category} \label{sec:MCFS} 
   
  In this section we  sketch, using our approach,
   an application of the previous considerations to the Fukaya-Seidel (FS)  categories. We assume familiarity with
   the conventional approach to the subject \cite{FSbook}, so our reminder in that respect will be minimal.

   \vskip .2cm
   
   \noindent {\bf A. Setup.} 
   Let $X$ be  a complex K\"ahler manifold $X$ of complex dimension $k$. The K\"ahler metric, denoted
   $\omega^{1,1}$, can be understood 
    as a family of complex valued Hermitian scalar products in the $\CC$-vector spaces $T_xX$. Globally,
    it gives a smooth isomorphism of $\CC$-vector bundles
    \be\label{eq:K-form}
    \omega^{1,1}: \ol TX \lra T^*X,
    \ee
    where $\ol TX$ is the conjugate bundle to the $\CC$-vector bundle $TX$. Separating the real and imaginary parts,
    we write
      $\omega^{1,1} = g+\omega\sqrt{-1}$. Thus $g$ is a Riemannian metric and $\omega$ is a symplectic form on $X$.
 We also assume that $X$ is equipped with a nowhere vanishing holomorphic volume form $\Omega^{k,0}$ (to be used later
   to obtain  $\ZZ$-gradings of various Floer complexes). 
   
   \vskip .2cm
   
   Let $W:X\to\CC$ be a holomorphic function with only non-degenerate (Morse)  critical points. 
   These points correspond to vacua of the Landau-Ginzburg theory associated to $W$. We denote the
   set of critical points by $\VV$ and assume it finite: $\VV=\{x_1, \cdots, x_r\}$. 
    We denote by $w_i=W(x_i)$ the critical value at $x_i$ and put $A=\{w_1, \cdots, w_r\}\subset\CC$. We assume that
   $A$ is in (affine) general position: no three points lie on a line. By identifying $\CC$ with $\RR^2$, we can
   apply the formalism of previous sections to this choice of $A$.  In particular, we have the polygon $Q=\Conv(A)$ and
   can speak about its triangulations, its secondary polytope etc.    
   
   \vskip .2cm

   The real part $\Ree(W): X\to\RR$ is a Morse function on $X$ (as a real $2k$-dimensional manifold), with the same
   critical points $x_1, \cdots, x_r$. Each $x_i$ is Morse with the same  signature $(k,k)$. As usual in Morse theory,
   by using the Riemannian metric $g$, we associate to $x_i$ the {\em unstable manifold} (also known in this
   context as the
   {\em Lefschetz thimble}) $T_i=T_i(W)\simeq\RR^k\subset X$. That is, $T_i$ is the union of downward gradient trajectories
   $\phi(t)$ originating from $x_i$, i.e., such that $\lim_{t\to -\infty} \phi(t) = x_i$. We recall the following well known
   
   \begin{prop}
   (a) The gradient flow of $\Ree(W)$ with respect to $g$ is equal to the Hamiltonian flow of $\Imm(W)$ with respect to
   $\omega$.
   
   (b) Each $T_i$ is Lagrangian with respect to $\omega$. \qed
   \end{prop}
   
   Since every Hamiltonian flow preserves the Hamiltonian, part (a) implies that the images $W(T_i)$ lie on horizontal
   half-lines $w_i+ \RR_{\leq 0}\subset\CC$. See Fig. \ref{fig:thimbles}, where we also indicated the hyperplane $H$
   describing the ``direction
   towards infinity",  cf. Fig. \ref{fig:half}. Its boundary is now positioned vertically,
   
   \bef
    \btp[scale=.4, baseline=(current  bounding  box.center)]
    \node (i) at (5,6){};
    \node (j) at (10,4){};
    \node (k) at (8,2){};
    \node (l) at (12, 0){};
    \fill (i) circle (0.15); 
       \fill (j) circle (0.15);    \fill (k) circle (0.15);    \fill (l) circle (0.15); 
    
 \draw[line width=0.4mm] (0,6) -- (5,6);
 \draw[line width=0.4mm]     (0,4) -- (10,4);
  \draw[line width=0.4mm] (0,2) -- (8,2);
   \draw[line width=0.4mm](0,0) -- (12,0);
   \filldraw[opacity=0.2] (-3,-2) rectangle (3,8);  
 \draw[line width=0.4mm, dotted]    (3,-2) -- (3,8); 
 \node at (-2,7){$H$}; 
 
 \node at (6,6){$w_i$};
  \node at (11,4){$w_j$};
   \node at (9,2){$...$};
    \node at (13,0){$w_l$};

    \etp
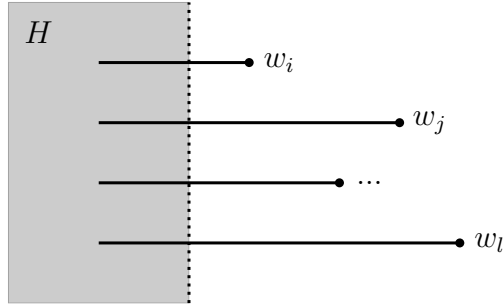
   \caption{Thimbles project to half-lines.}
   \label{fig:thimbles}
   \enf
   
   \begin{Defi}
   
   The set $A$ is called {\em horizontally generic}, if no two points of $A$ have the same imaginary part
   \footnote{In other words,  the extended set $\wt A = A\cup \{-\infty\}$ is in general position, see \S \ref{sec:rel-set} and
   Fig.\ref{fig:extended-polyt}.  }. For an angle $\theta$ we say that $A$ is $\theta$-generic, if $e^{\sqrt{-1}\theta} A$ is
   horizontally generic. 
   \end{Defi}
   
   We assume that $A$ is horizontally generic (which can be always assured by replacing $W$ with some $e^{ \sqrt{-1}\theta} W$).
   Then the half-lines $w_i+\RR_{\leq 0}$ are disjoint, and so the thimbles $T_i$ are disjoint closed Lagrangian submanifolds
   of $X$, which we assume numbered in order of increasing imaginary parts: $\Imm(w_1) < \cdots < \Imm(w_r)$. 
   We recall that $\k$ is our fixed base field of characteristic $0$. 
   The {\em Fukaya-Seidel category} $\FS=\FS(X,W)$ is a $\k$-linear $A_\infty$-category in which $T_1, \cdots, T_r$
   are objects, forming an exceptional collection:
   \[
   \Hom^\bullet_{\FS}(T_i, T_i) = \k, \quad \Hom^\bullet_{\FS}(T_i, T_j)=0, \,\, i>j. 
   \]
   We denote
   \[
   \RS = \RS(X,W) \,\,=\,\,\bigoplus_{i\leq j} \Hom^\bullet_{\FS}(T_i, T_j)
   \]
   the total $A_\infty$-algebra associated to $\FS$. The goal is to obtain $\RS$ by deforming an appropriate algebra
   $R$ as in \S \ref {sec:anal} with respect to an appropriate Maurer-Cartan element. 
   
   \begin{rem}
   There is no doubt that the formalism of Fukaya-Seidel categories generalizes to a situation much more general than
   the one described above: we must be able to allow $W$ to have non-Morse (and even non-isolated) critical points.
   In this case, each critical value $w_i$ (their number will still be finite in any algebraic situation) would give rise
   to a ``local Fukaya-Seidel category" $\Cc_i$. These local categories would then embed into a global category
   $\FS$ and form there a semi-orthogonal family. We expect that this more general setting to be still analyzable by the methods
   of this paper. In particular, the generality of extended systems of coefficients in \S \ref{sec:bim} is (intentionally) geared to this
   future context, the category $\Cc_i$ corresponding to the dg-algebra $S_i$, see Remark \ref{rem:Si-Fij}.

   \end{rem}
   
   \vskip .3cm
   
   \noindent {\bf B. The  system of coefficients.} The algebra $R$ comes from an extended coefficient system
   $(S_i, N^\bullet_{ij})$ as in \S \ref{sec:bim}. We put $S_i=\k$ (so the setting of \S \ref{sec:coeff-arbitr} is in fact 
   sufficient for this section). 
   This corresponds to the fact that each $x_i$
   is a Morse singular point for $W$. The graded vector spaces $N_{ij}^\bullet$ are obtained by analyzing
    those angles $\theta$
   for which $W$ is {\em not} $\theta$-generic. These are precisely
   \[
   \theta_{ij} = -\arg(w_i-w_j) \,\,\in\,\, S^1 = \RR/2\pi\ZZ, \quad i\neq j. 
   \] 
   We associate to them the unit vectors
   \[
   \zeta_{ij}\,\,=\,\, e^{\sqrt{-1}\theta_{ij}} \,\,=\,\, \biggl(\frac{w_i-w_j}{|w_i-w_j|}\biggr)^{-1}. 
   \]
   Note that for any $\theta\in S^1$ the function $\Ree(e^{\sqrt{-1}\theta} W): X\to\RR$ is again a Morse function
   with critical points $x_1, \cdots, x_n$ of signature $(k,k)$ and so gives rise to the ``$\theta$-rotated" thimbles
   $T_i^\theta = T_i(e^{\sqrt{-1}\theta} W)$. The image of $T_i^\theta$ under $W$ lies  now on the half-line
   originating from $w_i$ and having angle $-\theta$ with the real axis.  Following \cite{GMW} we introduce the
   following

   \begin{Defi}(a) By a  {\em $\zeta_{ij}$-soliton} we will mean a gradient trajectory $\phi=\phi(t)$ of $\Ree(\zeta_{ij} W)$
   originating (at $t=-\infty$)
   from $x_i$ and terminating (at $t=+\infty$) at $x_j$. Thus the image of a $\zeta_{ij}$-soliton under $W$
    is the straight interval $[w_i, w_j]$.
   
   (b) By a {\em gradient polygon} we mean a sequence  $\bphi = (\phi_{i_0i_1}, \phi_{i_1i_2}, \cdots,\phi_{i_mi_0})$,
  where $\phi_{i_\nu i_{\nu+1}}$ is a $\zeta_{i_\nu i_{\nu+1}}$-soliton. A gradient polygon is called {\em convex}, or
  a {\em cyclic fan of solitons}, if
  the intervals $[w_{i_\nu}, w_{i_{\nu+1}}]$ form the boundary of a convex polygon, in the counter-clockwise
  order. 
   \end{Defi}
   
   Thus the union of all  $\zeta_{ij}$-solitons is the intersection of the thimbles 
    $T_{ij}^{\theta_{ij}}$ and $T_j^{\theta_{ji}}$     Since existence of gradient trajectories between critical
   points of the same index is a codimension 1 phenomenon in the space of Morse functions,
   see, e.g., \cite{kapranov-saito},   the number of
   $\zeta_{ij}$-solitons is ``typically" finite. We assume this to be the case and denote this number by $n_{ij}=n_{ji}$. 
   Further, let $N_{ij}$ be
     the $\k$-vector space  having, as basis vectors, symbols $e_\phi$
   corresponding to the  $\zeta_{ij}$-solitons $\phi$.   So $\dim(N_{ij})=n_{ij}$.    
   
   We now explain how to introduce a $\ZZ$-grading on each $N_{ij}$. As usual with
   Floer complexes (of which the $N_{ij}$ are particular cases, see below), the grading
  is not fully canonical.  More canonical are the following two types of constraints, which
  admit a solution but not uniquely. 
  \vskip .2cm
  
  \noindent {\bf Duality:} $N_{ji}=N_{ij}^*$ as a graded space. Cf. the requirement (4) of Definition
  \ref{def:ex-sys-co} .

  \vskip .2cm
  
  \noindent {\bf Maslov index grading of cyclic products:} The grading in each tensor product
   of the form  $N_{i_0i_1}\otimes \cdots N_{i_{m-1} i_m}\otimes N_{i_mi_0}$,
  $m\geq 1$ (cf. Fig. \ref{fig:cyc-prod}, where all $S_i=\k$), is given topologically as follows.  
By definition, this tensor product has a basis $\{ e_{\bphi}\}$ labelled by gradient polygons $\bphi$,
and we define the degree of each $e_{\bphi}$. This is done by an instance of the Maslov index construction,
which we review here in the generality that we need, referring to   \cite{FSbook}, \S 11 for the general case. 

For a topological space $Y$ and $x,y\in Y$ we denote by $\Pi_{x,y}Y$ the set of homotopy classes
of paths in $Y$ beginning at $x$ an ending at $y$, i.e., the set of morphisms from $x$ to $y$ in the fundamental
groupoid of $X$.

  For a symplectic $\RR$-vector bundle $E$ over some base $B$ we denote by $\LG(E)$ the bundle of the Lagrangian
  Grassmannians of the fibers of $E$ over $B$. 
  
  Consider first the case when $B=\on{pt}$, so $E$ is a symplectic vector space. Let $L_1, L_2\in\LG(E)$ be two Lagrangian subspaces
  meeting transversely, i.e., $E=L_1\oplus L_2$. In this case there is a canonical homotopy class $\kappa_{L_1, L_2}
  \in \Pi_{L_1, L_2}\LG(E)$, defined uniquely by three properties:
  \begin{enumerate}
  \item[(1)] Naturality with respect to symplectic isomorphisms $E\to E'$.
  
  \item[(2)] Additivity in symplectic direct sums. That is,  if $E=E'\oplus E''$, and $L_\nu = L'_\nu\oplus L''_\nu$, then
  $\kappa_{L_1, L_2}$ is the image of $(\kappa_{L'_1, L'_2}, \kappa_{L''_1, L''_2})$ under the morphism of fundamental
  groupoids induced by $\LG(E')\times \LG(E'')\hookrightarrow\LG(E)$. 
  
  \item[(3)] Normalization in dimension 2. If $\dim_\RR(E)=2$, then $\LG(E)\simeq \RR ^1=S^1$ is a circle, equipped with orientation.
  In this case $\kappa_{L_1, L_2}$ is represented by the circle arc going from $L_1$ to $L_2$ in the positive direction. 
  
  \end{enumerate}
  
  \noindent Let us call $\kappa_{L_1, L_2}$ the {\em symplectic bridge} between $L_1$ and $L_2$.
  
  \vskip .2cm

  Returning now to our situation, 
each gradient polygon $\bphi$ gives a piecewise continuous
  closed path $\varpi_{\bphi}$ in $\LG(TX)$ obtained by associating to any internal point $p$ of the segment
  $\phi_{i_\nu i_{\nu+1}}$ the tangent space at $p$ to the thimble $T_{i_\nu}(\zeta_{i_\nu i_{\nu+1}}W)$. 
  At each point $i_\nu$, the path  $\varpi_{\bphi}$ has a simple discontinuity (jump), with the two limit
  values before and after the jump forming 
  a pair of transversal Lagrangian subspaces in $T_{x_\nu}X$. Joining each such pair by the symplectic
  bridge, we obtain a closed continuous path  $\varpi_{\bphi}^\wedge$ in $\LG(TX)$, defined
  uniquely up to homotopy.

  Taking the exterior power and applying the holomorphic volume form $\Omega^{k,0}$ gives a continuous map (``phase")
  \[
  \LG(TX) \buildrel \Lambda^k_\RR\over\lra \LG(\Lambda^k_\CC TX) \buildrel\Omega^{k,0}\over\lra \LG(\CC)=\RR P^1 = S^1.  
  \]
  Denote this composite map $\Theta$. Note that $\LG(\CC)=S^1$ has a canonical counterclockwise orientation. 
    
  \begin{Defi}  
  The Maslov index $d(\bphi)$ of a gradient polygon $\bphi$ is the winding number of the
  continuous closed path $\Theta(\varpi_{\bphi}^\wedge)$
   in $S^1=\LG(\CC)$. We assign to the basis vector $e_{\bphi}$
  the degree $d(\bphi)$. 
  
  \end{Defi}

   \vskip .3cm
   
   \noindent {\bf C. Transversality assumptions.} 
    Note that the intersection of the Lagrangian manifolds  $T_{ij}^{\theta_{ij}}$ and $T_j^{\theta_{ji}}$ is never transversal:
   if non-empty, it has dimension at least one. There are two ways of realizing $n_{ij}$ as the number
   of intersection points of two Lagrangian manifolds which can be (and, ``typically", are) transversal.
   
   \begin{enumerate}
   \item[(1)] Take an intermediate point $p$ on the interval $[w_i, w_j]$ and consider the (smooth) K\"ahler manifold
   $W^{-1}(p)$ of dimension $k-1$. The intersections
   \[
   V_i(p) = T_i^{\theta_{ij}}\cap W^{-1}(p), \quad  V_j(p) = T_j^{\theta_{ji}}\cap W^{-1}(p)
   \]
   are Lagrangian spheres (``vanishing cycles") in $W^{-1}(p)$ whose intersection points are in bijection with trajectories
   constituting $T_{ij}^{\theta_{ij}}$ and $T_j^{\theta_{ji}}$. See Fig. \ref{fig:vancycles}. We assume that $V_i(p)$ and $V_j(p)$
   intersect transversally in $W^{-1}(p)$.

   \bef
   \btp[scale=0.3]
   \node (i) at (-9,0){};
   \node (j) at (11,0){};
   \node (p) at (0,0){};
   \fill (i) circle (0.2);
    \fill (j) circle (0.2);
 \fill (p) circle (0.2);
 \node (xi) at (-9,9){};
 \node(xj) at (11,7){};
 \fill (xi) circle (0.2);
 \fill (xj) circle (0.2);
 
 \draw[ line width=.4mm] (0,7) ellipse (1cm and 2cm);
  \draw[ line width=.4mm] (0,10) ellipse (1cm and 2cm);

  \draw[line width=0.3mm] (-9,0) -- (11,0);
  
  \draw (-3,2) -- (-3,12) -- (3,16) -- (3,6) -- (-3,2); 
  
  \draw plot [smooth, tension=1]  coordinates {(0,12)  (-9,9)  (0,8)};  
  
    \draw plot [smooth, tension=1]  coordinates {(0,9)  (11,7)  (0,5)};  

\node at (-11,0) {\large$w_i$};
\node at (13,0) {\large$w_j$};
\node at (0,-2){\large$p$};
\node at (-11,9) {\large$x_i$};
\node at (13,7) {\large$x_j$};

\node at (2,17.5){\large$W^{-1}(p)$}; 

\node at (-7,12){$T_i^{\theta_{ij}}$};
\node at (9,9.5){$T_i^{\theta_{ij}}$};
\node at (1,13){\small$V_i(p)$}; 
\node at (-2.3,4.7){\small$V_j(p)$}; 

   \etp
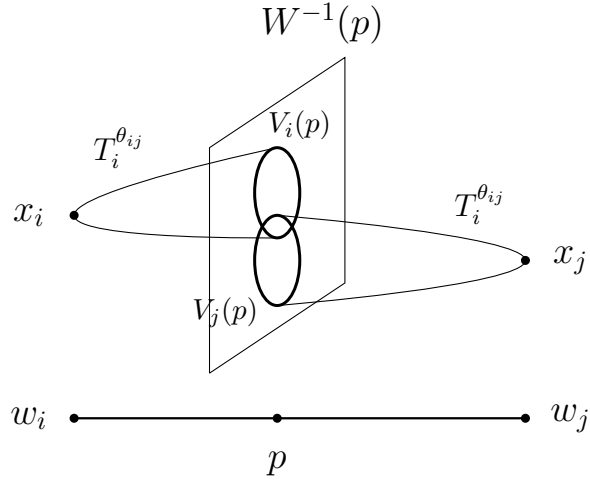
   \caption{Vanishing cycles over a midpoint.}
\label{fig:vancycles}
\enf
   
    \item[(2)] For small $\epsilon>0$ put $\theta'=\theta_{ij}-\epsilon$ and $\theta''=\theta_{ji}+\epsilon$. Then
   the thimbles $T_i^{\theta'}$ and $T_j^{\theta''}$ project  by $W$ onto two half lines which intersect in a point $q$ somewhere near
   the interval $[w_i, w_j]$, see Fig. \ref {fig:pert}. 
 The intersection points of  $T_i^{\theta'}$ and $T_j^{\theta''}$, all lying over $q$, 
  are again in bijection with trajectories constituting 
   $T_{ij}^{\theta_{ij}}\cap T_j^{\theta_{ji}}$. We assume (this assumption is equivalent to that in  (1)) that $T_i^{\theta'}$ and $T_j^{\theta''}$
   intersect transversally in $X$. 
   
   \bef
   \btp[scale=0.3]
   
    \node (i) at (-9,0){};
   \node (j) at (11,0){};
    \fill (i) circle (0.25);
    \fill (j) circle (0.25);
      \draw[line width=0.3mm] (-9,0) -- (11,0);
      \node at (-11,0) {\large$w_i$};
\node at (13,0) {\large$w_j$};

\node (q) at (1,2){};
\fill (q) circle (0.22){};

\draw[line width=0.4mm] (-9,0) -- (13,4.5); 
\draw[line width=0.4mm] (11,0) -- (-11,4.5); 
\node at (1,3){\large$q$}; 
\node at (12,6){$W(T_i^{\theta'})$};
\node at (-11,6){$W(T_j^{\theta''})$};
 \etp
 
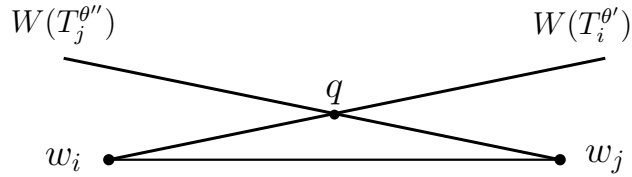
 \caption{Perturbed thimbles in projection to $\CC$.}
 \label{fig:pert}
   \enf
   
   \end{enumerate}

   \vskip .3cm
   
   \noindent {\bf D. The instanton equation.} Suppose that the K\"ahler manifold $X$ is exact, i.e., there exists a
   1-form $\lambda$ primitive for $\omega$, i.e., $d\lambda =\omega$. It is convenient to write $\lambda$
   symbolically as $\lambda = {\mathbf p}\, d{\mathbf q}$ which refers to a particular choice of a primitive
   in local Darboux coordinates ${\bf p} = (p_1, \cdots, p_k)$, ${\bf q} = (q_1, \cdots, q_k)$. 
   
   The $\zeta_{ij}$-solitons, being Hamiltonian flow curves for $H=\Imm(\zeta_{ij}W)$ are, according to elementary
   Hamiltonian mechanics, critical points of the ``phase space action" which is the functional
   \[
   h_{\zeta_{ij}}(\phi) \,\,=\,\, \int_\phi  {\mathbf p}\, d{\mathbf q} + H dt \,\,:=\,\,
   \int_{\bf R} \phi^{\ast}\lambda+\Imm(\zeta_{ij} W) dt. 
   \]
   It is defined on the space $\Xc_{ij}$ 
   of smooth maps    $\phi=\phi(t): {\bf R}\to X$ with the boundary conditions 
   \[
   \lim_{t\to -\infty}\phi(t)=x_i, \,\,\, \lim_{t\to +\infty}\phi(t)=x_j. 
   \]
  This interpretation is used in \cite{GMW} for constructing both $L_{\infty}$ and $A_{\infty}$ structures in terms of supersymmetric quantum mechanics.  That is, our transversality assumptions imply that $h_{\zeta_{ij}}$ is itself a Morse function on $\Xc_{ij}$,
  i.e., the $\zeta_{ij}$-solitons are non-degenerate critical points. The grading on $N_{ij}$ is, up to shift, given by
  the relative Morse (Floer) index of $h_{\zeta_{ij}}$. 
  
  The Riemannian metric $g=\Ree(\omega^{1,1})$ on $X$ defines the $L_2$-Riemannian metric on $\Xc_{ij}$:
  \[
  (\delta_1\phi, \delta_2\phi) \,\,=\,\, \int_\RR \bigl( \delta_1\phi(t), \delta_2\phi(t)\bigr)_g dt,
   \]
   and so we can speak about the gradient flow of $h_{\zeta_{ij}}$ on $\Xc_{ij}$. We will view a parametrized curve in $\Xc_{ij}$ 
   with parameter $s$ as a  map $\Phi=\Phi(t,s):\RR^2\to X$ or, equivalently, as an $X$-valued function $\Phi(\tau)$ of one
   complex variable $\tau=t+\sqrt{-1}s \in\CC$.    
   \begin{Defi} (a) The  {\em K\"ahler gradient} of $W$  (with respect to $\omega^{1,1}$) is 
    the vector field  $\on{grad}_{\omega^{1,1}}(W)$  on $X$ (i.e., a smooth section of $TX$)
    obtained from the section $dW$ of $T^*X$ by applying the composite
   isomorphism
   \[
   T^*X \buildrel(\omega^{1,1})^{-1}\over\lra \ol TX \buildrel v\to\ol v \over\lra TX,
   \]
   where $(\omega^{1,1})^{-1}$ is the inverse to \eqref{eq:K-form}, and  $v\to\ol v$ is the canonical $\CC$-antilinear identification
   between $\ol TX$ and $TX$. 
   
   (b) For $\zeta\in\CC^*$, the  {\em $\zeta$-instanton equation}, or the {\em Witten equation} 
   \cite{ruan}, is the following condition on an $X$-valued function $\Phi=\Phi(\tau):  \CC\to X$:
   \[
    \frac{ \partial\Phi }{\partial\ol\tau} \,\,=\,\, -\ol \zeta \on{grad}_{\omega^{1,1}} (W).
   \]
    \end{Defi}
    To explain the meaning of the equation, at any $\tau_0\in\CC$
    \[
  \frac{\partial}{\partial\ol\tau}\Phi (\tau_0)  \,\,=\,\, {{\partial}\over{\partial t}}\Phi(\tau_0)  +\sqrt{-1}{{\partial}\over{\partial s}}\Phi (\tau_0)
   \]
   is a $\CC$-linear combination of two tangent vectors to $X$ at $\Phi(\tau_0)$, so itself an element of $T_{\Phi(\tau_0)}X$.
   It is required that this element is equal to the value of $-\ol \zeta \on{grad}_{\omega^{1,1}}(W)$ at $\Phi(\tau_0)$. 
   
   The   $\zeta$-instanton equation
    can be seen as a complex $\ol\partial$-analog of
   the ordinary differential equation of downward gradient flow for $\Ree(\zeta W)$. In fact,  it reduces to that equation, 
   if we assume that $\Phi(t,s)=\phi(t)$ is independent on $s$. 
   
   \begin{prop}[\cite{GMW}]\label{prop:grad-witt}
   A curve $\Phi(s)$ in $\Xc_{ij}$ represented by a map $\Phi=\Phi(t,s):\RR^2\to\CC$, is a 
   gradient trajectory for $h_{\zeta_{ij}}$, 
   if and only if it satisfies the $\zeta_{ij}$-instanton equation. \qed
   \end{prop}

   \begin{exas} (a) If $X=\CC$ with coordinate $z$ and $\omega^{1,1}=dz\, d\ol z$ is the standard flat K\"ahler metric, then the
   equation has the form
   \[
   \frac{\partial\Phi}{\partial\ol\tau} \,\,=\,\,-\ol \zeta \, \ol{W'(\Phi)}. 
   \]
   (b) More generally, if $(z_1,\cdots, z_k)$ is a local holomorphic coordinate system on $X$, and $\omega^{1,1} = \sum\omega_{a, \ol b}^{1,1}
   dz_a \, d\ol z_b$, 
   then we can write
   $\Phi(\tau)$ as $(\Phi^1(\tau), \cdots,  \Phi^k(\tau))$ and the equation has the form
   \[
     \frac{\partial\Phi^a}{\partial\ol\tau} \,\,=
    \,\, \left({{\partial}\over{\partial t}}+\sqrt{-1}{{\partial}\over{\partial s}}\right)\Phi^a \,\,=
       \,\,-\ol \zeta\,\,  \,\sum_b\,   \eta^{a, \ol b} \frac{\partial \ol W}{\partial \ol \Phi^{\ol b}}.
   \]
   Here $\|\eta^{a, \ol b}\|$ is the inverse matrix to $\|\omega_{a, \ol b}^{1,1}\|$. 
    \end{exas}
    
    \vskip .3cm
    
    \noindent{\bf E. The Maurer-Cartan element.}
  Let $\gen=\gen_{A,N}$ be the $L_\infty$-algebra corresponding to our set $A$ and the coefficient system given by the
  $N_{ij}$ defined in \S B. Thus, see \S   \ref{sec:bim}, 
  \[
  \gen \,\,=\,\,\bigoplus_{Q'} N_{Q'} \otimes\orr(\Sigma(A')), \quad A'=Q'\cap A,
  \]
  where $Q'$ runs over all convex subpolygons with vertices on $A$.
  
  Fix one such $Q'$ and suppose that it has vertices $i_0, \cdots, i_m$
  (counter-clockwise). Then  the basis of $N_{Q'}$ is formed by the 
  vectors $e_{\bphi}$ for all gradient polygons $\bphi = (\phi_{i_0 i_1}, \cdots ,\phi_{i_m i_0})$ that project
  (under $W$) onto the boundary  $\partial Q'$.
  Let $\wen'$ be the web dual to $Q'$, i.e., $\wen'$ is the normal fan of $Q'$, see \S \ref {sec:GMW} and
  Fig. \ref{fig:inst}.
   Thus $\wen'$ has one vertex and the sectors $C_\nu$ of $\wen$ are in bijection with the vertices $i_\nu$ of $Q'$. 
    
    Fix further a gradient polygon $\bphi$ covering $Q'$ and let $\Mc_\zeta(\bphi)$ be the moduli space of solutions
    of the $\zeta$-instanton equations $\Phi: \CC\to X$ with the following asymptotic conditions:
    \begin{enumerate}
    \item[(1)] When $\tau$ approaches infinity deep inside $C_\nu$, then $\Phi(\tau)$ approaches $x_{i_\nu}$.
    
    \item[(2)] Consider an edge $\een$ of $\wen'$ separating  some $C_\nu$ and $C_{\nu+1}$,
    see Fig. \ref{fig:inst}. For a point $\tau\in\een$ let $l(\tau)$ be the line orthogonal to $\een$ passing through $\tau$.
    Then, as $\tau\to\infty$ on $\een$, the restriction $\Phi|_{l(\tau)}$ approaches the soliton $\phi_{i_\nu, i_{\nu+1}}$. 
     \end{enumerate}
     
     \bef
     \btp[scale=0.3, baseline=(current  bounding  box.center)]
     
     \node (0) at (0,0){};
     \fill (0) circle (0.15); 
     \draw[->, line width=0.3mm] (0,0) -- (-7,7); 
      \draw[->, line width=0.3mm] (0,0) -- (5,10); 
       \draw[->, line width=0.3mm] (0,0) -- (9,0); 
        \draw[->, line width=0.3mm] (0,0) -- (3,-6); 
         \draw[->, line width=0.3mm] (0,0) -- (-8,-4); 
         \draw (-8,2) -- (-2,8); 
     
     \node (tau) at (-5,5){};
     \fill (tau) circle (0.15); 
     \node at (-1,5){$C_\nu$}; 
      \node at (-5,1){$C_{\nu+1}$}; 
    \node at (5,4){$C_{\nu-1}$};   
    \node at (-7.5,7.5){\large$\een$}; 
    \node at (-4.2, 5){$\tau$}; 
    \node at (-2, 8.5){$l(\tau)$}; 
    \node at (0,12){\large$\wen'$};
     \etp
     \quad\quad\quad
       \btp[scale=0.4, baseline=(current  bounding  box.center)]
     \node (1) at (0,0){};
     \node (2) at (4,2){};
     \node (3) at (4,5){};
     \node (4) at (0,6){};
     \node (5) at (-3,4){};
     
     \fill (1) circle (0.15); 
        \fill (2) circle (0.15); 
           \fill (3) circle (0.15); 
              \fill (4) circle (0.15); 
                 \fill (5) circle (0.15); 
                 
  \draw (0,0) -- (4,2) -- (4,5) -- (0,6) -- (-3,4) -- (0,0);    
  \node at (1,3) {$Q'$};
  \node at (0,7) {$i_\nu$};
  \node at (4.5, 6){$i_{\nu-1}$}; 
  \node at (-4.5,4){$i_{\nu+1}$};

     \etp
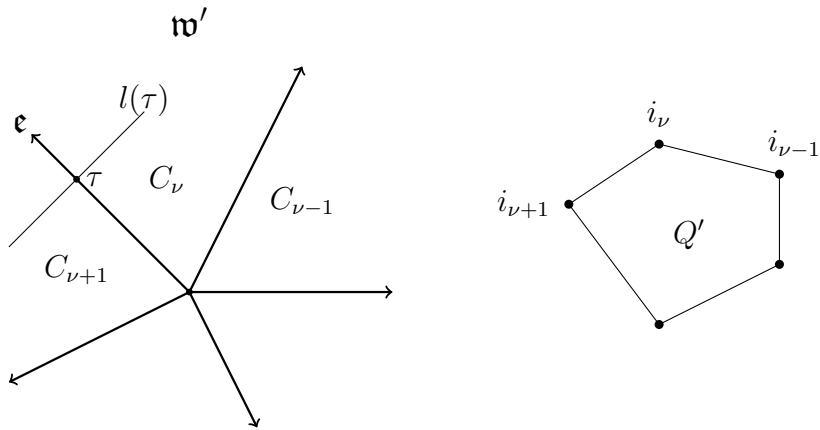
     \caption{The asymptotic condition for $\zeta$-instantons.}
     \label{fig:inst}
     \enf
     
     The proposal of \cite{GMW} can be mathematically
     summarized as follows.
     
     \begin{conj} (a) Under sufficient genericity assumptions, $\Mc_\zeta(\bphi)$ is a manifold of dimension $d(\bphi)-1$,
     equipped with a natural orientation.
     
     (b) In particular, for $\bphi$ such that $d(\bphi)=1$, i.e., $e_{\bphi}$ is a basis vector of $\gen^1$,
     we have a well defined {\em signed cardinality} $\gamma_\zeta(\bphi)\in\ZZ$ of $\Mc_\zeta(\bphi)$ (signs come from the
     orientation of the 0-dimensional manifold). The element
     \[
     \gamma = \gamma_\zeta = \sum_{d(\bphi)=1} \gamma_\zeta(\bphi) e_{\bphi} \,\,\in\,\,\gen^1
     \]
     is a Maurer-Cartan element.
     
     (c) The deformation of the $A_\infty$-algebra $R_\infty$ with respect to $\gamma_\zeta$ is
     identified with $\RS(\zeta W)$, the $A_\infty$-algebra associated to the Fukaya-Seidel category of $\zeta W$. 
     
     \end{conj}
     
     \begin{rems}
     (a) Since $N_{ij}$ is spanned by the critical points of $h_{\zeta_{ij}}$ on $\Xc_{ij}$, it carries
     a natural Floer differential. This differential is in fact a part of the datum provided by $\alpha$,
     corresponding to $m=1$, i.e., to gradient polygons with 2 edges, so that   $m=1$ and $i_0=i, i_1=j$. 
     
     \vskip .2cm
     
     (b) Suppose that  $m\geq 2$ and let $\bphi$ be a gradient polygon lying over the boundary 
     of a convex sub polygon $Q'$. In this case it seems likely that for any $\zeta$-instanton
     $\Phi\in\Mc_\zeta(\bphi)$, the composition
     \[
     \CC \buildrel \Phi\over\lra X \buildrel W\over\lra\CC
     \]
     maps $\CC$ bijectively onto the interior of $Q'$. If this is so, then one can think about $\zeta$-instantons
     as ``sections" of $W$ over $Q'$.

\vskip .2cm

(c) The argument  justifying the Maurer-Cartan condition for $\gamma$ proposed in \cite{GMW}
is based on ``counting the ends" of the 1-dimensional manifolds $\Mc_\zeta(\bphi)$ for $e_\bphi\in\gen^2$. More generally,
this type of argument suggests the existence of a compactification $\ol\Mc_\zeta(\bphi)$ of $\Mc_\zeta(\bphi)$
by adding strata which are  products of $\Mc_\zeta(\bpsi)$ for ``smaller" gradient polygons $\bpsi$.
More precisely, strata are labelled by data $\Pen$ consisting of:
\begin{enumerate}
\item[(1)] A regular polygonal subdivision $\Pc=\{Q''_\nu\}$ of the polygon $Q'$ (image of $\bphi$ under $W$).

\item[(2)] An assignment, to each intermediate edge of the subdivision,  of an instanton projecting onto this edge.
\end{enumerate}
Such a datum can be thought of as a lifting of $\Pc$ into $X$ and
 is equivalent to a ``colored web" in the terminology of \cite{GMW}. Each polygon $Q''_\nu$ is then lifted
 to a gradient polygon $\bpsi_\nu$, and 
 the stratum $\Mc_\zeta(\Pen)$ is the 
 product of the $\Mc_\zeta(\bpsi_\nu)$. In other words, the compactification $\ol\Mc_\zeta(\bphi)$  has  the strata structure
 and the factorization
 property, completely parallel to those of the secondary polytope $\Sigma(A\cap Q')$.

 Further,  for  a sequence of instantons  $(\Phi_n)$ in $\Mc_\zeta(\bphi)$,  the condition of approaching the stratum
 $\Mc_\zeta(\Pen)$, i.e., the condition that $\lim\Phi_n\in \Mc_\zeta(\Pen)$, 
   can be seen as a kind of  tropical
 degeneration. That is, fixing $\varepsilon>0$,
 consider the sets 
  \[
U_n ({\varepsilon}) \, :=\, \{\tau \in \CC: \,\, \| d_\tau\Phi_n\|^2\ge \varepsilon\}.
\] 
Then, as   $\varepsilon\to 0$  and $n\to\infty$ simultaneously  in a compatible way ($\forall \epsilon \, \exists n$),
 the sets $U_n(\varepsilon)$ converge,  in the Gromov-Hausdorff sense, to the web $\wen$ dual to $\Pc$.
The shape of each individual $U_n(\varepsilon)$ can be compared to the amoeba
of a plane algebraic curve. 

 \end{rems}

    \vfill\eject

 \section{Speculations and  directions for further work}\label{sec:spec}
 
 \noindent {\bf 1. From Stasheff polytopes to a categorical structure.}
 A $\zeta_{ij}$-instanton in \S \ref{sec:MCFS} is the $m=2$ case of the following ``catastrophe" (non-generic occurrence)
 for a Morse function: a chain of $m$ critical points of the same index connected by gradient trajectories in a sequential way.
 It was shown in \cite{kapranov-saito} that generic deformations of such a catastrophe are governed by the  
 Stasheff polytope $K_m$ whose vertices correspond to bracketings of $m+1$ factors. As well known, combinatorics
 of  the $K_{m}$ is at the basis of the formalism of $A_\infty$-algebras and categories. 
 It is therefore interesting to relate this appearance
 of  the $K_{m}$ to other categorical structures and to wall crossing formulas in 2 and more dimensions.
 This seems especially appealing since the $K_m$ also describe \cite{kapranov-saito} higher syzygies
 among the Steinberg relations for elementary matrices $e_{ij}(\lambda)\in GL_N$:
 \[
 e_{ij}(\lambda) e_{jk}(\mu) = e_{jk}(\mu) e_{ik}(\lambda+\mu) e_{ij}(\lambda). 
 \]
 Such relations, each having the shape of a pentagon $K_3$, appear naturally in the wall-crossing formulas. 
 
 \vskip .2cm
 
 \noindent {\bf 2.Higher TFT structures and $E_d$-algebras.} It seems certain that in $d$ dimensions the $L_\infty$-algebra
 $\gen$ from \S  \ref{sec:L-inf} can be refined to an $E_d$-algebra, and $\gen_\infty$ to an $E_{d-1}$-algebra.
  We plan to address this in a subsequent paper, currently in preparation.
 The possibility of such refinement is in agreement with the higher analog of the Deligne conjecture
 \cite{defoperads, lurie-fac, ginot} which says that the deformation complex of an $E_{d-1}$-algebra is naturally
 an $E_d$-algebra. The quasi-isomorphism $\Phi$ in the Universality Theorem should then be a morphism of $E_d$-algebras.
 Further, the full higher-dimensional analog of the concept of a coefficient system from \S \ref{sec:bim} should associate to any simplex of
 codimension $p$ with vertices in $A$, an $E_p$-algebra. 
  
   \vskip .2cm
   
   \noindent{\bf 3. Curvilinear theory.} It would be interesting to develop the analog of the formalism of  Fukaya-Seidel categories
   and of \cite{GMW} for the case when the potential $W$ takes values not in $\CC$ but in some Riemann surface $S$. 
   In this case the set $A$ of critical  values makes $S$ into a ``marked surface" in the sense of Teichm\"uller theory and
   one can use isotopy classes of curvilinear triangulations of $S$ with vertices in $A$ in order to analyze
   various categorical structures. In particular, geodesics of quadratic differentials as well as spectral networks 
   of Gaiotto-Moore-Neitzke \cite{GMN}
   seem like natural objects to appear in such a theory. 
   
   \vskip .2cm
   
   \noindent{\bf 4. Infinite-dimensional  case.} Much of the recent ``physical" work on Picard-Lefschetz theory and its
   generalizations was motivated by the infinite-dimensional example of the complexified Chern-Simons functional
   \cite{witten-CS}.
  Here $X$ is the (universal cover of the) space of all connections in a principal bundle on a 3-manifold $M$,
   whose structure group is a complex
  semi simple Lie group. The set $A$ of critical values of the Chern-Simons functional consists of complex numbers 
  of substantial arithmetic importance
  (involving volumes of hyperbolic manifolds, regulators of elements of $K_3(\ol\QQ)$ and the like). 
  The nature of polygons and triangulations that can be
  formed out of these numbers remains mysterious. It would be very interesting if Picard-Lefschetz theory
  imposed some constraints on the convex geometry of these numbers. 
  
  The case of the holomorphic Chern-Simons functional ($M$ is a 3-dimensional complex Calabi-Yau)
  is even more mysterious.

 \vfill\eject

\let\thefootnote\relax\footnote {
M.Ka.: Kavli Institute for Physics and Mathematics of the Universe (WPI), 5-1-5 Kashiwanoha, Kashiwa-shi, Chiba, 277-8583, Japan.
Email: {mikhail.kapranov@ipmu.jp}

\vskip .1cm

M.Ko.: Institute des Hautes \'Etudes Scientifiques, 	35 Route de Chartres, 91440 Bures-sur-Yvette, France. Email:
{maxim@ihes.fr}

\vskip .1cm
Y. So.: Department of Mathematics, Kansas State University,  Cardwell Hall, 
Manhattan, KS 66506 USA. Email: soibel@math.ksu.edu

}

\end{document}